\documentstyle[12pt,epic,eepic,amssymb]{amsart}

\numberwithin{equation}{section}

\newtheorem {Theorem}            	{Theorem} [section]
\newtheorem {Lemma}[equation]     	{Lemma}  	
\newtheorem {Proposition}[equation]	{Proposition}  

\theoremstyle{definition}
\newtheorem{Definition}[equation]{Definition}

\theoremstyle{remark}
\newtheorem{Remark}[equation]{Remark}
\newtheorem{Example}[equation]{Example}
\newtheorem {Corollary} [equation]	{Corollary}

\def	\R	{{\Bbb R}}
\def	\C	{{\Bbb C}}
\def	\Z	{{\Bbb Z}}
\def	\CP	{{\Bbb C}{\Bbb P}}
\def	\Cc	{{\Bbb C ^{\times}}}

\def    \t      {{\frak  t}}
\def    \g      {{\frak  g}}

\def	\ssminus{\smallsetminus}

\def	\to	{\longrightarrow}

\newcommand 	{\dd}[1]	{{\partial \over \partial #1}}

\def	\map	{\text{map}}

\def	\grad	{{\operatorname{grad}}} 

\def	\max	{{\operatorname{max}}} 
\def	\min	{{\operatorname{min}}} 

\def	\span	{{\operatorname{span}}}
\def	\tDH	{Duistermaat-Heckman\ }

\def	\tM	{{\tilde{M}}}
\def	\tC	{{\tilde{C}}}
\def	\tE	{{\tilde{E}}}
\def	\tF	{{\tilde{F}}}

\def	\tomega	{{\tilde{\omega}}}
\def	\tOmega	{{\tilde{\Omega}}}
\def	\tPhi	{{\tilde{\Phi}}}

\def	\ol	{\overline}

\def	\zbar	{\ol{z}}
\def	\wbar	{\ol{w}}
\def	\Bbar	{\ol{B}}
\def	\Cbar	{\ol{C}}

\def	\ss	{\scriptstyle}

\newcommand     {\brak}[1]        {\left< #1  \right>}

\def	\inv	{^{-1}}

\newcommand     {\mf}[1]        {}
\newcommand	{\comment}[1]	{}
\newcommand	{\printname}[1]	{}
\newcommand{\mute}[2] {}

\newcommand{\labell}[1]	{\label{#1}\printname{#1}}
\newcommand{\bibitemm}[2]{\bibitem[#1]{#2}\mf{#2}}

\begin{document}

\title [Periodic Hamiltonian flows on four manifolds]
{Periodic Hamiltonian flows on four dimensional manifolds}

\author{Yael Karshon}

\address{Institute of Mathematics, The Hebrew University
of Jerusalem, Giv'at Ram, Jerusalem 91904, Israel}  


\email{karshon@@math.huji.ac.il}


\date{April 1, 1998} 

\thanks{E-print: dg-ga/9510004.}

\begin{abstract}
We classify the periodic Hamiltonian flows on compact four 
dimensional symplectic manifolds 
up to isomorphism of Hamiltonian $S^1$-spaces. 
Additionally, we show that all
these spaces are K\"{a}hler, that every such space is obtained from a
simple model by a sequence of symplectic blowups, and that if the fixed
points are isolated then the space is a toric variety.
\end{abstract}

\maketitle

\tableofcontents

\section{Introduction}
\labell{sec:intro}
Take a compact symplectic manifold $(M,\omega)$
and a smooth function $\Phi:M \to \R$.
The corresponding Hamiltonian flow is
generated by the vector field $\xi_M$ that satisfies
\begin{equation} \labell{eq:hamilton}
 d\Phi = - \iota(\xi_M) \omega .
\end{equation}
We assume that this flow is {\em periodic\/}, with period $2\pi$,
i.e., that $\xi_M$ generates a circle action on $M$.
The function $\Phi$ is called the {\em moment map\/} for this action.
The triple $(M,\omega,\Phi)$ 
is called a {\em Hamiltonian $S^1$-space}, or a {\em
Hamiltonian circle action}, or a {\em periodic Hamiltonian flow}.
An {\em isomorphism\/} between two such spaces is a
diffeomorphism $F:M_1 \to M_2$ such that $F^*\omega_2 = \omega_1$
and $F^*\Phi_2=\Phi_1$; it follows that $F$ is $S^1$-equivariant.
We always assume that $M$ is connected and that the circle action 
is effective.

Hamiltonian actions of other Lie groups are defined in a similar
way. For a group $G$ we get a moment map $\Phi : M \to \g^*$ 
where $\g$ is the Lie algebra of $G$.
When $G$ is a torus, the {\em completely integrable\/} 
actions, i.e., those whose nonempty reduced spaces, $\Phi^{-1}(\alpha)/G$,
are single points, were classified by Delzant \cite{de}; 
they all turn out to be K\"ahler toric varieties.
(Also see \cite{g:toric}).  
The lowest dimensional Hamiltonian actions that are not completely
integrable are circle actions on 4-manifolds.

M. Audin \cite{au:paper,au:book} and K. Ahara and A. Hattori
\cite{ah-ha} proved that for every compact four dimensional 
Hamiltonian $S^1$-space, the underlying manifold and circle
action are obtained from a {\em minimal model\/} by a sequence of 
{\em equivariant blow ups\/} at fixed points, 
and they listed the minimal models.
This provides a full list of the four dimensional
manifolds and circle actions that admit symplectic forms
and moment maps. This result is strong and beautiful, but it only answers
a small part of the classification problem. For instance,
the above mentioned authors did not determine which 
different blowups produce spaces that are equivariantly diffeomorphic
(their list contains repetitions; see our Examples \ref{ex:sue},
\ref{ex:blow-min}, and \ref{ex:different}), 
nor did they specify which symplectic forms can be put on these spaces. 

The present paper answers these questions. 
We give a complete classification of the 
compact four dimensional Hamiltonian $S^1$-spaces.
The classification consists of a uniqueness part
(sections 2--4) and an existence part (sections 5--7).
The uniqueness part tells us how to determine
whether two spaces are isomorphic.
The existence part lists all the possible spaces.

In sections \ref{sec:graph}--\ref{sec:uniqueness} we associate a labeled
graph to every compact four dimensional Hamiltonian $S^1$-space, and we
show that two such spaces are isomorphic if and only if they have the same
graph. One bi-product is that a Hamiltonian $S^1$-space is determined up
to equivariant symplectomorphism by its underlying manifold, the $S^1$
action, and the cohomology class of the symplectic form; see Proposition
\ref{application}.

In section \ref{sec:toric} we show that if the fixed points are isolated
then the circle action extends to an action of a two dimensional torus
to yield a toric variety. This leads to a classification of the compact
four dimensional Hamiltonian $S^1$-spaces whose fixed points are isolated.

In section \ref{sec:blowup} we prove that a compact four dimensional
Hamiltonian $S^1$ space can be obtained from a minimal model by a sequence
of equivariant {\em symplectic\/} blow-ups. This can also be deduced
from the similar result of Audin, Ahara, and Hattori for the underlying
$S^1$-manifold, but our proof is simpler.

Unlike general symplectic blow-ups, the blow-up of a four dimensional
Hamiltonian $S^1$-space is unambiguous; the resulting Hamiltonian $S^1$-space
is determined up to isomorphism by the fixed point at which we blow up
and the amount by which we blow up. See Proposition \ref{same}.

To complete the classification it remains to determine
by which amounts it is possible to blow up a Hamiltonian $S^1$-space.
Equivalently, we need to determine which invariant symplectic forms 
can be put on our manifolds. This we do in section \ref{sec:existence}.
By the results of section \ref{sec:uniqueness},
it is enough to specify the cohomology classes of the
invariant symplectic forms.
Now, it is easy to state a necessary condition for a cohomology
class to represent an invariant symplectic form, and using
Nakai's criterion we show that this condition is also sufficient.
Moreover, a cohomology class that satisfies this condition
is represented by a compatible {\em K\"ahler\/} form.
Hence every compact four dimensional Hamiltonian $S^1$-space is K\"ahler!

These phenomena do not occur in higher dimensions.
S. Tolman has constructed a compact symplectic
6-manifold with a Hamiltonian action 
of a 2-dimensional torus, with isolated fixed points,
that does not admit a compatible K\"ahler structure. 
In particular, in dimension greater than four, having isolated fixed
points does not imply that the space is toric.
See \cite{tolman-example}; also see \cite{woodward}. 

Unfortunately, the methods in this paper do 
not tell whether two given Hamiltonian $S^1$-spaces are {\em
non-equivariantly\/} symplectomorphic. In particular,
given a symplectic 4-manifold, we do not give the list 
of all Hamiltonian $S^1$-actions on it.

The results of this paper for the case of isolated
fixed points already appeared in the author's thesis \cite{k:thesis}.
A five page summary of this paper appeared in \cite{k:newton}.

S. Tolman and the author are currently writing up a treatment
of Hamiltonian torus actions on higher dimensional manifolds 
where the dimension of the torus is one less than half
the dimension of the manifold \cite{k-t:deficiency}.

{\bf Acknowledgements.} \ 
I wish to thank the following people for useful discussions:
D.~Abramovich, M.~Audin, A.~Canas Da Silva,
V.~L.~Ginzburg, M.~Grossberg, V.~Guillemin, J.~Harris,
D.~McDuff,  E.~Lerman, R.~Sjamaar, S.~Sternberg,  S.~Tolman,
C.~Woodward.
The foundations of the subject lie in early papers by Guillemin
and Sternberg and by Atiyah \cite{at,g-s:convexity}.
I was influenced by the papers of M.\ Audin \cite{au:paper,au:book},
K.\ Ahara and A.\ Hattori \cite{ah-ha}, and T.\ Delzant \cite{de}.
M. Audin, E. Lerman, S. Singer, and S. Tolman, made helpful remarks 
on various versions of the manuscript.
The referee made excellent suggestions for improving the exposition,
which lead to a complete revision of the paper, for the better.

I proceed with more specific acknowledgements:

Many ideas come from \cite{au:paper,au:book,ah-ha}:
the chains of gradient spheres appear there,
Lemma \ref{exist-delzant'} is very similar to \cite[\S3.3]{au:paper} 
(the proof is different),
Proposition \ref{lem:two-branches} was claimed in 
\cite[Lemma 3.1.2]{au:paper} and was proved in \cite[\S6]{ah-ha}
(our proof is different), subsection \ref{subsec:down-to-minimal}
contain a new proof of results of \cite{au:paper,ah-ha}.
Section \ref{sec:uniqueness} is strongly influenced by the first 
section of Delzant's paper, \cite{de}.
D. McDuff helped me understand the work of Audin, Ahara and Hattori, 
she pointed out the problem with \cite[Prop.3.1.2]{au:paper},
and she indicated that one only needs to blow down and not to blow up, 
which made life and section \ref{sec:blowup} simpler.  
S. Tolman convinced me to study Hamiltonian circle actions
in dimension four before attempting to understand them in higher dimension;
if not for her advice I would still be staring in space. 
She clarified the important Example \ref{ex:sue},
and pointed at gaps in an earlier proof of Proposition \ref{semi-morphism}.
I learned of Nakai's criterion from D. McDuff and L. Polterovich,
and I could not have applied it in section \ref{sec:existence} without 
consulting D.\ Abramovich and J.\ Harris.
I learned the trick in the proof of Lemma \ref{kho} from a 
talk of A.\ Khovanskii, and I learned Lemma \ref{mg} from M.\ Grossberg.
Proposition \ref{prop:conj-del} was conjectured by T.\ Delzant and
V.\ Ginzburg, from whom I heard this conjecture.  

I wish to thank the Weizmann Institute of science for their hospitality
and support during the summer of 1992.
I was supported by an Alfred P.\ Sloan Dissertation Fellowship
in the academic year 1992-93.
I was supported by NSF grant DMS-9404404 during part of the work
on this paper.

{\bf Related works:}
Most important is
the work of Mich\`ele Audin \cite{au:paper,au:book} and of K.\ Ahara and 
A.\ Hattori \cite{ah-ha} that was mentioned earlier.
T.\ Delzant \cite{de} classified the completely integrable Hamiltonian 
torus actions on compact symplectic manifolds.
He also showed that a Hamiltonian circle action whose moment map
has exactly two singular values, one of which is non-degenerate,
is isomorphic to a standard action on $\CP^n$.
For Hamiltonian circle actions on compact symplectic four manifolds, 
already in 1989 E. Lerman wrote explicit formulas
describing the pre-image via the moment map 
of an interval of regular values \cite{le:notes}.

In other than the symplectic category,
locally smooth circle actions on four manifolds (no symplectic structure) 
were classified by Fintushel \cite{fin},
and holomorphic circle actions on complex projective surfaces 
were classified by Orlik and Wagreich \cite{or-wag},
who produced a list of surfaces and actions identical to Audin's list
in \cite{au:paper}.\footnote
{
Clearly, every complex projective surface is a symplectic four manifold.
However, non-isomorphic complex projective surfaces could be isomorphic
as Hamiltonian $S^1$-spaces; see Example \ref{ex:blow-min} of the
present paper.
In the other direction, every compact four dimensional Hamiltonian 
$S^1$ space is equivariantly symplectomorphic to a complex projective 
surface with a holomorphic circle action \cite{ah-ha,au:paper}.
Moreover, every compact four dimensional Hamiltonian $S^1$ space is 
K\"ahler; see section \ref{sec:existence} of the present paper. 
Of course, if the two-form is not integral, the space is not symplectomorphic
to a complex projective surface.
}

Hamiltonian actions of non-abelian groups are more difficult
than torus actions.
Completely integrable actions of non-abelian groups have been 
studied in recent years by Delzant, Guillemin, Knop, Sjamaar, de Souza, 
Woodward, and possibly others.

{\bf A technical remark: }
We use the following conventions. 
The circle group consists of the complex numbers of norm $1$.
Its Lie algebra is identified with $\R$ such that 
the exponential map is $t \mapsto e^{it}$ and its kernel is 
$\l=2\pi \Z$.
The dual of the Lie algebra is $\t^*=\R$,
and the weight lattice is $\l^*=\Z$,
so that $\brak{\l^*,\l} = 2\pi \Z$. 
Lebesgue measure on $\t^*$ is the standard 
measure on $\R$ so that the volume of $\t^*/\l^*$ is $1$.
The symplectic form on $\C$ is $r dr \wedge d\theta$
(in polar coordinates),
and the moment map for the standard circle action
is $\Phi = r^2/2$.
A disc around the origin of area $2\pi \lambda$ 
gets mapped to an interval of
length $\lambda$ (area$=\pi r^2$, $\lambda=r^2/2$).

\section{Graphs}
\labell{sec:graph}
In section \ref{subsec:graph} we associate a labeled graph to each compact
four dimensional Hamiltonian $S^1$ space.  In section \ref{subsec:delzant}
we describe the most important examples: Delzant spaces.
In section \ref{subsec:DH} we discuss the \tDH measure in relation to the
graph.

\subsection{The graph}
\labell{subsec:graph}
Let $(M,\omega)$ be a compact symplectic four-manifold 
with a Hamiltonian circle action and a moment map $\Phi:M \to \R$.
In Appendix \ref{sec:normal} we recall the following facts:

\begin{Lemma} \labell{fixed}
Each component of the fixed point set is either a single point
or a symplectic surface. 
The maximum and minimum of the moment map is each attained
on exactly one component of the fixed point set.
Fixed points on which the moment map is not extremal
are isolated. 
\end{Lemma}

We call a fixed point {\em extremal\/} ({\em maximal\/} or {\em minimal\/})
if it is an extremum for the moment map; otherwise, we call it an 
{\em interior\/} fixed point.

\begin{Lemma} \labell{Zk}
For each integer $k \geq 2$, consider the set of points
whose stabilizer is equal to the cyclic subgroup of $S^1$ of order $k$,
$$ Z_k = \{ \lambda \in S^1 \ | \ \lambda^k =1 \}.$$
Each connected component of the closure of this set is a closed 
symplectic two-sphere, on which the quotient circle, $S^1 / Z_k$, 
acts with two fixed points.  
\end{Lemma}

We call such a sphere a {\em $Z_k$-sphere}.

\begin{figure}
$$
\begin{picture}(107,99)(0,-10)
\put(30,40){\ellipse{60}{60}}
\put(30,70){\blacken\ellipse{6}{6}}
\put(30,10){\blacken\ellipse{6}{6}}
\spline(20,35) (45,35)(75,40) (75,45)(65,50)
\path(73.050,48.211)(65.000,50.000)(71.261,44.633)
\put(35,0){$\ss q$}
\put(35,75){$\ss p$}
\put(70,30){$\ss k \mbox{ \scriptsize times} $}
\end{picture}
$$
\caption{A $Z_k$-sphere}
\end{figure}
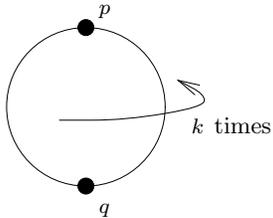

We now construct the graph associated to $(M,\omega,\Phi)$:\footnote
{
Earlier versions of this manuscript contained an equivalent,
but slightly more complicated, construction.
}

{\em
To every component of the fixed point set we assign a vertex,
and to every $Z_k$-sphere we assign an edge connecting the
corresponding vertices.
We label each edge by the isotropy weight, $k$, of the corresponding
$Z_k$-sphere. We label each vertex by the value of the moment map
on the corresponding fixed point set.
Additionally, to a vertex that corresponds to a symplectic surface, $B$,
we attach two additional labels: the genus of that surface, and its 
normalized symplectic area, defined to be $\frac{1}{2\pi} \int_B \omega$.
}

We call the labels of the vertices ``moment map labels",
``area labels", and ``genus labels", respectively.

\begin{Remark}
In our figures, we will often omit some of the labels.
The moment map labels will be indicated by the height
of a vertex in the plane. 
Those vertices that correspond to fixed surfaces will be drawn fatter
than those that correspond to isolated fixed points.
\end{Remark}

\begin{Example} \labell{ex:s2s2}
Let $M$ be the product of two spheres of radius $1$, each with the 
standard area form. Let the circle act by rotating the second sphere 
at twice the speed of the first:
$\lambda \cdot (\vec{u},\vec{v}) = (\lambda \vec{u}, \lambda^2 \vec{v})$,
where $(\vec{u},\vec{v}) \in S^2 \times S^2 \subset \R^3 \times \R^3$,
and where the action on $S^2$ is by rotations in the first two coordinates
of $\R^3$.
There are four fixed points: $(n,n)$, $(s,n)$, $(n,s)$, and $(s,s)$,
where $n$ and $s$ are the north and south poles of $S^2$.
There are two $Z_2$-spheres: $\{ n \} \times S^2$ and $\{ s \} \times S^2$.
The moment map is $\Phi(\vec{u},\vec{v}) = u_3 + 2v_3$.
The graph is shown in Figure \ref{graphs2s2}.
\end{Example}

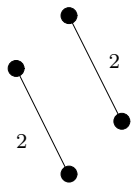
\begin{figure}
$$
\begin{picture}(46,81)(0,-10)
\put(23,63){\blacken\ellipse{6}{6}}
\put(3,43){\blacken\ellipse{6}{6}}
\put(43,23){\blacken\ellipse{6}{6}}
\put(23,3){\blacken\ellipse{6}{6}}
\path(3,43)(23,3)
\path(43,23)(23,63)
\put(3,13){$\ss 2$}
\put(38,43){$\ss 2$}
\end{picture}
$$
\caption{Graph for Examples \protect\ref{ex:s2s2} and \ref{ex:what}}
\labell{graphs2s2}
\end{figure}

Lemmas \ref{fixed} and \ref{Zk} force the graph to have a simple shape:
\begin{itemize}
\item 
there is a unique top vertex and a unique bottom vertex;
\item
the edges occur in a finite number of branches,
with the moment map labels increasing along each branch;
a branch needn't reach an extremal vertex;
\item
an extremal vertex is reached by at most two edges;
an extremal ``fat" vertex is not reached by any edge.
\end{itemize}

The isotropy weights at the fixed points can be read from the graph:
\begin{itemize}
\item
for $k \geq 2$, a fixed point has an isotropy weight $-k$ 
if and only if it is the north pole of a $Z_k$-sphere, 
and it has an isotropy weight $k$ if and only if it is the south pole
of a $Z_k$-sphere;
\item
an interior fixed point has one positive weight and one negative
weight, a maximal fixed point has both weights non-positive, a minimal
fixed point has both weights non-negative; 
\item 
a fixed point has a weight $0$ if and only if it lies on a fixed surface.
\end{itemize}

For example, in Figure \ref{graphs2s2}, the isotropy weights corresponding
to the left interior vertex are $\{-2,1\}$, and those corresponding to
the top vertex are $\{-2,-1\}$.

Note that the graph and its integer labels depend only 
on the manifold and circle action. The real labels 
are essentially determined by the cohomology class of the
symplectic form:

\begin{Lemma} \labell{cohomology}
The cohomology class of $\omega$ 
determines the moment map values at the fixed points
up to a simultaneous shift of all these values by the same amount,
and it determines the normalized symplectic areas of fixed surfaces.
\end{Lemma}

\begin{pf}
The second part is clear. For the first part,
let $p$ and $q$ be fixed points with $\Phi(p) > \Phi(q)$.
Choose any smooth path, $\gamma(t)$, $0 \leq t \leq 1$, 
from $p$ to $q$. Denote $\square := [0,2\pi] \times [0,1]$,
and define $f : \square \to M$
by $f(s,t) = e^{is} \cdot \gamma(t)$.
Then 
$\int_\square f^* \omega = 2\pi \int_0^1 \gamma^*( \iota(\xi_M) \omega )
 = 2\pi ( \Phi(p) - \Phi(q) )$.
Since $f$ defines a {\em cycle\/} in homology,
this integral depends only on the cohomology class of $\omega$, 
hence so does the difference $\Phi(p) - \Phi(q)$.
\end{pf}

In particular, if $p$ and $q$ are the north and south poles of a
$Z_k$-sphere, the difference $\Phi(p) - \Phi(q)$ is equal to the
symplectic area of the sphere times $k / 2\pi$.

\mute{reduced spaces read from graph}{
Recall that the {\em reduced spaces\/} are the quotients 
$\Phi\inv(\alpha) / S^1$ with $\alpha \in \R$.
Topologically, these are all surfaces 
(unless $\Phi\inv(\alpha)$ is an isolated extremum or is empty).
As observed by McDuff, these surfaces are all homeomorphic to each other;
a homeomorphism can be constructed out of the gradient flow
for the moment map with respect to any invariant metric
\cite{md:four-manifolds,au:book}.
Hence, topologically, the reduced surfaces are determined
by their genus, which is equal to the genus of any fixed surface, 
if one exists, and is zero if either the minimum or the maximum 
of the moment map is isolated.
If $\alpha$ is a regular value of the moment map then, 
differentiably, the reduced space is an orbifold; 
a singular point of order $k$ arises whenever the level set
$\Phi\inv(\alpha)$ meets a $Z_k$-sphere.
All this information, which determines the reduced spaces,
can be read from the graph.
}

\subsection{K\"ahler toric varieties}
\labell{subsec:delzant}
Take a compact symplectic manifold $(M,\omega)$ of dimension $2n$
with a Hamiltonian action of an $n$-torus, 
$T^n = S^1 \times \ldots \times S^1$,
and a moment map $\Phi: M \to \R^n$, meaning a map whose $n$
coordinates generate, via (\ref{eq:hamilton}), 
the actions of the $n$ circles.
Such a triple $(M,\omega,\Phi)$ is called a {\em Delzant space}.

By the convexity theorem \cite{g-s:convexity,at}, the image of the moment
map is a convex polytope. By Delzant's theorem \cite{de}, this
polytope determines the Hamiltonian space up to equivariant 
symplectomorphism, and the space is a {\em K\"ahler toric variety},
meaning that it admits a complex structure such that the torus $T$
acts holomorphically and the symplectic form $\omega$ is K\"ahler.
See \cite{de} and \cite{g:toric}.
The polytopes that arise in this way are called
{\em Delzant polytopes}. When $n=2$, these are exactly those polygons
in $\R^2$ that have the following properties:
\begin{equation} \labell{delzant-condition}
\parbox{4.0in}{
\begin{itemize}
\item
the slopes of the edges are rational or infinite;
\item
every two consecutive edges have integral outward normal vectors 
$(k,b)$ and $(k',b')$ with $kb' - bk' = 1$.
\end{itemize}
}
\end{equation}

The pre-image in the manifold of a vertex of the polygon 
is a fixed point for the torus action.
The pre-image of an edge with slope $k/b$ is an invariant
2-sphere whose stabilizer is the subgroup of $S^1 \times S^1$
consisting of the elements of the form $(\lambda^k,\lambda^{-b})$.
The pre-image of the interior of the polygon consists of free torus orbits.
These facts follow from the local normal form for Hamiltonian torus actions 
and from the connectedness of the level sets of the moment map,
and are explained in Delzant's paper \cite{de}.

If we now restrict the action to the sub-circle $ \{ e \} \times S^1$, 
we get a compact four dimensional Hamiltonian $S^1$-space. 
The moment map for the $S^1$-action is the $T$-moment map 
composed with the projection $\R^2 \to \R$ to the second coordinate.
The fixed surfaces are the pre-images, under the $T$-moment map,
of the horizontal edges of the Delzant polygon.
Such a surface has genus zero, and its normalized symplectic area
is equal to the length of the corresponding horizontal edge.
The isolated fixed points are the pre-images of those vertices of the 
polygon that do not lie on horizontal edges. 
The $Z_k$-spheres are the pre-images of edges with slope $\pm k/b$ 
in reduced form, where $b$ is relatively prime to $k$.  
With this information, it is easy to construct the graph for the $S^1$
space out of the Delzant polygon. 

\begin{Remark}
In our figures of Delzant polygons, the dots mark the weight lattice $\Z^2$
in $\R^2$.
\end{Remark}

\begin{Example}\labell{ex:what}
The Delzant spaces whose polygons are drawn in figure \ref{polys2s2} 
give Hamiltonian $S^1$-spaces with the same graph, which is 
the graph drawn in Figure \ref{graphs2s2}. 
These $S^1$-spaces are all isomorphic to the one described 
in Example \ref{ex:s2s2}.
The three polygons all correspond to $S^2 \times S^2$, with its 
standard torus action, composed with an automorphism of the torus 
that preserves the second sub-circle. 
\end{Example}

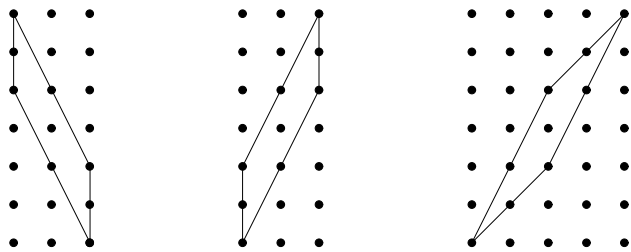
\begin{figure}
$$
\setlength{\unitlength}{0.0100in}
\begin{picture}(324,139)(0,-10)
\put(2,122){\blacken\ellipse{4}{4}}
\put(22,122){\blacken\ellipse{4}{4}}
\put(42,122){\blacken\ellipse{4}{4}}
\put(42,102){\blacken\ellipse{4}{4}}
\put(22,102){\blacken\ellipse{4}{4}}
\put(2,102){\blacken\ellipse{4}{4}}
\put(2,82){\blacken\ellipse{4}{4}}
\put(22,82){\blacken\ellipse{4}{4}}
\put(42,82){\blacken\ellipse{4}{4}}
\put(22,62){\blacken\ellipse{4}{4}}
\put(2,62){\blacken\ellipse{4}{4}}
\put(42,42){\blacken\ellipse{4}{4}}
\put(42,22){\blacken\ellipse{4}{4}}
\put(22,22){\blacken\ellipse{4}{4}}
\put(2,22){\blacken\ellipse{4}{4}}
\put(2,2){\blacken\ellipse{4}{4}}
\put(22,2){\blacken\ellipse{4}{4}}
\put(42,2){\blacken\ellipse{4}{4}}
\put(42,2){\ellipse{4}{4}}
\put(2,42){\blacken\ellipse{4}{4}}
\put(22,42){\blacken\ellipse{4}{4}}
\put(42,62){\blacken\ellipse{4}{4}}
\put(122,122){\blacken\ellipse{4}{4}}
\put(142,122){\blacken\ellipse{4}{4}}
\put(162,122){\blacken\ellipse{4}{4}}
\put(162,102){\blacken\ellipse{4}{4}}
\put(142,102){\blacken\ellipse{4}{4}}
\put(122,102){\blacken\ellipse{4}{4}}
\put(122,82){\blacken\ellipse{4}{4}}
\put(142,82){\blacken\ellipse{4}{4}}
\put(162,82){\blacken\ellipse{4}{4}}
\put(142,62){\blacken\ellipse{4}{4}}
\put(122,62){\blacken\ellipse{4}{4}}
\put(162,42){\blacken\ellipse{4}{4}}
\put(162,22){\blacken\ellipse{4}{4}}
\put(142,22){\blacken\ellipse{4}{4}}
\put(122,22){\blacken\ellipse{4}{4}}
\put(122,2){\blacken\ellipse{4}{4}}
\put(142,2){\blacken\ellipse{4}{4}}
\put(162,2){\blacken\ellipse{4}{4}}
\put(122,42){\blacken\ellipse{4}{4}}
\put(142,42){\blacken\ellipse{4}{4}}
\put(162,62){\blacken\ellipse{4}{4}}
\put(242,122){\blacken\ellipse{4}{4}}
\put(262,122){\blacken\ellipse{4}{4}}
\put(262,102){\blacken\ellipse{4}{4}}
\put(242,102){\blacken\ellipse{4}{4}}
\put(242,82){\blacken\ellipse{4}{4}}
\put(262,82){\blacken\ellipse{4}{4}}
\put(242,62){\blacken\ellipse{4}{4}}
\put(262,42){\blacken\ellipse{4}{4}}
\put(262,22){\blacken\ellipse{4}{4}}
\put(242,22){\blacken\ellipse{4}{4}}
\put(242,2){\blacken\ellipse{4}{4}}
\put(262,2){\blacken\ellipse{4}{4}}
\put(242,42){\blacken\ellipse{4}{4}}
\put(262,62){\blacken\ellipse{4}{4}}
\put(282,122){\blacken\ellipse{4}{4}}
\put(302,122){\blacken\ellipse{4}{4}}
\put(322,122){\blacken\ellipse{4}{4}}
\put(322,102){\blacken\ellipse{4}{4}}
\put(302,102){\blacken\ellipse{4}{4}}
\put(282,102){\blacken\ellipse{4}{4}}
\put(282,82){\blacken\ellipse{4}{4}}
\put(302,82){\blacken\ellipse{4}{4}}
\put(322,82){\blacken\ellipse{4}{4}}
\put(302,62){\blacken\ellipse{4}{4}}
\put(282,62){\blacken\ellipse{4}{4}}
\put(322,42){\blacken\ellipse{4}{4}}
\put(322,22){\blacken\ellipse{4}{4}}
\put(302,22){\blacken\ellipse{4}{4}}
\put(282,22){\blacken\ellipse{4}{4}}
\put(282,2){\blacken\ellipse{4}{4}}
\put(302,2){\blacken\ellipse{4}{4}}
\put(322,2){\blacken\ellipse{4}{4}}
\put(282,42){\blacken\ellipse{4}{4}}
\put(302,42){\blacken\ellipse{4}{4}}
\put(322,62){\blacken\ellipse{4}{4}}
\path(2,122)(2,82)(42,2) (42,42)(2,122)
\path(162,122)(162,82)(122,2) (122,42)(162,122)
\path(322,122)(282,82)(242,2) (282,42)(322,122)
\end{picture}
$$
\caption{Delzant polygons for Example \protect\ref{ex:what}}
\labell{polys2s2}
\end{figure}

\begin{Example} \labell{ex:polygons}
The five polygons in Figure \ref{polygons} give rise to the respective
graphs in Figure \ref{graphs}; we invite the reader to check this.
The first two of these are discussed further in Example \ref{ex:sue}.
\end{Example}

\begin{figure}

\newsavebox{\hlattice}
\savebox{\hlattice}{
\begin{picture} (0,0)
\multiput(0,0)(15,0){3}{\circle*{2}}
\end{picture}
}
\newsavebox{\slattice}
\savebox{\slattice}{
\begin{picture} (0,0)
\multiput(0,0)(15,0){2}{\circle*{2}}
\end{picture}
}
\begin{center}
\begin{picture}(15,60)(0,0)
\multiput(-4,0)(0,15){5}{\usebox{\slattice}}
\path(0,0)(0,45)(15,60)(15,15)(0,0)
\end{picture}
\hfill
\begin{picture}(15,60)(0,0)
\multiput(-4,0)(0,15){5}{\usebox{\slattice}}
\path(0,15)(0,45)(15,60)(15,0)(0,15)
\end{picture}
\hfill
\begin{picture}(30,60)(0,0)
\multiput(-4,0)(0,15){5}{\usebox{\hlattice}}
\path(0,0)(0,15)(15,45)(30,60)(30,45)(15,15)(0,0)
\end{picture}
\hfill
\begin{picture}(30,60)(0,0)
\multiput(-4,0)(0,15){5}{\usebox{\hlattice}}
\path(0,30)(0,60)(30,60)(30,0)(0,30)
\end{picture}
\hfill
\begin{picture}(30,60)(0,0)
\multiput(-4,0)(0,15){5}{\usebox{\hlattice}}
\path(0,0)(30,60)(30,30)(0,0)
\end{picture}
\end{center}
\caption{Polygons for Examples \protect\ref{ex:polygons} and
\protect\ref{ex:grad}}
\labell{polygons}

\begin{center}
\begin{picture}(20,70)(0,0)
\put(10,0){\circle*{4}}
\put(10,15){\circle*{4}}
\put(10,45){\circle*{4}}
\put(10,60){\circle*{4}}
\end{picture}
\hfill
\begin{picture}(20,70)(0,0)
\put(10,0){\circle*{4}}
\put(10,15){\circle*{4}}
\put(10,45){\circle*{4}}
\put(10,60){\circle*{4}}
\end{picture}
\hfill
\begin{picture}(20,70)(0,0)
\put(10,0){\circle*{4}}
\put(0,15){\circle*{4}}
\put(20,15){\circle*{4}}
\put(0,45){\circle*{4}}
\put(20,45){\circle*{4}}
\put(10,60){\circle*{4}}
\path(0,15)(0,45)
\put(-6,30){$\ss 2$}
\path(20,15)(20,45)
\put(14,30){$\ss 2$}
\end{picture}
\hfill
\begin{picture}(20,70)(0,0)
\put(10,60){\blacken\ellipse{30}{10}}
\put(10,30){\circle*{4}}
\put(10,0){\circle*{4}}
\end{picture}
\hfill
\begin{picture}(20,70)(0,0)
\put(20,0){\circle*{4}}
\put(0,30){\circle*{4}}
\put(20,60){\circle*{4}}
\path(20,0)(20,60)
\put(22,30){$\ss 2$}
\end{picture}
\end{center}
\caption{Graphs for Example \protect\ref{ex:polygons}}
\labell{graphs}

\end{figure}

\subsection{Push-forward measures}
\labell{subsec:DH}
{\em Liouville measure\/} on a $2n$ dimensional symplectic manifold 
$(M,\omega)$ is defined by integration of the volume form, 
$\omega^n / n!$, with respect to the symplectic orientation,
divided by the irrelevant factor of $(2\pi)^n$.
Take an effective action of a torus $G$ on $M$ with a moment map 
$\Phi: M \to \g^*$.
The push-forward of Liouville measure via the moment map $\Phi$
is a measure on $\g^*$ called the {\em \tDH\/} measure. 
If $M$ is compact, this measure is absolutely continuous 
with respect to Lebesgue measure, and the density function 
is piecewise polynomial \cite{d-h}.

\begin{Remark}
\labell{ex:guillemin}
Consider a K\"ahler toric variety of real dimension 4.
Its \tDH measure for the torus action
is equal to Lebesgue measure on the corresponding
polygon.
If we restrict the action to the circle subgroup
$\{e\} \times S^1$, the corresponding \tDH measure
is equal to the push-forward to $\R$ of Lebesgue measure
on the polygon via the projection $(x,y) \mapsto y$.

Consider, for instance, the Hamiltonian $S^1$-spaces 
that correspond to the first
three polygons in Figure \ref{polygons}.
The third space is not isomorphic to the first or the second --
its graph, in Figure \ref{graphs}, is different. 
However, all three spaces have the same \tDH
measure on $\R$; its density function (divided by $(2\pi)^2$) 
is illustrated in Figure \ref{density}.
Hence a compact four dimensional Hamiltonian $S^1$-space
is not determined by its push-forward measure.
This answers a question of V. Guillemin.
\end{Remark}

\begin{figure}
$$
\setlength{\unitlength}{0.0100in}
\begin{picture}(140,85)(0,-10)
\path(0,20)(140,20)
\path(30,20)(50,40)(90,40)(110,20)
\path(20,70)(20,0)
\path(15,40)(20,40)
\path(15,60)(20,60)
\put(5,35){$\ss 1$}
\put(5,55){$\ss 2$}
\end{picture}
$$
\caption{The \tDH measure for the $S^1$-spaces corresponding 
to the first three polygons in Figure \protect\ref{polygons}.
}
\labell{density}
\end{figure}
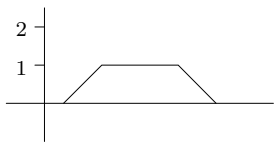

Guillemin, Lerman, and Sternberg \cite[\S3]{g-l-s:kostant}
wrote an explicit formula for the density function
of the \tDH measure in terms of fixed point data,
when the Lie group is a torus and the fixed points are isolated.
Canas da Silva and Guillemin extended this formula
to allow non-isolated fixed points and orbifolds \cite{C}.
We will use this formula for circle actions on four manifolds,
possibly with non-isolated fixed points.
Moreover, in section \ref{sec:existence}
we will need to apply this formula to closed invariant two-forms
that are not necessarily symplectic, i.e., that might be degenerate.
Let us clarify the relevant definitions.

Suppose that $M$ is a $2n$-dimensional compact oriented manifold
with a circle action, equipped with a closed invariant two-form $\omega$,
which is not necessarily symplectic.  The Liouville measure, which is
still defined by integration of $\omega^n/n!$ divided by $(2\pi)^n$,
is now a {\em signed
measure\/} on the manifold; its density function is negative wherever
$\omega^n$ is inconsisent with the given orientation.
Note that, whereas in the symplectic case we integrated with respect
to the symplectic orientation, here the orientation must be given
as an extra piece of structure.
Now suppose that a Lie group $G$ acts on $M$ and preserves the closed two-form
$\omega$.  For simplicity, let us assume that $G$ is a circle.
A moment map $\Phi$ can still be defined by Equation \eqref{eq:hamilton}.
As in the symplectic case, if a moment map exists then it is unique
up to translation, and the obstruction to the existence of a moment map
lies in the first de Rham cohomology of the manifold.
The cohomology class of $\omega$ still determines
the moment map values at the fixed points, up to a global translation;
the proof is identical to that of Lemma \ref{cohomology}.

Let us now fix a compact four manifold $S^1$-manifold, $M$,
with a symplectic form and a moment map; call its moment map
the {\em symplectic moment map}, to distinguish it from 
a moment map for another closed two-form, which will appear later.
Denote by $B_\min$ and $B_\max$ the extremal sets of the symplectic
moment map.
If $B_\min$ is a two dimensional surface, denote its self intersection
by $e_\min$.
If $B_\min$ is an isolated fixed point with isotropy weights $n$ and $n'$,
define $e_\min = -1/nn'$. 
Similarly, denote by $e_\max$ either the self intersection of $B_\max$,
if $B_\max$ is a surface, or the number $-1/mm'$ if $B_\max$ 
is an isolated fixed point with isotropy weights $-m$ and $-m'$.
Finally,
denote the isotropy weights at an interior fixed point $p$ 
by $-m_p$ and $n_p$.

Now let $\omega$ by {\em any\/} closed invariant two-form on $M$, 
not necessarily symplectic, and let $\Phi$ be a corresponding moment map.
Denote  $a_\min = \frac{1}{2\pi} \int_{B_\min} \omega$
and  $a_\max = \frac{1}{2\pi} \int_{B_\max} \omega$,
and denote the values of $\Phi$ at the fixed points by 
$y_\min = \Phi(B_\min)$, $y_\max = \Phi(B_\max)$, and $y_p = \Phi(p)$
for interior fixed points $p$.
Notice that the inequalities $y_\min < y_p < y_\max$ 
might no longer hold.
Denote
\begin{equation} \labell{eq:Theta}
H(x) = 	\begin{cases}	1 \quad \mbox{if $x \geq 0$} \\
	   	   	0 \quad \mbox{if $x < 0$}
 	\end{cases} 
\quad \mbox{and} \quad
\Theta(x) = \begin{cases}  x  \quad \mbox{if $x \geq 0$} \\
			   0  \quad \mbox{if $x < 0$}
	    \end{cases}.
\end{equation}

We now state the special case of the Guillemin-Lerman-Sternberg formula
that we need:

\begin{Lemma}
The density function for the \tDH measure is 
\begin{equation} \labell{eq:rho}
\begin{array}{rcl}
\rho(y) & = & a_\min H(y-y_\min) - e_\min \Theta(y-y_\min) \\
 &  & 		-  \sum_p \frac{1}{m_p n_p} \Theta(y-y_p)  \\
 &  & 	      	- e_\max \Theta(y - y_\max) 
	      	- a_\max H(y-y_\max).
\end{array}
\end{equation}
\end{Lemma}

\begin{pf*}{Outline of the proof}
%
The Fourier transform of the \tDH measure is the function 
\begin{equation} \labell{FT}
	t \mapsto \frac{1}{(2\pi)^2} \int_M e^{-it\Phi} \omega^2 / 2!  .	
\end{equation}
Applying the Atiyah-Bott-Berline-Vergne localization formula 
for equivariant differential forms, we can express the function \eqref{FT}
as a sum of contributions from the fixed points (see \cite{a-b:equivariant}). 
Guillemin, Lerman, and Sternberg showed how to perform the
inverse Fourier transform to this sum; this yields the formula \eqref{eq:rho}.  
\end{pf*}

\begin{Lemma} \labell{cor:DH}
Let $(M,\omega)$ be a symplectic four manifold with a circle
action and a moment map.   Then
\begin{equation} \labell{e+e}
   -e_\min - e_\max = \sum_p {1 \over m_p n_p},
\end{equation}
where we sum over the interior fixed points, and where
$e_\min$, $e_\max$, $m_p$, and $n_p$, are as defined above. 
In particular, if all the fixed points are isolated,
and if $-m$ and $-m'$ are the isotropy weights at the maximum
and $n$ and $n'$ are the isotropy weights at the minimum, then
\begin{equation} \labell{DHsum=0}
 { 1 \over mm'} + {1 \over nn'} = \sum_p {1 \over m_p n_p}.
\end{equation}
\end{Lemma} 

\begin{pf}
For $\omega \equiv 0$, $\Phi \equiv 0$, and $y > 0$,
\eqref{eq:rho} reads $0 = (-e_\min - \sum_p \frac{1}{m_p n_p} - e_\max) y$,
implying \eqref{e+e}.
\end{pf}

\begin{Lemma} \labell{lem:Euler}
Let $M$ be a compact symplectic four dimensional Hamiltonian
$S^1$-space.  Then the numbers $e_\min$, $e_\max$ 
and the density function $\rho(y)$ are determined by the associated graph.
\end{Lemma}

\begin{pf}
By substituting the expressions \eqref{eq:Theta}
for $H$ and $\Theta$ in \eqref{eq:rho}
when $y$ is very large, and moving all multiples of $y$ to the left, we get
$$ (e_\min + \sum \frac{1}{m_p n_p} + e_\max) y 
   = a_\min + e_\min y_\min + \sum \frac{y_p}{m_p n_p} + e_\max y_\max.$$
Since the left term is linear in $y$ and the right term is constant,
both terms must be zero; this gives two linear relations in $e_\min$
and $e_\max$, of the form
$ e_\min + e_\max = c_1 $ and $ y_\min e_\min + y_\max e_\max = c_2$.
The constants $c_1$ and $c_2$ are expressions in
$m_p$, $n_p$, $a_\min$, and $a_\max$, hence are determined by the graph.
The coefficient matrix,
$ \left( \begin{array}{cc} -1 & -1 \\ y_\min & y_\max
  \end{array} \right),$
is determined by the graph, and is nonsingular because $y_\min \neq y_\max$.
Hence we can solve for $e_\min$ and $e_\max$,
and substitute in \eqref{eq:rho} to find $\rho(y)$.
\end{pf}

\begin{Lemma} \labell{rem:convex}
The density function of the \tDH measure
for a compact four dimensional Hamiltonian $S^1$-space
is concave on its support.

On the same manifold and circle action,
take any other closed two-form and a corresponding moment map,
and let $\rho(y)$ be the density function for its \tDH measure.
Then, in the notation set above,  
$\rho(y)$ is still concave on the interval $y_\min < y < y_\max$.

Additionally, if $a_\min$ and $a_\max$ are non-negative,
and if $y_\min <  y_p  < y_\max$ for all $p$,
the function $\rho(y)$ is non-negative for all $y$.
\end{Lemma}

\begin{pf}
By \eqref{eq:rho}, the restriction of the function
$\rho(y)$ to the interval $y_\min < y < y_\max$ is 
continuous, is piecewise linear, and its slope is decreasing.
Hence this restriction is concave. This proves the second claim.
The first claim is a special case.

Again by \eqref{eq:rho}, the limits of $\rho(y)$
as $y$ approaches the endpoints of the interval 
$[y_\min , y_\max]$
from its interior are $a_\min$ and $a_\max$. 
These limits are non-negative by assumption.
This and the concavity of $\rho(y)$ imply that $\rho(y)$ is non-negative 
for $y$ in the interval $(y_\min , y_\max)$.  
If $y_\min \leq y_p \leq y_\max$ for all $p$,
\eqref{eq:rho} implies that $\rho(y)$ is zero for $y < y_\min$ 
and is linear in $y$ for $y > y_\max$.  Since $\rho(y)$ is compactly 
supported, it must also be zero for $y > y_\max$.
\end{pf}

\section{Metrics}

\labell{sec:metrics}

The key to understanding compact four dimensional Hamiltonian 
$S^1$-spaces is an idea due to Audin, Ahara and Hattori:
choose a compatible Riemann metric and consider the gradient
flow of the moment map. The flow lines that lead from one fixed
point to another typically form a sphere, called a  ``gradient sphere".
The manifold and circle action are essentially determined by the 
arrangement of the gradient spheres. 

Our main point here is that
{\em the arrangement of the gradient spheres depends on the metric.}
This leads to repetitions in Audin's list of manifolds and circle actions;
different arrangements of gradient spheres could occur in the same 
manifold with different metrics.
In order to use the gradient spheres to distinguish between manifolds, 
one must resolve this ambiguity. We do this by showing that the arrangement 
of gradient spheres does not change if we avoid a sparse set of 
``bad metrics".

\subsection{Gradient spheres}

Let $(M,\omega)$ be a symplectic manifold with an action of a compact
Lie group $G$. A {\em compatible metric\/} on $M$ is a $G$-invariant
positive definite Riemann metric for which the endomorphism 
$J : TM \to TM$ defined by $\brak{u,v} = \omega(u,Jv)$ is an almost 
complex structure, i.e., satisfies $J^2 = \text{identity}$; such a $J$ 
is called a {\em compatible almost complex structure}.

\begin{Lemma} \labell{lem:metrics}
Every symplectic $G$-manifold admits a compatible metric. Moreover,
given any compatible metric on a neighborhood of an invariant closed 
subset $K \subseteq M$, there exists a compatible metric on $M$ 
that coincides with the given metric on some smaller neighborhood of $K$.
\end{Lemma}

\begin{pf}
Apply Lemma \ref{standard} to the symplectic vector bundle $E=TM$.
\end{pf}

Let $(M,\omega,\Phi)$ be a compact Hamiltonian $S^1$-space, 
equipped with a compatible metric.
The gradient vector field of the moment map is
\begin{equation} \labell{gradPhi}
	\grad \Phi = - J \xi_M,
\end{equation}
where $J$ is the corresponding almost complex structure
and $\xi_M$ is the vector field that generates the circle action.
Because the metric is invariant, the gradient flow commutes with the
circle action, so they fit together into an action of 
$\R \times S^1 \cong \Cc$, generated by the vector fields 
$\xi_M$ and $J\xi_M$.  Since these vector fields span a symplectic
subspace of $T_pM$ for each $p$, each orbit of the $\Cc$-action
is either a fixed point or is a two-dimensional symplectic sub-manifold 
of $M$.   The proof of the following is an easy exercise:

\begin{Lemma} \labell{S1 times I}
Each non-trivial $\Cc$-orbit is symplectomorphic to $S^1 \times I$,
with $I$ an open interval, with the circle acting on the left,
with a moment map $(e^{i\theta},h) \mapsto h$, and with a symplectic 
form $dh \wedge d\theta$.
\end{Lemma}

Each non-constant gradient trajectory approaches a top limit point
and a bottom limit point, both of which are fixed points for the action.
Therefore, the closure of a $\Cc$-orbit is a topological sphere;
it is called a {\em gradient sphere}. 
The circle acts on this sphere by rotations, fixing its north
and south poles. All the other points on the sphere have the same 
stabilizer, which is a finite, possibly trivial, subgroup of $S^1$.
A {\em free gradient sphere\/} is a gradient sphere whose stabilizer 
group is trivial.

A gradient sphere might not be smooth at its poles; see
\cite[Lemma 4.9]{ah-ha}.

All but a finite number of gradient spheres are free gradient spheres
whose north and south poles are extrema for the moment map. These spheres
are boring.  The interesting spheres, which we'll call {\em non-trivial},
are those with a finite non-trivial stabilizer, 
and those whose north or south pole is an interior fixed point.

{
\begin{Example}[K\"ahler toric varieties] \labell{ex:grad}
In a K\"ahler toric variety, the pre-image of an edge of the polygon
under the moment map for the torus action
is a two-sphere which is complex and invariant. 
Therefore, when we view the space as a Hamiltonian $S^1$-space,
this two-sphere is either fixed by the action, or is a gradient sphere 
for the K\"ahler metric.
So the arrangement of the gradient spheres with respect to the
K\"ahler metric is given exactly by the arrangement of the 
non-horizontal edges of the polygon.  
These arrangements, for the spaces in Figure \ref{polygons},
are illustrated in Figure \ref{ext-graphs}.
\end{Example}

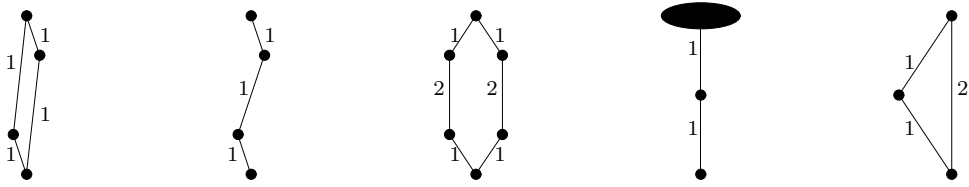
\begin{figure}
\begin{center}
\begin{picture}(20,70)(0,0)
\put(10,0){\circle*{4}}
\put(5,15){\circle*{4}}
\put(15,45){\circle*{4}}
\put(10,60){\circle*{4}}
\path(10,0)(5,15)(10,60)(15,45)(10,0)
\put(2,5){$\ss 1$}
\put(15,20){$\ss 1$}
\put(2,40){$\ss 1$}
\put(15,50){$\ss 1$}
\end{picture}
\hfill
\begin{picture}(20,70)(0,0)
\put(10,0){\circle*{4}}
\put(5,15){\circle*{4}}
\put(15,45){\circle*{4}}
\put(10,60){\circle*{4}}
\path(10,0)(5,15)(15,45)(10,60)
\put(1,5){$\ss 1$}
\put(5,30){$\ss 1$}
\put(15,50){$\ss 1$}
\end{picture}
\hfill
\begin{picture}(20,70)(0,0)
\put(10,0){\circle*{4}}
\put(0,15){\circle*{4}}
\put(20,15){\circle*{4}}
\put(0,45){\circle*{4}}
\put(20,45){\circle*{4}}
\put(10,60){\circle*{4}}
\path(10,0)(0,15)(0,45)(10,60)(20,45)(20,15)(10,0)
\put(-6,30){$\ss 2$}
\put(14,30){$\ss 2$}
\put(0,5){$\ss 1$}
\put(17,5){$\ss 1$}
\put(0,50){$\ss 1$}
\put(17,50){$\ss 1$}
\end{picture}
\hfill
\begin{picture}(20,70)(0,0)
\put(10,60){\blacken\ellipse{30}{10}}
\put(10,30){\circle*{4}}
\put(10,0){\circle*{4}}
\path(10,0)(10,30)(10,60)
\put(5,15){$\ss 1$}
\put(5,45){$\ss 1$}
\end{picture}
\hfill
\begin{picture}(20,70)(0,0)
\put(20,0){\circle*{4}}
\put(0,30){\circle*{4}}
\put(20,60){\circle*{4}}
\path(20,0)(20,60)(0,30)(20,0)
\put(22,30){$\ss 2$}
\put(2,15){$\ss 1$}
\put(2,40){$\ss 1$}
\end{picture}
\end{center}
\caption{Arrangements of gradient spheres for 
Example \protect\ref{ex:grad}}
\labell{ext-graphs}
\end{figure}
}


\subsection{Dependence on the metric}

\begin{Lemma} \labell{rigidity}
In a compact four dimensional Hamiltonian $S^1$-space
with a compatible metric:
\begin{enumerate}
\item every $Z_k$-sphere is a gradient sphere;
\item every non-free gradient sphere is a $Z_k$-sphere.
\end{enumerate}
\end{Lemma}

\begin{pf}
The $\Cc$-action preserves the set of points whose stabilizer is $Z_k$;
therefore, each $Z_k$-sphere is $\Cc$-invariant. 
Since the vector fields $\xi_M$ and $J\xi_M$ span the entire tangent
space of a $Z_k$-sphere at each point which is not a pole,
the $\Cc$-orbit of such a point is an open subset of the $Z_k$-sphere.
By connectedness, the whole $Z_k$-sphere minus its poles consists
of one $\Cc$-orbit. This implies part 1.

An orbit that is not fixed and not free is contained
in exactly one-$Z_k$ sphere and in exactly one gradient sphere.
By part 1, these spheres coincide. This proves part 2.
\end{pf}

The situation is quite different for free gradient spheres:

\begin{Lemma} \labell{non-rigidity} 
Fix a compact four dimensional Hamiltonian $S^1$-space 
and a compatible metric.  Let $C$ be a free gradient sphere
whose north and south poles, $p$ and $q$, are both interior fixed points.
Then there exists a smooth perturbation of the metric 
within the space of compatible metrics,
supported on an arbitrarily small neighborhood of an arbitrary
free orbit in $C$, such that for the perturbed metric, there
exists one free gradient sphere whose north pole is $p$ and whose
south pole is a minimum for the moment map, and another free gradient sphere 
whose south pole is $q$ and whose north pole is a maximum for the moment map,
and all other non-trivial gradient spheres are unchanged.
\end{Lemma}

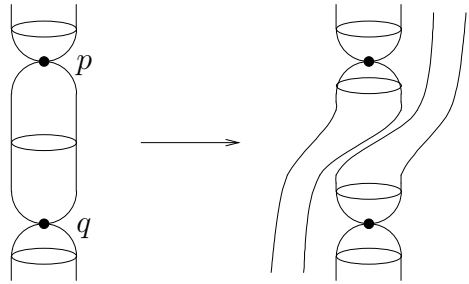
\begin{figure}
$$
\setlength{\unitlength}{0.0085in}
\begin{picture}(281,185)(0,-10)
\put(21.000,15.000){\arc{40.000}{3.1416}{6.2832}}
\put(21.000,55.000){\arc{40.000}{6.2832}{9.4248}}
\put(21.000,155.000){\arc{40.000}{6.2832}{9.4248}}
\put(21.000,115.000){\arc{40.000}{3.1416}{6.2832}}
\put(221.000,15.000){\arc{40.000}{3.1416}{6.2832}}
\put(221.000,55.000){\arc{40.000}{6.2832}{9.4248}}
\put(221.000,155.000){\arc{40.000}{6.2832}{9.4248}}
\put(221.000,115.000){\arc{40.000}{3.1416}{6.2832}}
\put(21,35){\blacken\ellipse{6}{6}}
\put(21,135){\blacken\ellipse{6}{6}}
\put(221,35){\blacken\ellipse{6}{6}}
\put(221,135){\blacken\ellipse{6}{6}}
\put(21,15){\ellipse{40}{10}}
\put(21,85){\ellipse{40}{10}}
\put(21,155){\ellipse{40}{10}}
\put(221,15){\ellipse{40}{10}}
\put(221,55){\ellipse{40}{10}}
\put(221,120){\ellipse{40}{10}}
\put(221,155){\ellipse{40}{10}}
\path(1,115)(1,55)
\path(41,115)(41,55)
\path(1,15)(1,0)
\path(41,15)(41,0)
\path(1,170)(1,155)
\path(41,170)(41,155)
\path(201,15)(201,0)
\path(241,15)(241,0)
\path(201,170)(201,155)
\path(241,170)(241,155)
\path(81,85)(141,85)
\path(133.000,83.000)(141.000,85.000)(133.000,87.000)
\path(241,115)	(241.000,111.734)
	(241.000,110.000)
\path(241,110)	(241.229,107.549) (241.000,105.000)
\path(241,105)	(239.243,102.028) (236.956,99.205) (234.214,96.513)
	(231.090,93.933) (227.660,91.447) (223.998,89.036)
	(220.178,86.681) (216.275,84.365) (212.364,82.069)
	(208.518,79.774) (204.813,77.462) (201.323,75.115)
	(198.122,72.713) (195.285,70.239) (192.886,67.674)
	(191.000,65.000)
\path(191,65)	(188.877,60.857) (187.062,56.107) (186.261,53.449)
	(185.526,50.573) (184.853,47.458) (184.239,44.080)
	(183.680,40.419) (183.172,36.452) (182.712,32.158)
	(182.296,27.513) (182.104,25.053) (181.921,22.497)
	(181.747,19.843) (181.582,17.088) (181.425,14.229)
	(181.276,11.263) (181.135,8.187) (181.000,5.000)
\path(201,115)	(201.000,111.734) (201.000,110.000)
\path(201,110)	(201.152,107.528) (201.000,105.000)
\path(201,105)	(198.745,100.253) (197.139,97.758) (195.278,95.200)
	(193.213,92.595) (190.992,89.958) (188.664,87.302)
	(186.278,84.643) (183.883,81.995) (181.529,79.372)
	(179.265,76.790) (177.139,74.264) (173.498,69.434)
	(171.000,65.000)
\path(171,65)	(169.444,60.917) (168.035,56.209) (167.375,53.566)
	(166.742,50.701) (166.130,47.591) (165.536,44.216)
	(164.956,40.553) (164.388,36.579) (163.827,32.274)
	(163.269,27.615) (162.990,25.146) (162.711,22.580)
	(162.431,19.915) (162.149,17.147) (161.866,14.275)
	(161.580,11.295) (161.292,8.204) (161.000,5.000)
\path(201,55)	(201.000,58.266) (201.000,60.000)
\path(201,60)	(200.771,62.451) (201.000,65.000)
\path(201,65)	(202.757,67.972) (205.044,70.795) (207.787,73.487)
	(210.910,76.067) (214.340,78.553) (218.002,80.964)
	(221.822,83.319) (225.725,85.635) (229.636,87.931)
	(233.482,90.226) (237.187,92.538) (240.677,94.885)
	(243.878,97.287) (246.715,99.761) (249.114,102.326)
	(251.000,105.000)
\path(251,105)	(253.123,109.143) (254.938,113.893) (255.739,116.551)
	(256.474,119.427) (257.147,122.542) (257.761,125.920)
	(258.320,129.581) (258.828,133.548) (259.288,137.842)
	(259.704,142.487) (259.896,144.947) (260.079,147.503)
	(260.253,150.157) (260.418,152.912) (260.575,155.771)
	(260.724,158.737) (260.865,161.813) (261.000,165.000)
\path(241,55)	(241.000,58.266) (241.000,60.000)
\path(241,60)	(240.848,62.472) (241.000,65.000)
\path(241,65)	(243.255,69.747) (244.861,72.242) (246.722,74.800)
	(248.787,77.405) (251.008,80.042) (253.336,82.698)
	(255.722,85.357) (258.117,88.005) (260.471,90.628)
	(262.735,93.210) (264.861,95.736) (268.502,100.566)
	(271.000,105.000)
\path(271,105)	(272.556,109.083) (273.966,113.791) (274.625,116.434)
	(275.259,119.299) (275.871,122.409) (276.465,125.784)
	(277.044,129.447) (277.613,133.421) (278.174,137.726)
	(278.732,142.385) (279.010,144.854) (279.289,147.420)
	(279.569,150.085) (279.851,152.853) (280.134,155.725)
	(280.420,158.705) (280.708,161.796) (281.000,165.000)
\put(41,130){$p$}
\put(41,30){$q$}
\end{picture}
$$
\caption{Breaking a free gradient sphere}
\end{figure}

\begin{pf}
Fix a free orbit in the gradient sphere $C$.
By Corollary \ref{stab=e}, on a neighborhood $U$ of this orbit
there exist coordinates
\begin{equation} \labell{coord}
	(e^{i\theta},h,x,y) : U \to S^1 \times I \times D^2,
\end{equation}
in which the moment map is $\Phi=h$, 
the $S^1$-action is multiplication on the left factor,
and the symplectic form is $\omega = dh \wedge d\theta + dx \wedge dy$,
and such that the points of $U \cap C$ are given by $x=y=0$.

We first perturb the gradient flow:
let $\rho$ be an $S^1$-invariant non-negative function
which is positive on our free orbit and is compactly supported in $U$.
Define 
$$ 
	\zeta_s = -J \xi_M + s \rho {\dd{x}} 
$$ 
on $U$ and $\zeta_s = -J \xi_M$ on $M \ssminus U$.
Then $\zeta_s$ is a smooth perturbation of the vector-field 
$\zeta_0 = \grad \Phi$ on $M$.
For $s > 0$, the flow-lines of $\zeta_s$ whose top limit point is $p$
are disjoint from those flow-lines whose bottom limit point is $q$. 
Moreover, if we choose the neighborhood $U$ to be small enough
so that it does not intersect any non-trivial gradient sphere except for $C$,
then for each flow-line of $\zeta_s$, if its top limit point is $p$, then
its bottom limit point is a minimum for $\Phi$, and if its bottom limit point 
is $q$, its top limit point is a maximum for $\Phi$.
Otherwise, the arrangements of the gradient spheres for 
$\grad \Phi = \zeta_0$ and for $\zeta_s$ are the same.

It remains to realize $\zeta_s$ as the gradient flow for some compatible 
metric.

Since $\omega(\zeta_s,\xi_M) = \omega(-J\xi_M,\xi_M) \neq 0$ on $U$,
$\span \{ \zeta_s , \xi_M \}$ is a symplectic sub-bundle of the tangent
bundle $TU$.   Let $E_s$ be its symplectic ortho-complement. 
Let $\eta_s, \mu_s$ be vector-fields on $U$ 
that form an oriented ortho-normal frame for $E_s$ (with respect 
to the compatible metric given by $J$ and $\omega$).
Define $J_s : TM \to TM$ by sending $J_s : \zeta_s \mapsto \xi_M$,
$\xi_M \mapsto -\zeta_s$, $\eta_s \mapsto \mu_s$, 
and $\mu_s \mapsto -\eta_s$ on $TU$, and by $J_s = J$ on $M \ssminus U$.
Then $J_s$ is an almost complex structure, it is compatible with $\omega$,
and its corresponding gradient flow is generated by $\zeta_s$.
\end{pf}

\begin{Corollary} \labell{generic}
Lemma \ref{non-rigidity} implies that for a generic compatible metric
there exists no free gradient sphere whose north and south poles are 
both interior fixed points. More precisely, the set of metrics
with this property is open and dense in the space of all compatible
metrics with the $C^\infty$ topology.
\end{Corollary}

\begin{Lemma}
Take a compact four dimensional Hamiltonian $S^1$-space.
For a generic compatible metric, the arrangement of gradient spheres
is determined by the graph.
\end{Lemma}

\begin{pf}
By Lemma \ref{rigidity} and by the construction of the graph,
two fixed points are connected by a non-free gradient sphere
if and only if their corresponding vertices in the graph
are connected by an edge.
By Corollary \ref{generic}, for a generic metric,
each interior fixed point whose positive isotropy weight
is $1$ is connected to the maximum of the moment map
by a free gradient sphere, each interior fixed point whose
negative isotropy weight is $-1$ is connected to the minimum
of the moment map by a free gradient sphere, and all other free 
gradient spheres are boring.
\end{pf}

\begin{Example} \labell{ex:sue}
Consider the K\"ahler toric variety that corresponds
to the second polygon in Figure \ref{polygons}.
View it as a Hamiltonian $S^1$-space, with its K\"ahler metric. 
As indicated in Figure \ref{ext-graphs}, this space
has a free gradient sphere whose north and south poles
are interior fixed points.  (See Example \ref{ex:grad}.)
If we perturb the K\"ahler metric as in Lemma \ref{non-rigidity},
we get a new arrangement of gradient spheres, 
which is the same as the one shown 
in Figure \ref{ext-graphs} on the left.
Hence the same Hamiltonian $S^1$-space can admit
two different compatible metrics for which the gradient spheres
are arranged differently.
\end{Example}

\section{Uniqueness: Graph determines space}

\labell{sec:uniqueness}

Let $(M,\omega,\Phi)$ and $(M',\omega',\Phi')$ be
two compact four dimensional Hamiltonian $S^1$ spaces. 
Any equivariant symplectomorphism $F : M \to M'$ that respects the 
moment maps induces an isomorphism on the corresponding graphs.
In this section we prove the converse:

\begin{Theorem}[Uniqueness Theorem] \labell{thm:uniqII}
Let $(M,\omega,\Phi)$ and $(M',\omega',\Phi')$ be two compact four
dimensional Hamiltonian $S^1$ spaces. Then any isomorphism between
their corresponding graphs is induced by an equivariant symplectomorphism.
\end{Theorem}

Audin made a beautiful observation, that a neighborhood
of the gradient spheres can be obtained by the plumbing
of disk bundles.  Since the result of plumbing 
is determined up to diffeomorphism by certain data,
which can be read from our graph, 
a neighborhood $U$ of the gradient spheres is determined by the graph.
To show that the graph determines the whole manifold, one
``sweeps" the neighborhood $U$ along the gradient flow
to fill almost the entire manifold, and then ``glues" the top.
See \cite{au:paper,au:book}.

Our first result, in \S\ref{sec:plumbing}, is that
if two manifolds correspond to the same graph,
there exists an orientation preserving equivariant diffeomorphism 
between them that also respects the moment maps.
This should require only a slight modification of Audin's argument.  
We nevertheless prove it from scratch,
because Audin left out differential topological details that
we felt should have been provided:
first, Audin's argument relies on the 
``uniqueness of equivariant plumbing", which involves standard
techniques in differential topology, but for which we were not
able to find a satisfactory reference; moreover, we need our plumbing
to keep track of additional information -- the moment map.
Second, we did not see how to avoid a certain variant of a theorem 
of Smale, on the diffeomorphisms of $S^2$, which Audin does not use
explicitly.
Ahara and Hattori refer to Smale's theorem in \cite[\S 10]{ah-ha},
but we believe that one needs an additional argument, which we provide 
in appendix \ref{sec:smale}. 
Third, the plumbing argument does not seem to work if the moment map has
an isolated minimum and there are more than two ``chains of gradient
spheres". (In contrary to \cite[Lemma 3.1.2]{au:paper}, this situation
can occur; see our Example \ref{ex:X} and \ref{ex:X2}).

Having constructed in \S\ref{sec:plumbing} an equivariant diffeomorphism 
that respects the moment map, we modify it in \S\ref{sec:uniquenessII},
using Moser's method, into one that also respects the symplectic form.

Before proceeding with the proof of Theorem \ref{thm:uniqII},
we give one application:

\begin{Proposition} \labell{application}
Let $M$ be a manifold with a circle action, let $\omega$ and $\omega'$
be equivariant symplectic forms for which the action is Hamiltonian,
and suppose that these forms represent the same de-Rham cohomology class.
Then $(M,\omega)$ is equivariantly symplectomorphic to $(M,\omega')$.
\end{Proposition}

\begin{pf}
By Lemma \ref{cohomology}, the integrals of $\omega$ and $\omega'$
over the fixed surfaces are the same, and the moment maps $\Phi$, $\Phi'$ 
can be chosen to have the same values at the fixed points. It follows 
that $(M,\omega,\Phi)$ and $(M',\omega',\Phi')$ have the same graph,
so, by the Uniqueness Theorem \ref{thm:uniqII}, they are isomorphic.
\end{pf}

\begin{Remark}
The condition on the action to be Hamiltonian for both $\omega$ 
and $\omega'$ is a simple condition: 
by McDuff \cite[Proposition 2]{md:four-manifolds}, 
a symplectic circle action on a symplectic four-manifold is Hamiltonian
if and only if it has fixed points.
\end{Remark}

\subsection{Building an equivariant diffeomorphism
that respects the moment maps} \labell{sec:plumbing}

\begin{Proposition} \labell{semi-morphism}
Let $(M,\omega,\Phi)$ and $(M',\omega',\Phi')$ be
compact four dimensional Hamiltonian $S^1$ spaces
whose graphs are isomorphic; fix an isomorphism between their graphs.
Then there exists an equivariant orientation preserving diffeomorphism 
$F: M \to M'$ that satisfies $\Phi' \circ F = \Phi$
and that induces the given isomorphism on the graphs.
\end{Proposition}

Note that
we {\em do not\/} yet require the map $F$ to respect the symplectic forms.

Throughout this section we fix the two spaces
and the isomorphism between their graphs. 
To each isolated fixed point, fixed surface, and $Z_k$-sphere
in $M$, there corresponds then an isolated fixed point, a fixed surface,
or a $Z_k$-sphere in $M'$.

\begin{pf*}{Preview of the proof of Proposition \ref{semi-morphism}}
By the local normal form, there exists an equivariant symplectomorphism 
on a neighborhood of the fixed point set. 
We extend it, first to neighborhoods of $Z_k$-spheres, then to a 
neighborhoods of non-trivial free gradient spheres,
while respecting the moment maps but not necessarily the symplectic forms.
After restricting to an open subset that is preserved under the 
descending gradient flow, we use this flow to extend our diffeomorphism 
to the complement of the maximal set of the moment map.  We then ``glue" it 
with another diffeomorphism, which is only defined near the maximum, 
using Smale's theorem if the maximum is isolated.
\end{pf*}

We now carry out this proof in a series of lemmas.

\begin{Lemma} \labell{near fixed points}
There exist open neighborhoods $U$ and $U'$ of the fixed point 
sets in $M$ and in $M'$ and an equivariant symplectomorphism 
$F : U \to U'$ that respects the moment maps.
\end{Lemma}

\begin{pf}
Since the isotropy weights are determined by the graph,
corresponding isolated fixed points have the same isotropy weights.
Hence, by Corollary \ref{stab=S1}, they have isomorphic neighborhoods.

Corresponding fixed surfaces are isomorphic, because a symplectic
surface is determined by its genus and total area, which are encoded
in the graph. Their neighborhoods are isomorphic, because, by Corollary
\ref{normal-near-B}, these neighborhoods are determined by the self
intersections of the surfaces, and, by Lemma \ref{lem:Euler}, these
self intersections, $e_\min$ and $e_\max$, are determined by the graph.
\end{pf}

\begin{Lemma} \labell{near Zk spheres}
There exist open subsets $U \subset M$ and $U' \subset M'$,
containing the fixed points and the $Z_k$-spheres,
and an equivariant diffeomorphism $F : U \to U'$
that respects the moment maps.
\end{Lemma}

\begin{pf}
Let $F'$ be the map of Lemma \ref{near fixed points},
defined near the fixed points.
Choose compatible metrics in $M$ and in $M'$ in such a way
that $F'$ is an isometry on a neighborhood of the fixed point set
and that a neighborhood of each isolated fixed point
is isometric to a disk in $\C^2$ as in \eqref{normalI};
this is possible by Lemma \ref{lem:metrics}.
Let $C$ be a $Z_k$-sphere in $M$, and let $C'$ be the corresponding
$Z_k$-sphere in $M'$.
Use Lemma \ref{CI} to identify neighborhoods 
of these spheres minus their poles with the model
$$
  S^1 \times_{Z_k} (I \times D^2).
$$
The map $F'$, where defined, has the form 
$$F'([\lambda,h,z]) = [a(h,z)\lambda,h,F_h(z)],$$
where $a(h,z) \in S^1$ and $F_h$ is a diffeomorphism between neighborhoods 
of the origin in $D^2$.
Since $F'$ commutes with the gradient flow,
and since the gradient flow only effects the $I$ coordinate,
$a(h,z)$ and $F_h(z)$, when defined, are locally independent of $h$.
Since the intersection of each gradient trajectory with a neighborhood
of the north pole of $C$ is connected (as follows by the model 
\eqref{normalI}), there exists a neighborhood of the origin in $D^2$
such that $a(h,z)$ and $F_h(z)$ are defined for all $z$ in this neighborhood
and all $h$'s which are sufficiently close to the endpoints of $I$.
The $F_h$'s are then equal to one diffeomorphism, $F_1$,
when $h$ is near the top of $I$, and to another, $F_2$,
when $h$ is near the bottom of $I$.
By standard differential topology, we can smoothly deform each of $F_1$
and $F_2$ to linear maps, and we can deform these linear maps into
each other. In this way we obtain a smooth family of diffeomorphisms
$F_h$, which coincide with our previous 
$F_h$'s for $h$ near the end-points of $I$. We extend $a(h,z)$ 
in a similar way.  The formula
$F([\lambda,h,z]) = [a(h,z)\lambda,h,F_h(z)]$ then defines an equivariant
diffeomorphism from a neighborhood of $C$ to a neighborhood of $C'$, 
which respects the moment maps, and which coincides with $F'$
on neighborhoods of the north and south poles of $C$.
\end{pf}

The proof of Lemma \ref{near Zk spheres} does not work for free
gradient spheres, for the following reason.
Let $C$ be a non-trivial free gradient sphere in $M$,
and let $C'$ be the corresponding sphere in $M'$.
Suppose that their north poles are interior fixed points
and that their south poles are minimal for the corresponding moment maps.
The map $F$ constructed in Lemma \ref{near fixed points} is defined on 
neighborhoods of the poles.
Since $F$ is an isometry, it interwines the gradient flows,
so it sends a neighborhood of the north pole in $C$ to a neighborhood
of the north pole in $C'$. However, a neighborhood of the south pole in $C$
will be sent to a neighborhood of the south pole in some free gradient 
sphere in $M'$, which might be different from $C'$.
This is because near the minimum of the moment map, all free gradient
spheres ``look alike";
it is impossible to tell which ones came from interior fixed points.

We return to our proof of Proposition \ref{semi-morphism}:

\begin{Lemma} \labell{near free spheres}
There exist generic metrics on $M$ and on $M'$, and there exist 
open subsets $U \subset M$ and $U' \subset M'$
that contain the fixed points, the $Z_k$-spheres,
and those free gradient spheres whose north poles are interior fixed points,
and there exists an equivariant diffeomorphism $F : U \to U'$
that respects the moment maps.
\end{Lemma}

\begin{pf}
Let $F'$ be the map of Lemma \ref{near Zk spheres},
defined near the fixed points and the $Z_k$-spheres.
Choose generic compatible Riemann metrics on $M$ and $M'$,
such that $F'$ is an isometry (possibly after shrinking its 
domain of definition), and such that neighborhoods of the 
isolated fixed points are isometric to the model \eqref{normalI}.

Let $p_1, \ldots, p_s$ be those interior fixed points in $M$
that have an isotropy weight $-1$, and let $p'_1, \ldots, p'_s$
be the corresponding points in $M'$.  Let $C_j$ and $C'_j$ be 
the free gradient spheres whose north poles are $p_j$ and $p'_j$.

Choose a real number $\beta$ greater than the minimum of $\Phi$
and close enough to it so that the interval $I' := (\min \Phi,\beta)$ 
contains no critical values of $\Phi$ and its pre-image 
is contained in the domain of definition of $F'$.

We will soon change $F'$ on the set $\Phi\inv(I')$ 
into an equivariant moment-map-preserving diffeomorphism $F$ 
that sends $C_j \cap \Phi\inv(h)$ to
$C'_j \cap (\Phi')\inv(h)$ for all $h$ near the top of $I'$
and that coincides with $F'$ on the intersection of $\Phi\inv(I')$
with a neighborhood of the $Z_k$-spheres.
Having done that, we can proceed exactly as in the proof
of Lemma \ref{near Zk spheres} to obtain an equivariant
moment-map-preserving diffeomorphism $F_j$ from a neighborhood
$U_j$ of $C_j \cap \{ \Phi > \beta - \epsilon \}$ onto a neighborhood $U'_j$
of $C'_j \cap \{ \Phi' > \beta - \epsilon \}$
and that coincides with $F$ on the intersection of $U_j$
with the set $\{ \beta - \epsilon < \Phi < \beta - \epsilon/2 \}$,
for some very small $\epsilon$.
The diffeomorphisms $F$ and $F_j$ fit together into a diffeomorphism 
with the desired properties, defined on the union of
$\{ \Phi < \beta - \epsilon/2\}$ with the $U_j$'s and with 
neighborhoods of the $Z_k$-spheres.

It remains to define $F$ on $\Phi\inv(I')$.
Since $F'$ interwines the gradient flow on $M$ with that on $M'$,
we can identify both $\Phi\inv(I')$ and $(\Phi')\inv(I')$ 
with a product $P \times I'$, where $P$ is a manifold with a locally
free circle action, in such a way that $F$ becomes
the identity map, the moment maps become the projection to $I'$,
and the gradient flows only effects the $I'$ coordinate.
The intersections of $C_j$ with the level sets $\Phi\inv(h)$, $h\in I'$, 
determine a single orbit $O_j$ in $P$. Similarly, $C'_j$ determines 
an orbit $O'_j$ in $P$, which coincides with $O_j$ exactly if $F'$ sends 
a neighborhood of the south pole in $C_j$ to a neighborhood of the
south pole in $C'_j$.
By a simple differential topological argument, which we provide
in Lemma \ref{P} below, there exists an equivariant diffeotopy
$F_h : P \to P$, $h \in I'$, supported away from the non-free orbits,
that deforms the identity map on $P$ into a map that sends each $O_j$ 
to the corresponding $O'_j$. We finish by defining $F(p,h) = (F_h(p),h)$.
\end{pf}

In the above proof we used the following standard 
differential topological lemma:

\begin{Lemma} \labell{P}
Let $P$ be a manifold with a locally free circle action.
Let $O_j$ and $O'_j$, $j=1,\ldots,s$, be free orbits in $P$,
with $O_i \neq O_j$ and $O'_i \neq O'_j$ for all $i \neq j$.
Then there exists a smooth family of equivariant maps, $F_h : P \to P$, 
parametrized by $h \in I'$, such that 
\begin{itemize}
\item for $h$ near the bottom of $I'$, the map $F_h$ is the identity map;
\item each non-free orbit in $P$ has a neighborhood on which 
all the $F_h$'s are the identity map;
\item for $h$ near the top of $I'$, the map $F_h$ sends $O_j$ to $O'_j$.
\end{itemize}
\end{Lemma}

\begin{pf}
Let $P'$ be a connected invariant open subset of $P$
which is equivariantly diffeomorphic to $S^1 \times B$
for some two dimensional manifold $B$.
(Such a $P'$ is obtained by removing the closures of small tubular 
neighborhoods of any finite non-empty set of orbits in $P$ that includes 
all the non-free orbits.)
Let $q_j$ and $q'_j$ be the points of $B$ whose preimages in $P'$
are the orbits $O_j$ and $O'_j$.  For each $j$, let $\gamma_j$ be 
a path from $q_j$ to $q'_j$, and let $U_j$ be a neighborhood of this
path, such that the $U_j$'s are disjoint from each other.
By standard differential topology (see Thom's isotopy embedding
theorem in \cite[chap.9]{BJ}), there exists a diffeotopy of $B$, 
compactly supported in the union of the $U_j$'s, 
that connects the identity map to a map that sends each $q_j$ to $q'_j$.
This trivially lifts and extends to a diffeotopy of $P$ 
with the required properties.
\end{pf}

We would like to restrict $F$ to an open set that is invariant
under the downward gradient flow of the moment map. Denote the gradient
flow by $g_t : M \to M$, i.e., 
$\frac{d}{dt}g_t = \grad \Phi \circ g_t$.

\begin{Lemma} \labell{flow down invariant}
Let $U \subset M$ be an invariant open subset that contains 
the minimal set of the moment map and that contains all gradient spheres
whose north pole is an interior fixed point.
Then there exists an invariant open subset $U'$ of $U$, with these same 
properties, for which $g_{-t} (U') \subseteq U'$ for all $t \geq 0$.
\end{Lemma}

\begin{pf}
We construct the set $U'$ inductively, ``from the bottom up".
Let $p_1, \ldots, p_s$ be the interior fixed points in $M$, ordered
in such a way that if there is a gradient sphere with south pole $p_i$
and north pole $p_j$, then $i < j$.
We will construct open sets 
$U_0 \subseteq U_1 \subseteq \ldots \subseteq U_s$,
such that $g_{-t}(U_j) \subseteq U_j$ for all $t \geq 0$,
and such that $U_j$ contains $p_1, \ldots, p_j$.

Take $U_0 = \Phi\inv(-\infty,\beta)$, where $\beta > \min \Phi$
is sufficiently small so that $U_0 \subseteq U$.

Suppose that we have already constructed $U_{j-1}$.
If $p_j$ is contained in $U_{j-1}$, set $U_j = U_{j-1}$.
Otherwise, let $C_j$ be the gradient sphere descending from $p_j$,
and let $q_j$ be its south pole.
By the local normal form of Lemma \ref{CI}, 
there exists a neighborhood $W$ of $C_j \ssminus \{ p_j,q_j \}$
in which each descending gradient trajectory approaches the level set 
$\Phi\inv (\Phi (q_j))$. Since $U_{j-1}$ is open and contains $q_j$,
it also contains a neighborhood of $q_j$ in its level set, and so
we can shrink $W$ so that
the bottom of each gradient trajectory in $W$ will intersect $U_{j-1}$.
The union $U'_j = U_{j-1} \cup W$ satisfies 
$g_{-t} (U'_j) \subseteq U'_j$ for all $t > 0$ and contains 
$C_j \ssminus \{p_j\}$.

By the local normal form \eqref{normalI} for a neighborhood of $p_j$, 
there exists a disk $W'$ around $p_j$ in which each gradient trajectory 
either approaches $p_j$ or intersects the set 
$\{ \Phi < \Phi(p_j) \}$.
By shrinking the disk $W'$, we can arrange that each gradient trajectory
either approaches $p_j$ or intersects the set $U'_j$ constructed above. 
The open set $U_j = U'_j \cup W'$ satisfies the required properties.
\end{pf}

\begin{Lemma} \labell{minus maximum}
There exists an equivariant diffeomorphism, respecting the moment maps,
from the subset $\{ \Phi < \max \Phi \}$ of $M$
onto the subset $\{ \Phi' < \max \Phi' \}$ of $M'$.
\end{Lemma}

\begin{pf}
Let $F':U \to U'$ be the map of Lemma \ref{near free spheres}.
Pull back the metric on $U'$ to a metric on $U$, and extend it
to a metric on $M$, possibly after shrinking $U$ and $U'$. 
With this new generic metric, $F'$ is an isometry, and $U$ and $U'$
still contain all free gradient spheres whose north poles are interior
fixed points.

By Lemma \ref{flow down invariant}, we can shrink $U$ and $U'$ further 
so that they would be invariant under the downward gradient flow.
Then $F'$ has a unique extension to an equivariant diffeomorphism $F$
on $\{ \Phi < \max \Phi \}$ that interwines the gradient flow on $M$ 
with that on $M'$ and satisfies $\Phi' \circ F = \Phi$.
\end{pf}

\begin{pf*}{Proof of Theorem \ref{thm:uniqII}}
Let $U_- \subset M$ and $U'_- \subset M'$ be the complements of the 
maximal sets of the moment maps. In Lemma \ref{minus maximum}
we constructed an equivariant diffeomorphism
$F_- : U_- \to U'_-$ which respects the moment maps.

First assume that the maximum of $\Phi$ is attained on a surface.
The normal form of Corollary \ref{normal-near-B} implies that
there is a principal $S^1$-bundle $P \to B$ such that
neighborhoods of the maximum in $M$ and in $M'$ are both given by a model 
\begin{equation} \labell{near max}
  P \times_{S^1} D^2,
\end{equation}
in which the moment map is $[p,z] \mapsto \max \Phi - |z|^2$.
The map $F_-$ is defined on this model minus the zero section. 
Assume that we have chosen our metric to be standard 
on this model.  Then the facts that $F_-$ preserves the moment map 
and interwines the gradient flows implies that it is induced
by an equivariant diffeomorphism of $P$, and therefore extends
to the whole neighborhood \eqref{near max}. 

Now assume that $\Phi$ has an isolated maximum.
The local normal form near the maximum \eqref{normalI}
gives an isomorphism $F_+ : U_+ \to U'_+$
between the sets of points whose moment image
is greater than $\max(\Phi) - \epsilon$ for some $\epsilon >0$.
To glue $F_+$ and $F_-$, it is enough to find an equivariant 
diffeomorphism $F: U_+ \cap U_- \to U'_+ \cap U'_-$
that respects the moment maps,
that coincides with $F_+$ on a neighborhood of the maximum,
and that coincides with $F_-$ on a neighborhood of 
$\Phi\inv( \{ \max(\Phi) - \epsilon \} )$.
As in the proof of Lemma \ref{near free spheres},
we use the gradient flows and the map $F_-$ to
identify $U_+ \cap U_-$ and $U'_+ \cap U'_-$
with the model $P \times I$
where $P$ is a level set of $\Phi$ and where 
$I$ is the interval $(\max(\Phi) - \epsilon , \max(\Phi))$.
On this model, $F_-$ is an identity map, and $F_+$
is given by an equivariant diffeomorphism of $P$ that is
independent of the $I$ coordinate. 
So it is enough to show that each equivariant diffeomorphism
$f$ of $P$ that sends each non-free orbit to itself
can be connected to the identity by an equivariant diffeotopy.

Let $\ol{f}$ be the induced diffemorphism on the quotient orbifold,
$B := P / S^1$.  As in \cite[\S 10]{ah-ha}, we can diffeotope
$\ol{f}$ into a map that is the identity on a neighborhood 
of each orbifold singularity. 
The orbifold $B$ is a two-sphere with at most two singular points;
this follows from the local model \eqref{normalI}.
Fix a homeomorphism from $B$ to $S^2$ that is smooth outside 
the singular points;
the images of the singular points in $B$ are marked points in $S^2$. 
The diffeomorphisms of $B$ that are the identity on neighborhoods of the 
singular points become the diffeomorphisms of $S^2$ that are fixd on
neighborhoods of the marked points.
If the number of marked points is one or zero, Smale's theorem \cite{s} 
essentially tells us that every such diffeomorphism
is diffeotopic, again through such diffeomorphisms, to the identity.
(See appendix \ref{sec:smale}.) 
If there are two marked points, this is no longer true.
However, every diffeomorphism of $S^2$ that fixes neighborhoods
of these marked points is diffeotopic to the identity through
diffeomorphisms that act on neighborhoods of these points
by rotations (again, see appendix \ref{sec:smale}).
This is enough because we can identify $B$ with $S^2$ 
in such a way that a rotation of a neighborhood of a marked point in $S^2$ 
becomes a rotation (in particular, is smooth) 
on a neighborhood of the corresponding singular point in $B$.
In this way we obtain a diffeotopy $\ol{f}_t$ of $B$,
with $\ol{f}_0 = \ol{f}$ and $\ol{f}_1 = $identity.
This lifts to an equivariant diffeotopy $f_t$ of $P$
with $f_0 = f$ and such that $f_1$ rotates each fiber of $P \to B$
by an element of $S^1$, given by a smooth function 
$h : B \to S^1$. Since $h$ can be smoothly deformed to the constant map
$1$ (because $B$ is a 2-sphere), the map $f_1$ can be smoothly deformed 
to the identity.
\end{pf*}

\subsection{Building an isomorphism} \labell{sec:uniquenessII}

Proposition \ref{semi-morphism} reduces Theorem \ref{thm:uniqII}
to the following proposition, which was conjectured by Delzant
and Ginzburg:

\begin{Proposition} \labell{prop:conj-del}
Take a compact manifold $M$ with a circle action and 
a smooth map $\Phi : M \to \R$.
Let $\omega$ and $\omega'$ be two invariant symplectic forms for which
$\Phi$ is the moment map, i.e., such that (\ref{eq:hamilton}) is
satisfied for both $\omega$ and $\omega'$.
Suppose that the integrals of $\omega$ and $\omega'$ over any fixed surface
are the same.  Then there exists an equivariant 
diffeomorphism $F: M \to M$ such that $F^* \omega = \omega'$
and $\Phi' \circ F = \Phi$.
\end{Proposition}

The proof of the proposition will use the following Lemma,
which is also interesting in itself:

\begin{Lemma} \labell{same-coh}
Under the assumptions of Proposition \ref{prop:conj-del},
the de Rham cohomology classes of $\omega$ and $\omega'$ 
are the same.
\end{Lemma}

\begin{pf}
The moment map $\Phi$ is a perfect Bott-Morse function 
\cite[\S5]{g-s:convexity}.  The gradient spheres, with respect
to any invariant metric, form {\em completing cycles\/} for the
critical points. Therefore, the gradient spheres,
together with the minimal set of $\Phi$ if that happens to be a surface,
form a basis of $H_*(M,\Z)$.
The assumptions of Proposition \ref{prop:conj-del} imply that
the integrals of $\omega$ and $\omega'$ on these basis elements 
are the same. (The symplectic area of a gradient sphere 
with isotropy $k$ is ${1 \over k} (\Phi(p) - \Phi(q))$
where $p$ and $q$ are the north and south poles;
see the proof of Lemma \ref{cohomology}.)
\end{pf}

\begin{pf*}{Proof of Proposition \ref{prop:conj-del}}

Let $\omega_t = (1-t) \omega + t \omega'$, for $0 \leq t \leq 1$.
By Lemma \ref{same-coh}, the cohomology classes $[\omega_t]$
are all the same.  By Moser's method \cite[lecture 5]{w:lectures},
to prove Proposition \ref{prop:conj-del} it is sufficient
to prove that $\omega_t$ is nondegenerate for all $t$.
This can be checked pointwise.
There are two cases:

\subsubsection*{Case 1: At a point with a discrete stabilizer} 

%
Let $\xi_M$ be the vector field on $M$ that generates the circle action,
and suppose that $\xi_M |_p \neq 0$.
In the vector space $T_pM$, choose an oriented basis of
the form $x, y, u, v$ where $\Phi_*x =1$, $y=\xi_M|_p$, and
$\Phi_* u = \Phi_*v =0$.
Since $\omega_t(\zeta,y) = \Phi_*\zeta$ for all $\zeta$,
the matrix representing $\omega_t|_{T_pM}$ has the form
$$\left(\begin{array}{rrrr}
  0 &  1 &  a_t & b_t \\
 -1 &  0 &  0 & 0 \\
 -a_t &  0 &  0 & c_t \\
 -b_t &  0 & -c_t & 0 
\end{array}\right).$$
We have $\omega_t \wedge \omega_t = 2c_t$ times 
a generator of $\wedge^4 T^*_p M$.
Since $[\omega] = [\omega']$, the symplectic forms 
$\omega$ and $\omega'$ induce the same orientation on $M$,
so $c_0$ and $c_1$ must have the same sign.
This further implies that $c_t = (1-t) c_0 + t c_1$ is nonzero,
hence $\omega_t \wedge \omega_t \neq 0$.

\subsubsection*{Case 2: At a fixed point} 
Let $p$ be a fixed point.  The Hessian of $\Phi$ at $p$
provides a quadratic moment map for the linear symplectic
form $\omega_p$ on $T_pM$ with the linear isotropy action.
Our Proposition hence follows from the following lemma:

\begin{Lemma} \labell{omega-unique}
Let $S^1$ act on $\C^2$ by
$\lambda \cdot (z,w) = (\lambda^m z, \lambda^n w)$
where $m,n$ are integers, not both zero.
Let $\omega$ and $\omega'$ be $S^1$-invariant 
linear symplectic forms on $\C^2$, compatible with the complex
orientation, and with the same quadratic moment map $\Phi: \C^2 \to \R$.
Then $\omega_t := (1-t) \omega + t \omega'$ is nondegenerate for all 
$0 \leq t \leq 1$.
\end{Lemma}

\begin{pf}
The circle action is generated by the vector field
$$\xi_M = i m (z {\partial \over \partial z} 
              - \ol{z} {\partial \over \partial \ol{z}}) 
       + i n (w {\partial \over \partial w} 
              - \ol{w} {\partial \over \partial \ol{w}}).$$ 
Any real 2-form can be written as 
$$\omega = {i \over 2}
	      ( A dz \wedge d\zbar + 
		B dz \wedge dw - \Bbar d\zbar \wedge d\wbar +
		C dz \wedge d\wbar - \Cbar d\zbar \wedge dw +
		D dw \wedge d\wbar)$$ 
where $A,D \in \R$ and $B,C \in \C$. 
The $S^1$-invariance amounts to 
\begin{equation} \labell{eq:equiv-cond}
\lambda^{m+n} B = B \quad \mbox{and} \quad \lambda^{m-n} C = C \quad 
\mbox{for all } \lambda \in S^1.
\end{equation}
We will consider four sub-cases.
In the first three sub-cases below, we will show that the symplectic form
$\omega$ is determined from the function $\Phi$ and from the integers
$m$ and $n$. This implies that all the $\omega_t$'s are the same and 
are therefore nondegenerate.
In the fourth sub-case, we might have $\omega \neq \omega'$, but 
$\omega_t$ will still be nondegenerate.

\subsubsection*{Sub-case 1: $m = n \neq 0$.}  
Then (\ref{eq:equiv-cond}) implies that $B=0$, so
$$
	\omega = {i \over 2} ( A dz \wedge d\zbar + 
	  C dz \wedge d\wbar - \Cbar d\zbar \wedge dw + 
	  D dw \wedge d\wbar). 
$$
The moment map is 
$$
	\Phi(z,w) = -{m \over 2} (A |z|^2 + D |w|^2 
			+ C z \wbar + \Cbar \zbar w).
$$ 
From this we can extract $A$, $D$, and $C$.

\subsubsection*{Sub-case 2: $m=-n \neq 0$.} 
Then (\ref{eq:equiv-cond}) implies $C=0$, so
$$
	\omega = {i \over 2} ( A dz \wedge d\zbar 
	+ B dz \wedge dw - \Bbar d\zbar \wedge d\wbar 
	+ D dw \wedge d\wbar). 
$$
The moment map is 
$$
	\Phi(z,w) = -{m \over 2} (A |z|^2 - D |w|^2 
			+ B z w + \Bbar \zbar \wbar). 
$$
From this we can extract $A$, $D$, and $B$.

\subsubsection*{Sub-case 3: $m \neq \pm n$, both $\neq 0$.}
Then (\ref{eq:equiv-cond}) implies $B=C=0$, so 
$$
	\omega = {i \over 2} (A dz \wedge d\zbar + D dw \wedge d\wbar).
$$
The moment map is 
$$
	\Phi(z,w) = -{1 \over 2}(m A |z|^2 + n D |w|^2).
$$
From this we can extract $A$ and $D$.

\subsubsection*{Sub-case 4: $m=0$, $n \neq 0$.}  
Again, (\ref{eq:equiv-cond}) implies $B=C=0$, so 
$$
	\omega = {i \over 2} (A dz \wedge d\zbar + D dw \wedge d\wbar).
$$
The moment map is 
$$
	\Phi(z,w) = -{1 \over 2} n D |w|^2.
$$
From this we can extract $D$.
The coefficients $A$ and $A'$ in $\omega$ and in $\omega'$
might be different, but their signs must both be equal to the
sign of $D$. So
$\omega_t = {i \over 2} (A_t dz \wedge d\zbar + D dw \wedge d\wbar)$,
with $A_t = (1-t) A + t A' \neq 0$, is nondegenerate.
\end{pf}
This completes the proof of Proposition \ref{prop:conj-del},
and hence of Theorem \ref{thm:uniqII}.
\end{pf*}

\begin{Example} \labell{hirz}
The Hamiltonian $S^1$-spaces
corresponding to the first two polygons in Figure \ref{polygons} 
have the same graph (see Example \ref{ex:polygons}), hence are isomorphic.

This shows that a Hamiltonian $S^1$-space can admit two nonisomorphic
compatible K\"ahler metrics (nonisomorphic because their gradient spheres 
are arranged differently; see Example \ref{ex:grad}),
and that a Hamiltonian $S^1$-action can extend to two nonisomorphic
toric actions (their Delzant polygons are different).

We note that the first space is $S^2 \times S^2$, and the 
second is a Hirzebruch surface, $M=P \times_{S^1} S^2$, 
where $P \to S^2$ is a principal circle bundle 
with Chern number $2$.
\end{Example}

\section{Isolated fixed points implies toric variety}

\labell{sec:toric}

We have shown that a compact four dimensional Hamiltonian $S^1$-space
is determined by its graph. To complete the classification, it remains
to determine which such spaces (or graphs) occur. 

In this section we treat the
case of isolated fixed points.  We show that these are all K\"ahler
toric varieties, considered as Hamiltonian $S^1$-spaces as in 
\S\ref{subsec:delzant}.
The Uniqueness Theorem \ref{thm:uniqII} reduces the proof to combinatorics: 
we only need to show 
that the graph of every such a space comes from a Delzant polygon.
This result has two important consequences, Corollaries 
\ref{cor:isol-torus} and \ref{cor:isol-kahler}:
first, the circle action extends to a Hamiltonian action of a 
two dimensional torus;
second, the space is K\"ahler, i.e., it admits a complex
structure that is compatible with the given symplectic structure
and is invariant under the circle action.

\begin{Definition}
Let $(M,\omega,\Phi)$ be a compact four dimensional Hamiltonian
$S^1$-space, equipped with a compatible metric. A {\em chain of
gradient spheres\/} is a sequence of gradient spheres, $C_1, \ldots, C_l$,
such that the south pole of $C_1$ is a minimum for the moment map,
the north pole of $C_{i-1}$ is the south pole of $C_i$ for each 
$1 < i \leq l$,
and the north pole of $C_l$ is a maximum for the moment map.
A chain is {\em non-trivial\/} if it contains more than one sphere,
or if it contains one sphere whose stabilizer is non-trivial.
\end{Definition}

\begin{Lemma} \labell{lem:conditions} 
Let $C_1, \ldots, C_l$ be a chain of gradient spheres,
and let $k_1, \ldots, k_l$ be the orders of their stabilizers.
(See Figure \ref{branch}.) Then 
\begin{equation} \labell{conditions}
\parbox{5.0in}{
\begin{enumerate}
\item
$\gcd(k_i,k_{i+1})=1$ for $i=1,\ldots,l-1$, and
\item
$(k_{i-1} + k_{i+1})/k_i$ is an integer for $i=2,\ldots,l-1$.
\end{enumerate}
}
\end{equation}
Moreover, $-(k_{i-1} + k_{i+1})/k_i$ is the self intersection of $C_i$. 
\end{Lemma}

\begin{figure}
$$
\setlength{\unitlength}{0.0100in}
\begin{picture}(46,129)(-10,-10)
\put(30,102){\blacken\ellipse{4}{4}}
\put(20,72){\blacken\ellipse{4}{4}}
\put(20,32){\blacken\ellipse{4}{4}}
\put(30,2){\blacken\ellipse{4}{4}}
\path(30,102)	(28.809,98.922) (27.710,96.056) (26.700,93.390)
	(25.775,90.914) (24.167,86.487) (22.856,82.688)
	(21.814,79.427) (21.010,76.617) (20.000,72.000)
\path(20,72)	(19.506,67.999) (19.153,63.141) (19.030,60.481)
	(18.942,57.713) (18.889,54.875) (18.871,52.000)
	(18.889,49.125) (18.942,46.287) (19.030,43.519)
	(19.153,40.859) (19.506,36.001) (20.000,32.000)
\path(20,32)	(21.010,27.383) (21.814,24.574) (22.857,21.313)
	(24.168,17.513) (25.776,13.086) (26.700,10.610)
	(27.710,7.945) (28.809,5.078) (30.000,2.000)
\put(-7,82){$k_{i+1}$}
\put(0,47){$k_i$}
\put(-5,7){$k_{i-1}$}
\put(30,72){$\ss p_i$}
\put(30,32){$\ss p_{i-1}$}
\end{picture}
$$
\caption{A chain of gradient spheres}
\labell{branch}
\end{figure}
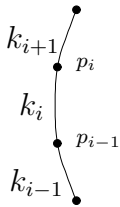

\begin{pf}
Denote by $p_i$ the north pole of $C_i$, which is also the south pole
of $C_{i+1}$.  The isotropy weights at $p_i$ are $k_{i+1}$ and $-k_i$.
Since the circle action is effective (on $M$, and hence on a
neighborhood of $p$), these integers are relatively prime.
This proves part 1.

The rest of the lemma follows from a certain fact about $S^1$-equivariant 
line bundles over $S^2$ that we prove in Lemma \ref{mg} below.
The normal bundle of $C_i$ in $M$ can be viewed 
as an $S^1$-equivariant complex line bundle (Corollary \ref{s=c}).
The circle rotates the two-sphere $k_i$ times, and 
acts on the fibers over the north and south poles
with weights $k_{i+1}$ and $-k_{i-1}$.
By Lemma \ref{mg}, the self intersection of the zero section is 
$- (k_{i-1} + k_{i+1} ) / k_i$.
\end{pf}

\begin{Lemma}  \labell{mg}
Let $S^1$ act on $S^2$ by rotating it $k$ times while fixing the north 
and south poles. Suppose that the action lifts to a complex line bundle.
Then $S^1$ acts linearly on the fibers over the north and south poles;
let $m$ and $n$ be the weights for these actions.  Then 
$$
	m - n = -ek
$$ 
where $e$ is the self intersection of the zero section.
\end{Lemma}

\begin{pf}
See \cite[lemma 4.3]{ah-ha}.
It is also not hard, and is a worthwhile exercise, to prove this directly:
if we decompose $S^2$ into the upper and lower hemispheres and trivialize
the bundle over each hemisphere, the self intersection $e$ is equal to the
winding number of the gluing map, $\psi: \text{equator} \to \C^\times$,
which one can compute.
\end{pf}

\begin{Corollary}
Any gradient sphere whose north and south poles are interior fixed
points has a negative self intersection.
\end{Corollary}

{\em Fans\/} are central ingredients in the theory of toric varieties.
We will use them not directly to construct toric varieties,
but as an aid in our computations.
The following definition of cones and fans in $\R^2$ is less general 
than the standard definition, but is sufficient for our needs:

\begin{Definition}
A {\em smooth cone\/} $c$ in $\R^2$ is the positive span of
two integer vectors $u = (k,b)$ and $u' = (k',b')$ 
satisfying $kb'-bk'=1$; we say that $c$ is {\em generated\/} 
by $u$ and $u'$.
A {\em smooth fan\/} in $\R^2$ is a collection of smooth cones in which
the intersection of any two different cones is either the origin or is 
a common ray.
A smooth fan in $\R^2$ is {\em complete\/} if the union of its cones
is $\R^2$.  
We say that
a complete smooth fan $c_1,\ldots,c_N$ in $\R^2$ is {\em generated\/} by
$u_1, \ldots, u_N$ if $c_i$ is generated by $u_i$ and $u_{i+1}$ for all
$i$, with the indices taken cyclically. 
\end{Definition}

\begin{Lemma} \labell{mivi}
Let $k_1,\ldots,k_l$ be positive integers that satisfy \eqref{conditions}.
Then there exist integers $b_1,\ldots,b_l$ such that
\begin{equation} \labell{eq:det=1}
k_i b_{i+1} - b_i k_{i+1} =1 \quad \mbox{for all $i$}.
\end{equation}
\end{Lemma}

\begin{pf}
By subtracting the equations
\begin{eqnarray*}
k_{i-1} b_i - b_{i-1} k_i & = & 1	\\
k_{i} b_{i+1} - b_i k_{i+1} & = & 1	,
\end{eqnarray*}
we get the requirement
\begin{equation}\labell{yi}
b_{i+1} = -b_{i-1} + { k_{i+1} + k_{i-1} \over k_i} b_i.
\end{equation}
Choose $b_1,b_2 \in \Z$ such that $k_1 b_2 - b_1 k_2 =1$;
this is possible because $\gcd(k_1,k_2)=1$.  
Construct the other $b_i$'s recursively by (\ref{yi}).
\end{pf}

\begin{Lemma} \labell{kho}
Let $k_1,\ldots,k_l$ be positive integers that satisfy \eqref{conditions}.  
Then 
$$ {1 \over k_1 k_2 } + \ldots + {1 \over k_{l-1} k_l} 
 = {d \over k_1 k_l}$$
for some positive integer $d$.
\end{Lemma}

\begin{pf}
Consider the vectors $u_i = (k_i,b_i)$,
where $b_i$ are as in Lemma \ref{mivi}.
Denote by $c_i$ the cone in $\R^2$ generated by $u_{i}$ and $u_{i+1}$. 
These cones are arranged in counterclockwise order 
(because $\det(u_i u_{i+1}) = k_i b_{i+1} - b_i k_{i+1} >0$),
and they do not complete a full turn around the origin and eventually 
overlap, because they remain in the right half plane (because $k_i > 0$).
Therefore they fit together to form a smooth fan in $\Z^2$, 
as in Figure \ref{fan}.
An elementary computation shows that
$$
	\int_{c_i} e^{-x} dx \wedge dy = 
 	{1 \over k_{i} k_{i+1}} \det (u_i u_{i+1}) = 
	{1 \over k_{i} k_{i+1}}.
$$
Similarly, using the fact that the counterclockwise angle
from $u_1$ to $u_l$ is less than $180^\circ$,
$$\int_{c_1 \cup \cdots \cup c_{l-1}} e^{-x} dx \wedge dy 
  = {1 \over k_1 k_l} \det(u_1 u_l).$$
The desired equality follows from the additivity of the integral.
\end{pf}

\begin{figure}
$$
\setlength{\unitlength}{0.0050in}
\begin{picture}(240,315)(0,-10)
\texture{c0c0c0c0 0 0 0 0 0 0 0 
	c0c0c0c0 0 0 0 0 0 0 0 
	c0c0c0c0 0 0 0 0 0 0 0 
	c0c0c0c0 0 0 0 0 0 0 0 }
\shade\path(0,200)(200,300)(240,200)(0,200)
\shade\path(0,200)(240,80)(200,0)(0,200)
\texture{c0c0c0c0 0 0 0 0 0 0 0 
	c0c0c0c0 0 0 0 0 0 0 0 
	c0c0c0c0 0 0 0 0 0 0 0 
	c0c0c0c0 0 0 0 0 0 0 0 }
\shade\path(0,200)(240,200)(240,80)(0,200)
\thicklines
\path(0,200)(80,120)
\path(0,200)(80,120)
\path(65.858,128.485)(80.000,120.000)(71.515,134.142)
\path(0,200)(160,120)
\path(0,200)(160,120)
\path(143.900,123.578)(160.000,120.000)(147.478,130.733)
\path(0,200)(80,200)
\path(0,200)(80,200)
\path(64.000,196.000)(80.000,200.000)(64.000,204.000)
\path(0,200)(160,280)
\path(0,200)(160,280)
\path(147.478,269.267)(160.000,280.000)(143.900,276.422)
\put(135,280){$\ss u_4$}
\put(75,205){$\ss u_3$}
\put(155,125){$\ss u_2$}
\put(85,120){$\ss u_1$}
\put(155,225){$c_3$}
\put(190,145){$c_2$}
\put(170,65){$c_1$}
\end{picture}
$$
\caption{A fan.}
\labell{fan}
\end{figure}

\begin{Lemma} \labell{lem:DH} 
Let $(M,\omega,\Phi)$ be a four dimensional compact Hamiltonian 
$S^1$-space, equipped with a compatible metric.
Suppose that the fixed points are isolated, and that the isotropy
weights are $-m,-m'$ at the maximum and $n,n'$ at the minimum.
Consider those chains of gradient spheres that consist of two or more
spheres; let $m_j$ and $n_j$ be the orders of the stabilizers
of the top and bottom spheres in the $j$'th such chain.
Then there exist positive integers $d_j$ such that
\begin{equation} \labell{L=R}
 {1 \over m m'} + { 1 \over n n'} = \sum_j {d_j \over m_j n_j } ,
\end{equation}
summing over those chains that contain two or more spheres.
\end{Lemma}

\begin{pf}
Let $k_i^j$ be the order of the stabilizer of the $i$th gradient sphere
in the $j$th chain.  The isotropy weights at the $i$th fixed point
in this chain are $ k_{i+1}^j$ and $-k_i^j $.
By Lemma \ref{cor:DH},
$
	{ 1 \over mm'} + {1 \over nn'} = 
 		\sum_j \sum_i {1 \over k_i^j k_{i+1}^j }.
$
By Lemma \ref{kho}, 
$ 
	\sum_i { 1 \over k_i^j k_{i+1}^j} = { d_j \over m_j n_j}.
$
\end{pf}

\begin{Proposition} \labell{lem:two-branches}
Let $(M,\omega,\Phi)$ be a four dimensional compact Hamiltonian $S^1$-space
with isolated fixed points, equipped with a compatible metric.
Then it contains at most two non-trivial chains of gradient spheres.
\end{Proposition}

\begin{Remark}
It would be enough for our purposes to show that there {\em exists\/}
a compatible metric with at most two non-trivial chains.
However, we show that {\em every\/} compatible metric 
has at most two non-trivial chains.
\end{Remark}

\begin{pf*}{Proof of Proposition \ref{lem:two-branches}}
In each chain of gradient spheres,
the top weight is either $1$, $m$, or $m'$, and 
if $m$ or $m'$ is greater than one, it occurs exactly once
as the top weight in such a chain.
A similar situation holds at the minimum.
There are two cases, shown in Figure \ref{branches}:

\smallskip

\noindent {\em Case 1:}
Suppose that there exists a non-trivial chain of gradient spheres
with top weight equal to 1;
denote its bottom weight by $n_1$, and suppose that there exists
a different non-trivial chain with bottom weight equal to 1; denote its 
top weight
by $m_2$. Then the right hand side of \eqref{L=R} is at least
$\frac{d_1}{n_1} + \frac{d_2}{m_2}$.
Assume, without loss of generality, that $m \geq m' \geq 1$ and that
$n \geq n' \geq 1$.   Then this right hand side is at least
$\frac{1}{n} + \frac{1}{m}$, and it is greater than this number
if there exist additional chains of length greater than one.
Since the left hand side of \eqref{L=R}
is no larger than $\frac{1}{n} + \frac{1}{m}$, they both must be equal
to $\frac{1}{n} + \frac{1}{m}$.   This equality for the right hand side 
implies that $n=n_1$, $m=m_2$ and that there exist no additional 
chains of length greater than one.
The equality for the left hand side of implies that 
$m' = n' =1$, so there cannot be additional non-trivial chains of length one.
Therefore there are exactly two non-trivial chains of gradient spheres.

\smallskip

\noindent {\em Case 2:}
Now suppose that there are no two different non-trivial chains 
with the top isotropy weight being equal to $1$ in one chain
and the bottom isotropy weight being equal to $1$ in the other chain,
and suppose that there are three or more non-trivial chains.
Since at most two can have bottom weight different than $1$,
our assumption implies that there cannot be two chains
whose top isotropy weight is $1$. So there have to be exactly 
three non-trivial chains of gradient spheres, one with top weight $1$,
and two with top weights greater than $1$; in particular, $m$ and $m'$
are both greater than one.  A similar situation holds
at the bottom. Moreover, the $1$'s at the top and bottom must match: 
there must be a non-trivial
chain of gradient spheres in which both the top and bottom isotropy
weights are $1$. This chain contributes a positive integer to the
right hand side of \eqref{L=R}. 
Since the integers $m$, $m'$, $n$, and $n'$ are all greater than $1$,
the left hand side of \eqref{L=R} is less than $1$. This cannot happen.
\end{pf*}

\begin{figure}
$$
\begin{picture}(53,101)(0,-10)
\put(23,83){\blacken\ellipse{4}{4}}
\put(3,63){\blacken\ellipse{4}{4}}
\put(43,63){\blacken\ellipse{4}{4}}
\put(23,3){\blacken\ellipse{4}{4}}
\dashline{4.000}(3,63)(3,23)
\dashline{4.000}(43,63)(43,23)
\path(3,63)(23,83)(43,63)
\path(3,23)(23,3)(43,23)
\put(5,5){$\ss n_1$}
\put(0,73){$\ss 1$}
\put(35,73){$\ss m_2$}
\put(35,5){$\ss 1$}
\end{picture}
\quad
\begin{picture}(53,101)(0,-10)
\put(23,83){\blacken\ellipse{4}{4}}
\put(3,63){\blacken\ellipse{4}{4}}
\put(43,63){\blacken\ellipse{4}{4}}
\put(73,63){\blacken\ellipse{4}{4}}
\put(23,3){\blacken\ellipse{4}{4}}
\dashline{4.000}(3,63)(3,23)
\dashline{4.000}(43,63)(43,23)
\dashline{4.000}(73,63)(73,23)
\path(3,63)(23,83)(43,63)
\path(23,83)(73,63)
\path(3,23)(23,3)(43,23)
\path(23,3)(73,23)
\put(5,10){$\ss 1$}
\put(23,10){$\ss n$}
\put(53,10){$\ss n'$}
\put(8,73){$\ss 1$}
\put(21,73){$\ss m$}
\put(50,73){$\ss m'$}
\end{picture}
$$
\caption{Possible and impossible graphs}
\labell{branches}
\end{figure}
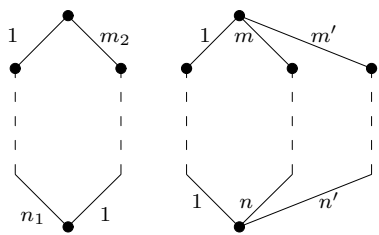

We now want to show that the graph of any compact four dimensional
Hamiltonian $S^1$-space with isolated fixed points comes from a Delzant
polygon. By Proposition \ref{lem:two-branches} it is enough to prove
the following, stronger, statement:

\begin{Proposition} \labell{exist-delzant}
Let $(M,\omega,\Phi)$ be a compact four dimensional Hamiltonian
$S^1$-space, equipped with a compatible metric.  Suppose that it contains
no more than two non-trivial chains of gradient spheres, and suppose
that any fixed surface, if exists, has genus zero.  Then there exists
a Delzant polygon which, by the recipe of section \ref{subsec:delzant},
gives the graph corresponding to this space.  
\end{Proposition}

We will prove an even stronger statement, that will be useful later.
Consider the graph corresponding to a compact four dimensional Hamiltonian
$S^1$-space.  Define an {\em extended graph\/} to be a graph obtained
from it by adding edges labeled $1$, in such a way that every interior
vertex is contained in exactly two edges, and in such a way that 
the moment map labels remain monotone along each branch of edges.
For instance, we can get an extended graph by choosing a compatible metric
(not necessarily generic) and marking a new edge for each non-trivial
free gradient sphere, and possibly for some trivial ones too.
(A trivial gradient sphere is a free gradient sphere that goes 
all the way from the bottom to the top.)
The branches in the graph
then correspond to the non-trivial chains of gradient spheres.

If we are given an extended graph with no more than two branches,
we can choose smooth invariant two-spheres in our manifold whose 
arrangement reflects the graph.  More specifically, 
for every new edge we can choose an invariant two-sphere connecting the
corresponding fixed points, such that the tangent spaces to the sphere
at its north and south poles are symplectic.
This guarantees that over the north and south poles, the circle acts
on the fiber of the normal bundle to the two-sphere with the weight
equal to the corresponding isotropy weight at that fixed point. 
It is easy to construct such a sphere by starting from a curve and 
"sweeping" it by the circle action.
The only possible problem is that if we have
an isolated extremum with both isotropy weights different than $\pm 1$, 
an invariant two-sphere that reaches this extremum will not be 
smooth unless it's a $Z_k$-sphere. (See \cite[Lemma 4.9]{ah-ha}.)
However, this problem does not come up if the extended graph has no more
than two branches.

We can now state our more general variant of \ref{exist-delzant}.

\begin{Proposition} \labell{exist-delzant'}
Let $(M,\omega,\Phi)$ be a compact four dimensional Hamiltonian
$S^1$-space, equipped with a compatible metric. Consider its graph.
Suppose that there exists an extended graph with exactly two branches.
Also suppose that any fixed surface, if exists, has genus zero.
Then there exists a Delzant polygon which, by the recipe of section
\ref{subsec:delzant}, gives the graph corresponding to this space.
\end{Proposition}

\begin{Remark}
If there exists an extended graph with no more than two branches,
then there exists one with exactly two branches: simply add 
trivial branches, each consisting of a single new edge going
all the way from the top to the bottom.  
\end{Remark}

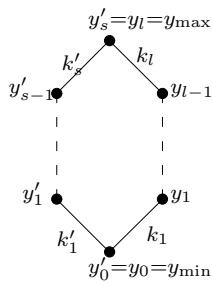
\begin{figure}
$$
\begin{picture}(53,101)(-12,-10)
\put(23,83){\blacken\ellipse{4}{4}}
\put(3,63){\blacken\ellipse{4}{4}}
\put(43,63){\blacken\ellipse{4}{4}}
\put(43,23){\blacken\ellipse{4}{4}}
\put(3,23){\blacken\ellipse{4}{4}}
\put(23,3){\blacken\ellipse{4}{4}}
\dashline{4.000}(3,63)(3,23)
\dashline{4.000}(43,63)(43,23)
\path(3,63)(23,83)(43,63)
\path(3,23)(23,3)(43,23)
\put(3,5){$\ss k'_1$}
\put(5,72){$\ss k'_s$}
\put(33,75){$\ss k_l$}
\put(37,7){$\ss k_1$}
\put(15,-5){$\ss y'_0 = y_0 = y_\min$}
\put(-10,23){$\ss y'_1$}
\put(46,23){$\ss y_1$}
\put(-15,63){$\ss y'_{s-1}$}
\put(46,63){$\ss y_{l-1}$}
\put(15,89){$\ss y'_s = y_l = y_\max$}
\end{picture}
$$
\caption{A graph with two branches}
\labell{two-branches}
\end{figure}

\begin{Remark}
Audin, in \cite[\S 3.3]{au:paper}, has a similar result
for the underlying $S^1$ manifold.
\end{Remark}

\begin{pf}
Take two chains of gradient spheres that correspond to the 
two branches of the extended graph.\footnote
{
We will call those spheres that correspond to the edges of the 
extended graph ``gradient spheres", even if the extended graph 
did not come from a compatible metric.
}
Let $k_1,\ldots,k_l$ be the isotropy weights along the first chain
and $k'_1,\ldots,k'_s$ the isotropy weights along the second chain.
Let  $y_0, y_1, \ldots, y_l$ be the values of the moment map at the
fixed points along the first chain, and $y'_0, y'_1, \ldots, y'_s$ 
those along the first chain. See Figure \ref{two-branches}.
Then $y_0 = y'_0 = y_\min$ and
$y_l = y'_s = y_\max$ are the extremal values of the moment map.

We will construct the edges of the polygon inductively, 
from the bottom up. The edges on the right will form
a polygonal path in $\R^2$ between a sequence of points
$p_0 , \ldots, p_l$, and those on the left will form a polygonal
path between a sequence of points $p'_0 , \ldots, p'_s$.
The $y$ coordinates of these points will be the corresponding 
moment map values. To determine these points,
it is enough to specify their $y$ coordinates, the slopes of
the edges $\ol{p_{i-1} p_i}$ and $\ol{p'_{i-1} p'_i}$, 
and the initial points, $p_0$ and $p'_0$.

Let us first specify the normal vectors to the edges.
If the minimum is isolated, let 
$u_1 = (k_1,b_1)$ and $u'_1 = (-k'_1,b'_1)$ be integer vectors 
such that $\det (u'_1 u_1) =1$; there exist such vectors because
$k_1$ and $k'_1$ are relatively prime.
If the minimum is a surface, let $u_0 = (0,-1)$,
and let $u_1 = (1,b_1)$ and $u'_1 = (-1,b'_1)$ be such that
$b_1 + b'_1 = -e_\min$, 
where $e_\min$ is the self intersection of the surface.

Let the next vectors be of the form $u_i = (k_i,b_i)$,
$i = 2, \ldots, l$, with $\det (u_{i-1} u_i) =1$.
As in the proof of Lemma \ref{mivi}, $u_i$ are integral.
Similarly, $u'_i = (-k'_i , b'_i)$, $i=2,\ldots,s$,
with $\det (u'_{i-1} u'_i) = -1$.  

Let $a_\min$ be the normalized symplectic area of the minimal
set of the moment map; if the minimum is isolated, this number is $0$.
Let $p_0 = (x_0,y_\min)$ and $p'_0 = (x'_0,y_\min)$
with $x'_0 = x_0 + a_\min$ and $x_0$ arbitrary.
Determine the points $p_1, \ldots, p_l$ 
by demanding the $y$ coordinate of $p_i$ to be the moment map 
value $y_i$ and the interval $\ol{p_{i-1} p_i}$ to be perpendicular
to $u_i$. Similarly, determine $p'_1, \ldots, p'_l$
by demanding the $y$ coordinate of $p'_i$ to be the moment map
value $y'_i$ and the interval $\ol{p'_{i-1} p'_i}$ to be perpendicular
to $u'_i$. 

We need to show that the polygonal paths
$\ol{ p_0 p_1 \ldots p_l}$ and $\ol{p'_0 p'_1 \ldots p'_s}$,
together with, possibly, horizontal edges 
$\ol{p'_0 p_0}$ and $\ol{p'_s p_l}$,
form a Delzant polygon whose corresponding graph coincides
with the graph of our manifold.
For this we need to check two things:
\begin{enumerate}
\item
the closed polygonal path
$\ol{p'_0 p_0 p_1 \ldots p_l p'_s p'_{s-1} \ldots p'_0}$
satisfies the conditions \eqref{delzant-condition} of Delzant, and
\item
if the maximum of the moment map is isolated,
the endpoints of the two paths coincide, and if the maximal
set is a surface, the interval $\ol{p'_s p_l}$ has length
equal to the normalized symplectic area of the surface;
\end{enumerate}

Let us parametrize the two paths by $(x(t),t)$ and $(x'(t),t)$
for $y_\min \leq t \leq y_\max$. 
The horizontal distance function, $f(t) = x(t) - x'(t)$, 
has the following properties:
\begin{enumerate}
\item[(i)]
The function $f(t)$ is continuous for $y_\min \leq t \leq y_\max$,
and is piecewise linear. Its slope, $\dd{t}f(t)$, only changes 
as $t$ crosses the points $y_i$ and $y'_j$. 
\item[(ii)]
At the bottom, $f(y_\min) = a_\min$.
If the moment map has an isolated minimum with isotropy weights
$n$ and $n'$, the slope of $f(t)$ right above the minimum
is $1 / nn'$. If the minimal set of the moment map is a surface
of self intersection $e_\min$, the slope of $f(t)$ right above
the minimum is $-e_\min$.
\item[(iii)]
As $t$ crosses the value $y_i$, the slope of $f(t)$
decreases by $\frac{1}{k_{i-1} k_i}$, and as $t$ crosses
the value $y'_j$, the slope of $f(t)$ decreases by
$\frac{1}{k'_{j-1} k'_j}$.
(Of course, if $y_i$ happens to be equal to $y'_j$, 
their contributions add up.)
\end{enumerate}

Conditions (i) and (ii) follow immediately from our construction.
Let us prove condition (iii). 
For $y_{i-1} < t < y_i$, $\frac{dx(t)}{dt} = -\frac{b_i}{k_i}$,
hence, as we cross the point $y_i$, the slope of $f(t)$
increases by $(-\frac{b_i}{k_i}) - (-\frac{b_{i-1}}{k_{i-1}})$
$= -\frac{k_{i-1} b_i - b_{i-1} k_i}{k_{i-1} k_i}$
$= -\frac{ \det(u_{i-1} u_i)  }{k_{i-1} k_i}$
$ = -\frac{1}{k_{i-1}k_i}$. 
Similarly, for $y'_{j-1} < t < y'_j$, 
$\frac{dx'(t)}{dt} = \frac{b'_j}{k'_j}$,
hence, as we cross the point $y'_j$, the slope of $f(t)$
increases by $(\frac{b'_j}{k'_j}) - (\frac{b'_{j-1}}{k'_{j-1}})$
$= -\frac{k'_j b'_{j-1} - b'_j k'_{j-1}}{k'_{j-1} k'_j}$
$= -\frac{ \det(u'_j u_{j-1})  }{k'_{j-1} k'_j}$
$ = -\frac{1}{k'_{j-1}k'_j}$. 

The density function for the \tDH measure, $\rho(t)$,
whose formula is given in \eqref{eq:rho}, satisfies these 
same Conditions (i)--(iii).
(To check this, notice that the $y_p$'s in \eqref{eq:rho}
are the same as our $y_i$'s and $y'_j$'s, and that for each $p$,
$(m_p,n_p)$ is either $(k_{i+1}, -k_i)$ or $(k'_{j+1}, -k'_j)$.)
Since conditions (i)--(iii) determine a function uniquely,
$x(t) - x'(t) = \rho(t)$. This implies part 2.

It remains to check the Delzant conditions.
Our construction of the normal vectors guarantees these
conditions are satisfied, except, perhaps, at the top. 
If the maximal set for the moment map is a surface, 
we have a top horizontal edge with an outward normal
$u_\max = (0,1)$, and the required conditions,
$\det (u_l u_\max) = \det (u_\max u'_s) =1$,
follow immediately from the fact that 
$u_l = (k_l,b_l)$ and $u'_s = (-k'_s,b'_s)$ with $k_l = k_s = 1$.
If the moment map has an isolated maximum, we have
$u_l = (m,b_l)$ and $u'_s = (-m',b'_s)$ 
where $-m,-m'$ are the isotropy weights at the maximum.
By \eqref{eq:rho}, the slope of the density function $\rho(t)$
for $t$ right below the maximum is $-\frac{1}{mm'}$.
For such $t$, the slope of $x(t)$ is $-\frac{b_l}{m}$
and the slope of $x'(t)$ is $\frac{b'_s}{m'}$.
Since $\rho(t) = x(t) - x'(t)$, we have
$-\frac{1}{mm'} = -\frac{b_l}{m} - \frac{b'_s}{m'} = 
-\frac{b_l m' + m b'_s}{mm'} = -\frac{\det(u_l u'_s)}{mm'}$,
which implies the Delzant condition $\det(u_l u'_s) =1$.
\end{pf}

We now reach the main result of this section:

\begin{Theorem} \labell{thm:toric}
Every four dimensional, compact Hamiltonian $S^1$-space
with isolated fixed points comes from a {\em K\"ahler toric variety}
by restricting the action to a sub-circle. 
\end{Theorem}

\begin{pf}
By Proposition \ref{lem:two-branches}, the space contains
at most two chains of gradient spheres (with respect to any 
compatible metric). By Proposition \ref{exist-delzant},
this implies that there exists a Delzant polygon which gives
rise to the graph of the space. Finally, by the Uniqueness
Theorem \ref{thm:uniqII}, the space is $S^1$-equivariantly
symplectomorphic to the K\"ahler toric variety that corresponds
to this Delzant polygon.
\end{pf}

Theorems \ref{thm:toric} and \ref{thm:uniqII} provide a 
{\em classification\/}
of compact four dimensional Hamiltonian $S^1$-spaces with isolated fixed
points: a complete list of spaces (up to isomorphism) is provided by the
list of Delzant polygons, and two Delzant polygons give isomorphic spaces
if and only if their graphs, as constructed in \S\ref{subsec:delzant},
coincide.

\begin{Corollary} \labell{cor:isol-torus}
In every four dimensional compact Hamiltonian $S^1$ space
with isolated fixed points, the circle action extends
to an effective Hamiltonian action of a two dimensional torus.
\end{Corollary}

\begin{Corollary} \labell{cor:isol-kahler}
Every four dimensional compact Hamiltonian $S^1$ space
with isolated fixed points admits a compatible K\"ahler structure.
\end{Corollary}

{
For completeness, let us state the classification 
more explicitly:

\begin{Theorem} \labell{thm:toric-classif}
The Hamiltonian $S^1$ spaces with isolated fixed points are classified 
by a set of polygons in $\R^2$, modulo an equivalence relation.
The set consists of those Delzant polygons 
(defined in \S\ref{subsec:delzant}) that have the following 
two additional properties.  
\begin{enumerate}
\item[(i)]
there are no horizontal edges;
\item[(ii)]
if the slope of an edge is $1/b$ with $b$ an integer,
the edge reaches either the top or bottom vertex.
\end{enumerate} 
The equivalence relation is the following.  Two polygons are equivalent
if they differ by an affine transformation of $\R^2$ 
of the form $(x,y) \mapsto (a + x + my, y)$ 
or of the form $(x,y) \mapsto (a - x + my, y)$ with $m \in \Z$.
\end{Theorem}

\begin{pf}
Consider the following maps:
$$ \begin{array}{ccccc}
\left\{ \parbox{1.13in}
 {\small Delzant polygons\\ that satisfy\\ (i) and (ii)}
\right\}
&
\stackrel{map_1}{\to}
&
\left\{ \parbox{1.50in}{\small Hamiltonian $S^1$-spaces\\ with\\ isolated
fixed points} \right\}
&
\stackrel{map_2}{\to}
&
\left\{ \parbox{0.45in}{\small labeled\\ graphs} \right\}
\end{array}$$
The maps $\map_2$ and $\map_2 \circ \map_1$ were described explicitly in
\S\ref{sec:graph}. 
The Uniqueness Theorem \ref{thm:uniqII} showed that $\map_2$ is one to one.
Proposition \ref{exist-delzant} showed that the image of $\map_2$ 
is equal to the image of $\map_2 \circ \map_1$. It follows that $\map_1$ 
is onto. 
The composition $\map_2 \circ \map_1$ is invariant under affine
transformations described above, hence so is $\map_1$. 
(These transformations correspond to automorphisms of the acting torus,
$S^1 \times S^1$, that preserve the second circle.)
Finally, by examining the proof of Proposition \ref{exist-delzant},
it follows that two Delzant polygons give the same graph {\em if and 
only if\/} they differ by an affine transformation of the above form.
(Indeed, the only choices involved in the construction of the Delzant
polygon were choosing which chain is on the left and which is on the right,
choosing the value of $x_0$, and choosing the normal vectors $u_1$ and
$u'_1$. Different choices amount to an affine transformation of the
polygon of the form described above.)
\end{pf}
}

Theorem \ref{thm:toric} enables us to determine when the circle
action extends to a Hamiltonian torus action:

\begin{Proposition} \labell{extend-action}
For a compact four dimensional Hamiltonian $S^1$-space,
the following are equivalent:
\begin{enumerate}
\item
each fixed surface has genus $0$, and
each non-extremal level set for the moment map contains at most 
two non-free orbits; 
\item
each fixed surface has genus $0$, and
there exists a metric for which there are no more than 
two non-trivial chains of gradient spheres,
\item
The circle action extends to a Hamiltonian action of a 2-torus.
\end{enumerate}
\end{Proposition}

\begin{pf}
It is easy to see that the second condition follows from the third
(take the K\"ahler metric of the toric variety) and implies the first.
It remains to show that the first condition implies the third.

The first condition implies that there exists an extended graph
with two or less branches that corresponds to the space.
(Construct the extended graph by connecting each interior
vertex from which no edge is coming up with the lowest vertex
above it from which no edge is coming down.)
The third then follows by Proposition \ref{exist-delzant'}.
\end{pf}

\begin{Remark}
Despite Part 2 of Proposition \ref{extend-action}, 
there could be three or more non-trivial
chains of gradient spheres for certain compatible metrics;
see Example \ref{ex:four} and Figure \ref{four}.  However,
this cannot happen if the fixed points are isolated, because of
Proposition \ref{lem:two-branches}.
\end{Remark}

\begin{Example} \labell{ex:four}
The toric variety shown in Figure \ref{four} 
provides a Hamiltonian $S^1$ space
in which there are four non-trivial chains of gradient spheres
with respect to a generic metric,
but in which the action extends to a torus action. 
\end{Example}

\begin{figure}
$$
\begin{picture}(250,85)(0,-10)
\put(220,65){\blacken\ellipse{60}{10}}
\put(220,5){\blacken\ellipse{60}{10}}
\put(210,25){\blacken\ellipse{6}{6}}
\put(230,45){\blacken\ellipse{6}{6}}
\put(195,45){\blacken\ellipse{6}{6}}
\put(245,25){\blacken\ellipse{6}{6}}
\path(195,65)(195,5)
\path(210,65)(210,5)
\path(230,65)(230,5)
\path(245,65)(245,5)
\path(20,65)(60,65)(80,45)(80,25)(60,5)(20,5)(0,25)(0,45)(20,65)
\put(0,65){\blacken\ellipse{2}{2}}
\put(20,65){\blacken\ellipse{2}{2}}
\put(40,65){\blacken\ellipse{2}{2}}
\put(60,65){\blacken\ellipse{2}{2}}
\put(80,65){\blacken\ellipse{2}{2}}
\put(0,45){\blacken\ellipse{2}{2}}
\put(20,45){\blacken\ellipse{2}{2}}
\put(40,45){\blacken\ellipse{2}{2}}
\put(60,45){\blacken\ellipse{2}{2}}
\put(80,45){\blacken\ellipse{2}{2}}
\put(0,25){\blacken\ellipse{2}{2}}
\put(20,25){\blacken\ellipse{2}{2}}
\put(40,25){\blacken\ellipse{2}{2}}
\put(60,25){\blacken\ellipse{2}{2}}
\put(80,25){\blacken\ellipse{2}{2}}
\put(0,5){\blacken\ellipse{2}{2}}
\put(20,5){\blacken\ellipse{2}{2}}
\put(40,5){\blacken\ellipse{2}{2}}
\put(60,5){\blacken\ellipse{2}{2}}
\put(80,5){\blacken\ellipse{2}{2}}
\end{picture}
$$
\caption{A toric variety with four non-trivial chains of
gradient spheres for a generic metric}
\labell{four}
\end{figure}
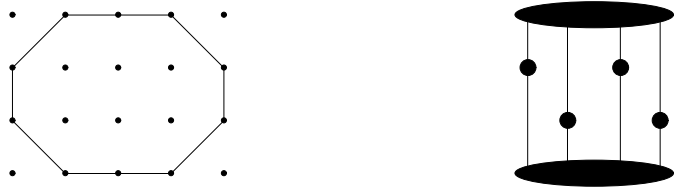

\begin{Remark}
Proposition \ref{extend-action} is very close 
to Audin's claim, that a circle action extends to a torus 
action if and only if each level set for the moment map 
meets at most two gradient spheres \cite[Thm.3.5.1]{au:paper}.
Audin must have had in mind non-free gradient spheres;
otherwise, Example \ref{ex:four} provides a counterexample.
Alternatively, she might have implicitly assumed that the
metric would be preserved by the action of the whole torus.
This assumption is also necessary for the second part 
of \cite[Thm.3.5.1]{au:paper}, which states that 
the torus action is determined by the circle action;
otherwise, Example \ref{hirz} provides a counterexample.
\end{Remark}

\section{Blowing-up}
\labell{sec:blowup}
In sections \ref{sec:graph} and \ref{sec:uniqueness} we showed
how to distinguish between compact four dimensional Hamiltonian
$S^1$-spaces. To complete the classification we need to determine
which such spaces occur.
In this section we show that every such space can be obtained 
from a {\em minimal space\/} by a sequence of equivariant symplectic 
blowups at fixed points. 
Moreover, the result of such a blowup is completely determined 
by its ``size", and in section \ref{sec:existence} we give 
the precise restrictions on these sizes.

This section essentially reproduces the results of Ahara, Hattori, 
and Audin.  We took the liberty to prove them from scratch, because 
our proof is shorter,
and because we provide a more precise result, which we need for
section \ref{sec:existence}.

\subsection{Equivariant symplectic blow-ups and blow-downs}
\labell{subsec:blowdown}

Symplectic blowups and blowdowns are explained carefully
in \cite{md:blowups}. Let us now recall what they are,
while also keeping track of a circle (or torus) action.

Take $\C^2$ with its standard symplectic structure, remove a ball of
radius $r$, and collapse the boundary along the Hopf fibration.
The resulting topological space is the union of the manifolds
$\{ |z|^2+ |w|^2 > r^2 \}$ and $\{ |z|^2 + |w|^2 = r^2 \} / S^1$,
each of which inherits a symplectic structure from $\C^2$.
To obtain a smooth and symplectic structure on this space, 
one can realize it as the 
symplectic quotient of $\C^2 \times \C$ with respect to the 
circle action $(az,aw,a^{-1}u)$ at the 
level set ${1 \over 2} (|z|^2 + |w|^2 - |u|^2) = {1 \over 2} r^2$;
for more details see \cite{g-s:birational}. 
This is a special case of Lerman's {\em symplectic cutting} 
\cite{l:cutting}.

The {\em size\/}, or {\em amount\/}, of this blowup is  
$\lambda = r^2/2$.
The submanifold $\{ |z|^2 + |w|^2 = r^2 \} / S^1$ is the 
{\em exceptional divisor;} its area is $2\pi \lambda$.
Any two integers $m,n$ determine a Hamiltonian action of the circle group 
on $\C^2$, given by $(\lambda^m z, \lambda^{-n} w)$;
this circle action and its moment map descend to the blow up.
The exceptional divisor is a $Z_{m+n}$-sphere; the circle
action on it can be written in homogeneous coordinates as
$[\lambda^m z, \lambda^{-n} w] = [\lambda^{m+n} z,w]$.

\begin{figure} 
$$
\begin{picture}(60,75)(0,-10)
\path(0,60)(0,20)(20,0)(60,0)
\dashline{4.000}(0,20)(0,0)(20,0)
\end{picture}
$$
\caption{Blowing up ${\protect \Bbb C}^2$}
\labell{blowupC2}
\end{figure}
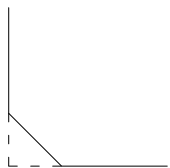

As a toric variety, the blow-up of $\C^2$ 
corresponds to the positive quadrant minus a triangle; 
see Figure \ref{blowupC2}.
The positive quadrant is the image of $\C^2$ under the moment map
$(z,w) \mapsto ({1 \over 2} |z|^2, {1 \over 2} |w|^2)$,
and the triangle is the image of the ball $\{ |z|^2 + |w|^2 < r^2 \}$.
Transforming this picture by an appropriate element of  
SL$(2,\Z)$, we can get a picture that corresponds
to an action of $S^1 \times S^1$ in which the second $S^1$ acts by the
integers $(m,n)$ as described earlier. Applying the recipe of 
\S\ref{subsec:delzant}, we get the graph corresponding to this 
$S^1$-space.
Figures \ref{blowup-A}--\ref{blowup-D}, together with
Figures \ref{blowup-B}--\ref{blowup-D} turned up side down, 
give all the possible graphs of $\C^2$ and its blow up
(with the addition of some edges labeled $1$ which correspond 
to free gradient spheres).
In these figures, $\alpha$ is a moment map label,
$\lambda = r^2/2$ is the amount by which we blow up,
$v$ and $v-\lambda$ are area-labels, and $g$ is a genus.

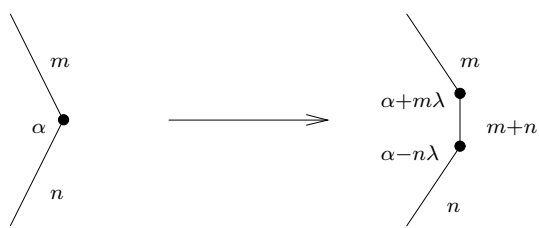
\begin{figure}
$$
\begin{picture}(204,99)(0,-10)
\put(22,42){\blacken\ellipse{4}{4}}
\put(172,52){\blacken\ellipse{4}{4}}
\put(172,32){\blacken\ellipse{4}{4}}
\path(2,82)(22,42)(2,2)
\path(62,42)(122,42)
\path(114.000,40.000)(122.000,42.000)(114.000,44.000)
\path(152,82)(172,52)(172,32)(152,2)
\put(182,37){$\ss m+n$}
\put(10,37){$\ss \alpha$}
\put(142,47){$\ss \alpha + m\lambda$}
\put(142,27){$\ss \alpha - n\lambda$}
\put(17,62){$\ss m$}
\put(17,12){$\ss n$}
\put(167,7){$\ss n$}
\put(172,62){$\ss m$}
\end{picture}
$$
\caption{Blowing up at an interior fixed point}
\labell{blowup-A}
\end{figure}

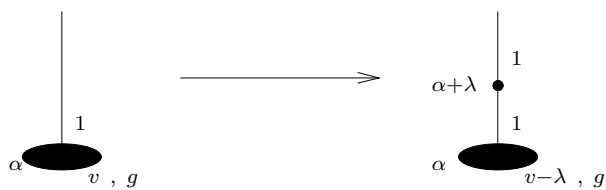
\begin{figure}
$$
\begin{picture}(219,90)(0,-10)
\put(195,15){\blacken\ellipse{30}{10}}
\put(30,15){\blacken\ellipse{30}{10}}
\put(195,42){\blacken\ellipse{4}{4}}
\path(75,45)(150,45)
\path(142.000,43.000)(150.000,45.000)(142.000,47.000)
\path(195,70)(195,15)
\path(30,70)(30,15)
\put(200,50){$\ss 1$}
\put(200,25){$\ss 1$}
\put(35,25){$\ss 1$}
\put(10,10){$\ss \alpha$}
\put(40,5){$\ss v \ , \ g$}
\put(205,5){$\ss v - \lambda \ , \ g$}
\put(170,10){$\ss \alpha$}
\put(170,40){$\ss \alpha + \lambda$}
\end{picture}
$$
\caption{Blowing up at a point on $B_{\protect \operatorname{min}}$.}
\labell{blowup-B}
\end{figure}

\begin{figure} 
$$
\setlength{\unitlength}{0.0095in}
\begin{picture}(239,97)(0,-10)
\put(32,2){\blacken\ellipse{4}{4}}
\put(37,-3){$\ss \alpha$}
\put(197,2){\blacken\ellipse{4}{4}}
\put(202,-3){$\ss \alpha + n\lambda$}
\put(167,32){\blacken\ellipse{4}{4}}
\put(170,33){$\ss \alpha +m\lambda $}
\path(2,62)(32,2)(72,42)
\dottedline{5}(2,82)(2,62)
\dottedline{5}(72,82)(72,42)
\path(167,62)(167,32)(197,2)(237,42)
\dottedline{5}(167,82)(167,62)
\dottedline{5}(237,82)(237,42)
\path(92,27)(147,27)
\path(139.000,25.000)(147.000,27.000)(139.000,29.000)
\put(5,22){$\ss m$}
\put(52,12){$\ss n$}
\put(217,12){$\ss n$}
\put(150,47){$\ss m$}
\put(150,12){$\ss m-n$}
\end{picture}
$$
\caption{Blowing up at a minimum with isotropy weights $n<m$.}
\labell{blowup-C}
\end{figure}
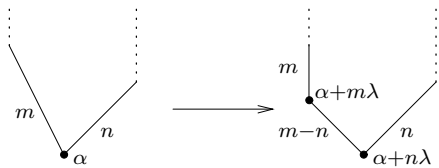

\begin{figure}
$$
\begin{picture}(163,62)(0,-10)
\put(142,5){\blacken\ellipse{30}{10}}
\put(105,2){$\ss \alpha+\lambda$}
\put(160,0){$\ss \ v=\lambda \ , \ g=0$}
\put(22,5){\blacken\ellipse{4}{4}}
\put(12,2){$\ss \alpha$}
\path(2,45)(22,5)(42,45)
\path(132,45)(132,10)
\path(152,45)(152,10)
\path(57,25)(107,25)
\path(99.000,23.000)(107.000,25.000)(99.000,27.000)
\put(157,20){$\ss 1$}
\put(122,20){$\ss 1$}
\put(37,15){$\ss 1$}
\put(5,15){$\ss 1$}
\end{picture}
$$
\caption{Blowing up at a minimum with isotropy weights $1,1$.}
\labell{blowup-D}
\end{figure}
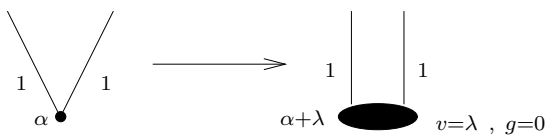

We can now perform the same construction on a manifold:
an {\em equivariant Darboux chart\/} on a symplectic manifold
with a group action is an equivariant symplectomorphism from
an open subset of the manifold onto an open subset of $\C^2$
with a linear action. If its image contains a closed ball of radius
$r$, we can use this chart to define an equivariant symplectic
blow-up of the manifold by the amount $\lambda = r^2 /2$.
By the local normal form (namely, Corollary \ref{stab=S1}),
there exists an equivariant Darboux chart on the neighborhood
of any fixed point.  However, in general it is difficult to determine
how large it can be.

We now show that the result of an equivariant symplectic
blow-up on our spaces is independent of the choice of a Darboux chart. 
This is in contrast to the non-equivariant case, when the uniqueness 
of the symplectic blow-up is a difficult issue (see \cite{md:uniqueness}).

\begin{Proposition} \labell{same}
Let $M$ be a compact four dimensional Hamiltonian $S^1$-space
and $p\in M$ a fixed point. Then any two equivariant
symplectic blow-ups of $M$ at $p$ by the same amount $\lambda$
yield isomorphic Hamiltonian $S^1$ spaces. 
\end{Proposition}

\begin{pf}
The graphs of the blow-ups are the same;
they are given by Figures \ref{blowup-A}--\ref{blowup-D}.
By the Uniqueness Theorem \ref{thm:uniqII},
the blown-up spaces are the same.
\end{pf}

{

We now use blow-ups to construct 
a compact four dimensional Hamiltonian $S^1$-space 
with an isolated minimum and with three non-trivial chains
of gradient spheres:

\begin{Example}
\labell{ex:X}
Consider the Delzant polygon in Figure \ref{exX} on the left,
whose edge vectors, starting from the low left vertex 
and proceeding counterclockwise, are
$(1,0)$, $(1,1)$, $(7,14)$, $(-8,-12)$, $(-1,-2)$, and $(0,-1)$.
Take the corresponding toric variety, and perform on it an 
$S^1$-equivariant blow-up at a point on the minimal surface 
of the moment map.
The resulting space has three non-trivial chains of gradient spheres.
Its graph is drawn in Figure \ref{exX} on the right.
\end{Example}

\begin{Remark}
Example \ref{ex:X} provides a counterexample to Proposition 3.1.2
of \cite{au:paper}, where it was implicitly assumed that the 
gradient spheres all have non-trivial stabilizers.
\end{Remark}

Example \ref{ex:X} also provides an example of a compact four dimensional
Hamiltonian $S^1$-space that contains a non-smooth non-trivial
gradient sphere:
the free gradient sphere that reaches the maximum is not smooth
at its north-pole. This follows, as is shown in \cite[Lemma 4.9]{ah-ha},
from the fact that both isotropy weights at the maximum have an absolute
value greater than $1$.

\begin{Remark}
In Example \ref{ex:X}, if the blow-up is by an amount $\lambda$
different than $1$, the resulting space admits {\em some\/} compatible
metric (necessarily non-generic) for which there are
only {\em two\/} chains of gradient spheres,
and the space is a toric variety;
this follows from Proposition \ref{extend-action}.
This is no longer true if $\lambda =1$;
see Example \ref{ex:X2}.
\end{Remark}

\begin{figure} 
$$
\setlength{\unitlength}{0.0050in}
\begin{picture}(184,319)(0,-10)
\path(0,0)(20,0)(40,20)(180,300)(20,60)(0,20)(0,0)
\put(180,300){\blacken\ellipse{10}{10}}
\put(40,20){\blacken\ellipse{10}{10}}
\put(20,60){\blacken\ellipse{10}{10}}
\put(0,20){\blacken\ellipse{10}{10}}
\put(0,0){\blacken\ellipse{10}{10}}
\put(20,0){\blacken\ellipse{10}{10}}
\end{picture}
\quad \quad \quad
\setlength{\unitlength}{0.0095in}
\begin{picture}(77,162)(0,-10)
\put(45,145){\blacken\ellipse{4}{4}}
\put(15,85){\blacken\ellipse{4}{4}}
\put(75,45){\blacken\ellipse{4}{4}}
\put(15,45){\blacken\ellipse{4}{4}}
\put(45,45){\blacken\ellipse{4}{4}}
\put(45,5){\blacken\ellipse{30}{10}}
\path(35,5)(55,5)(75,45)(45,145)(15,85)(15,45)(35,5)
\path(45,5)(45,145)
\put(65,95){$\ss 2$}
\put(70,15){$\ss 1$}
\put(15,110){$\ss 3$}
\put(0,55){$\ss 2$}
\put(15,15){$\ss 1$}
\put(35,65){$\ss 1$}
\put(35,15){$\ss 1$}
\end{picture}
$$
\caption{Blowing up at the minimum for Example \protect\ref{ex:X}}
\labell{exX}
\end{figure}
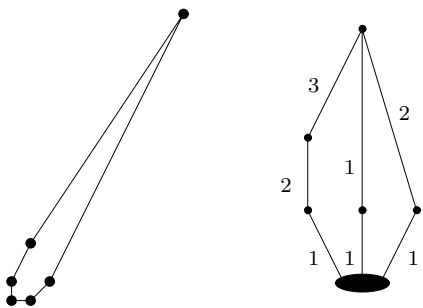
}

{\em Symplectic blowing down} is the reverse procedure to symplectic
blowing up: let $(M,\omega,\Phi)$ be a four dimensional Hamiltonian $S^1$
space, and let $C \subseteq M$ be an invariant symplectic sphere with
self intersection $-1$ and total area $\pi r^2$.  By the equivariant
tubular neighborhood theorem (specifically, Corollaries \ref{normal-near-B}
and \ref{nbhd-Zk}), a neighborhood of $C$ in $M$ is equivariantly
symplectomorphic to a neighborhood of the exceptional divisor in the
$\lambda$-blowup of $\C^2$, where $\lambda = r^2/2$, with the action being
given by some integers $(m,n)$.  To blow down, we remove a neighborhood
of $C$ and glue in a standard ball.

{
\subsection{Blowing down to a minimal space}
\labell{subsec:down-to-minimal}
In this section we show that every compact four dimensional 
Hamiltonian $S^1$ space is obtained by a sequence of 
equivariant symplectic blow-ups from a minimal space.
The minimal space is either $\CP^2$, or a Hirzebruch
surface, or a space with two fixed surfaces and no
interior fixed points, which we later identify as a ruled
manifold.

\begin{Theorem} \labell{thm:neat}
Every compact four dimensional Hamiltonian $S^1$-space with two fixed 
surfaces can be obtained from a space with no interior fixed points
by a sequence of equivariant symplectic blow-ups
at fixed points that are not minima for the moment map.
\end{Theorem}

\begin{pf}
It is enough to show that any non-trivial 
chain of gradient spheres contains a gradient sphere with self 
intersection $-1$, which is not the bottom sphere; 
such a sphere can be blown down, and the theorem will follow by induction
on the number of interior fixed points.

Consider a non-trivial chain of gradient spheres.
Denote their isotropy weights, from the bottom up, 
by $k_1, \ldots, k_l$. 
If all the $k_i$'s are equal to $1$,
the top sphere has self intersection $-1$ 
(by Lemma \ref{mg} with $m=0$, $n=-1$, and $k=1$),
and we are done.

If not all the $k_i$'s are equal to $1$, let $k_j$ be 
the largest weight in the chain. Then $k_j>1$.
Since $k_1 = k_l = 1$, the $j$'th sphere is not the first
nor last in its chain.  Since $k_j$ is maximal,
\begin{equation} \labell{geq}
k_j \geq k_{j-1} \quad \mbox{and} \quad k_j \geq k_{j+1}.
\end{equation}
Since any two consecutive weights in the chain are relatively prime
(Lemma \ref{lem:conditions}),
the inequalities in \eqref{geq} must be strict; adding them up
and dividing by $k_j$, we get
\begin{equation} \labell{gt}
\frac{k_{j-1} + k_{j+1}}{k_j} < 2.
\end{equation}
The left term in \eqref{gt}, being a positive integer
(Lemma \ref{lem:conditions}) and less than $2$, must be equal to $1$.
The self intersection of the $j$th sphere is the negative of this
integer (Lemma \ref{lem:conditions}), so it is equal to $-1$.
\end{pf}
}


Next, we deal with spaces with one fixed surface.

\begin{Theorem} \labell{thm:one}
Any compact four dimensional Hamiltonian $S^1$-space with one fixed
surface can be obtained from a K\"ahler toric variety by a sequence of
equivariant symplectic blow-ups.
\end{Theorem}

\begin{pf}
Let us assume that the fixed surface is the minimal set of the moment
map; a maximum can be treated similarly.  Fix a compatible metric.
Figure \ref{one-surface} then shows the arrangement of gradient spheres
and their isotropy weights. In a chain of gradient spheres, if all
isotropy weights are equal to one, we can blow down the bottom sphere.
Otherwise, if some isotropy weight is greater than both its neighbors,
then, arguing as in the proof of Theorem \ref{thm:neat}, we can blow
down that gradient sphere. Therefore, by a sequence of equivariant
symplectic blow-downs,
we can arrive at a space in which the isotropy weights
along each chain of gradient spheres are strictly increasing.  This graph
must have one of the three forms in Figure \ref{three}; in particular, it
can have at most two non-trivial chains of gradient spheres. By Proposition
\ref{exist-delzant}, the space is a toric variety.
\end{pf}

\begin{figure}
$$
\setlength{\unitlength}{0.0090in}
\begin{picture}(147,162)(0,-10)
\put(45,145){\blacken\ellipse{4}{4}}
\put(5,105){\blacken\ellipse{4}{4}}
\put(45,105){\blacken\ellipse{4}{4}}
\put(85,105){\blacken\ellipse{4}{4}}
\put(85,55){\blacken\ellipse{4}{4}}
\put(45,55){\blacken\ellipse{4}{4}}
\put(5,55){\blacken\ellipse{4}{4}}
\put(45,5){\blacken\ellipse{40}{10}}
\put(145,105){\blacken\ellipse{4}{4}}
\put(145,55){\blacken\ellipse{4}{4}}
\put(105,80){\blacken\ellipse{2}{2}}
\put(115,80){\blacken\ellipse{2}{2}}
\put(125,80){\blacken\ellipse{2}{2}}
\path(5,55)(30,5)
\path(45,55)(45,5)
\path(85,55)(55,5)
\path(45,145)(5,105)
\path(45,145)(45,105)
\path(45,145)(85,105)
\dottedline{5}(5,105)(5,55)
\dottedline{5}(45,105)(45,55)
\dottedline{5}(85,105)(85,55)
\dottedline{5}(145,105)(145,55)
\path(45,145)(145,105)
\path(145,55)(65,5)
\put(30,25){$\ss 1$}
\put(5,25){$\ss 1$}
\put(60,25){$\ss 1$}
\put(115,25){$\ss 1$}
\put(0,115){$\ss m$}
\put(35,115){$\ss n$}
\put(77,115){$\ss 1$}
\put(130,115){$\ss 1$}
\end{picture}
$$
\caption{A Hamiltonian $S^1$-space with one fixed surface}
\labell{one-surface}
\end{figure}
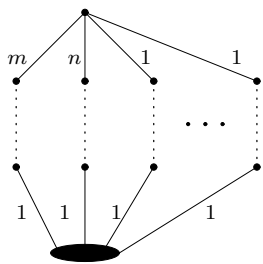

\begin{figure}
$$
\setlength{\unitlength}{0.0100in}
\begin{picture}(350,134)(0,-10)
\put(175,85){\blacken\ellipse{4}{4}}
\put(175,35){\blacken\ellipse{4}{4}}
\put(190,5){\blacken\ellipse{40}{10}}
\put(190,110){\blacken\ellipse{4}{4}}
\put(330,5){\blacken\ellipse{40}{10}}
\put(330,110){\blacken\ellipse{4}{4}}
\put(35,85){\blacken\ellipse{4}{4}}
\put(35,35){\blacken\ellipse{4}{4}}
\put(50,5){\blacken\ellipse{40}{10}}
\put(65,85){\blacken\ellipse{4}{4}}
\put(65,35){\blacken\ellipse{4}{4}}
\put(50,110){\blacken\ellipse{4}{4}}
\dottedline{5}(175,85)(175,35)
\path(190,110)(175,85)
\dottedline{5}(35,85)(35,35)
\dottedline{5}(65,85)(65,35)
\path(50,110)(35,85)
\path(50,110)(65,85)
\put(170,95){$\ss m$}
\put(30,95){$\ss m$}
\put(65,95){$\ss n$}
\end{picture}
$$
\caption{These spaces cannot be blown down along a gradient sphere}
\labell{three}
\end{figure}
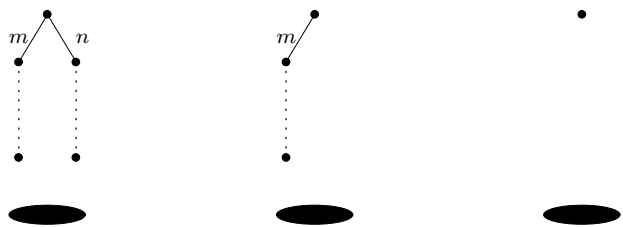

Equivariant symplectic blowing up and down, 
with the exclusion of Proposition \ref{same},
can be performed on symplectic manifolds of arbitrary dimension
with symplectic actions of arbitrary compact Lie groups.
Blowing up a four dimensional symplectic toric variety
corresponds to ``chopping off a corner" of the Delzant
polygon (Figure \ref{blowupC2}); blowing down corresponds
to ``filling in" a corner (Figure \ref{fig:corner}).

\begin{Example}[Blow ups at different minimal points]
\labell{ex:blow-min}
If we blow up a compact four dimensional Hamiltonian $S^1$-space 
at different minimal points for the moment map,
we get isomorphic spaces; the graphs all change in the same way
(Figure \ref{blowup-B}).
This is no longer true if we keep track of more structure, for instance, 
a toric action, or a complex structure:

On the left of Figure \ref{blow-min} we see $S^2 \times S^2$ 
blown up at a point as a toric variety
(a square with a corner chopped off)
and as a Hamiltonian $S^1$-space (a graph).
On the right we see two non-isomorphic toric varieties,
obtained from it by blowing up at two different points
on the bottom surface.
They remain non-isomorphic
even if we just keep track of the complex structure and the circle 
action, and not the torus action,
because their $S^1$-invariant complex curves are arranged differently. 
However, as Hamiltonian $S^1$-spaces they are isomorphic; 
their graph is drawn further on the right. 
Note that 
the K\"ahler metric on the first of the two toric varieties on the right 
is not generic in the sense of Corollary \ref{generic}.
\end{Example}

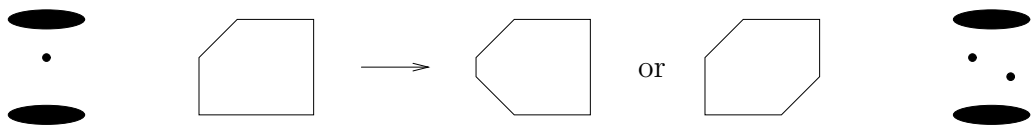
\begin{figure}
$$
\setlength{\unitlength}{0.0100in}
\begin{picture}(535,75)(0,-10)
\path(160,55)(120,55)(100,35)(100,5)(160,5)(160,55)
\path(305,55)(265,55)(245,35)(245,25)(265,5)(305,5)(305,55)
\path(425,55)(385,55)(365,35)(365,5)(405,5)(425,25)(425,55)
\put(515,5){\blacken\ellipse{40}{10}}
\put(515,55){\blacken\ellipse{40}{10}}
\put(505,35){\blacken\ellipse{4}{4}}
\put(525,25){\blacken\ellipse{4}{4}}
\put(20,5){\blacken\ellipse{40}{10}}
\put(20,35){\blacken\ellipse{4}{4}}
\put(20,55){\blacken\ellipse{40}{10}}
\path(185,30)(220,30)
\path(212.000,28.000)(220.000,30.000)(212.000,32.000)
\put(330,25){\mbox{or}}
\end{picture}
$$
\caption{Two blow-ups, for Example \ref{ex:blow-min}} 
\labell{blow-min}
\end{figure}

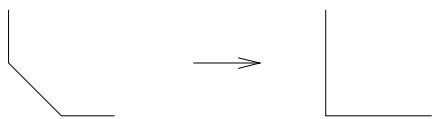
\begin{figure} 
$$
\begin{picture}(160,55)(0,-10)
\path(70,20)(95,20)
\path(87.000,18.000)(95.000,20.000)(87.000,22.000)
\path(0,40)(0,20)(20,0)(40,0)
\path(120,40)(120,0)(160,0)
\end{picture}
$$
\caption{Blowing down a toric variety}
\labell{fig:corner}
\end{figure}

A four dimensional symplectic toric variety can be blown down along
a two-sphere corresponding to any edge whose inward (or outward)
primitive normal 
vector is the sum of the inward primitive normal vectors of the 
neighboring edges. (In Figure \ref{fig:corner}, these normal edges 
are $(1,1) = (0,1) + (1,0)$; every other blow down
is obtained by transforming this picture by an element of 
SL$(2,\Z)$.)
In particular, in order to tell whether we can blow down,
we don't need to know the Delzant polygon;
we only need to know its primitive inward normal vectors.
These form a complete and smooth fan.
Blowing down amounts to removing a vector that was equal
to the sum of its neighbors.
A sequence of exercises in \cite[\S2.5]{fulton} shows that
if a fan cannot be blown down, it must be one of the two fans 
in Figure \ref{min-fans}, or their image by an element of SL$(2,\Z)$.
The spaces corresponding to these fans are $\CP^2$ and Hirzebruch 
surfaces.   This proves:

\begin{Lemma} \labell{thm:fulton}
Every compact four dimensional symplectic toric variety
with one or no fixed surfaces
can be obtained from either $\CP^2$ or a Hirzebruch surface
by a sequence of equivariant symplectic blow-ups.
\end{Lemma}

\begin{Remark}
If we just consider the underlying topological spaces,
this theorem was also proved in \cite{yavin}.
\end{Remark}

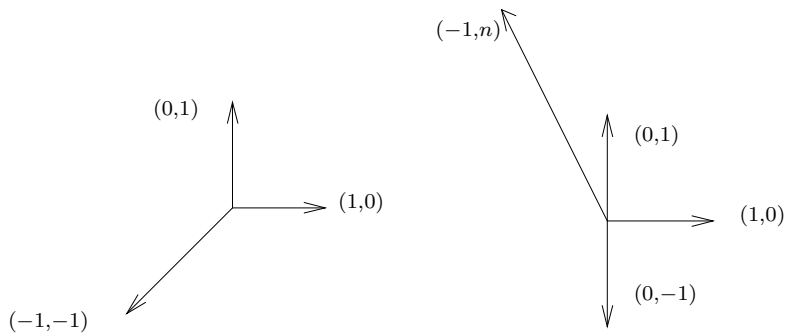
\begin{figure}
$$
\begin{picture}(150,127)(0,-10)
\path(87.000,77.000)(85.000,85.000)(83.000,77.000)
\path(85,85)(85,45)
\path(85,45)(120,45)
\path(112.000,43.000)(120.000,45.000)(112.000,47.000)
\path(85,45)(45,5)
\path(49.243,12.071)(45.000,5.000)(52.071,9.243)
\put(55,80){$\ss (0,1)$}
\put(0,0){$\ss (-1,-1)$}
\put(125,45){$\ss(1,0)$}
\end{picture}
\quad
\begin{picture}(178,96)(0,-10)
%
%
\path(65,40)(65,80)
\path(67.000,72.000)(65.000,80.000)(63.000,72.000)
\put(75,70){$\ss (0,1)$}
\path(65,40)(65,0)
\path(67,8)(65,0)(63,8)
\put(75,10){$\ss (0,-1)$}
\path(65,40)(105,40)
\path(97.000,38.000)(105.000,40.000)(97.000,42.000)
\put(115,40){$\ss (1,0)$}
\path(65,40)(25,120)
\path(26.789,111.950)(25,120)(30.367,117.739)
\put(0,110){$\ss (-1,n)$}
%
%
\end{picture}
$$
\caption{Minimal fans}
\labell{min-fans}
\end{figure}

\begin{Theorem} \labell{big}
Every compact four dimensional Hamiltonian $S^1$ space
can be obtained by a sequence of $S^1$-equivariant symplectic
blow-ups from 
\begin{enumerate}
\item
a space with two fixed surfaces and no interior fixed points, or 
\item
$\CP^2$ or a Hirzebruch surface, with a symplectic form
and a circle action that come from a K\"ahler form and a toric action,
and in which the K\"ahler metric is generic in the sense of 
Corollary \ref{generic}.
\end{enumerate}
\end{Theorem}

\begin{pf}
For a space with two fixed surfaces, this was proved
in Theorem \ref{thm:neat}. If the space has one or no
fixed surfaces, by Theorems \ref{thm:toric} and \ref{thm:one} 
it is obtained from a K\"ahler toric variety by a 
(possibly empty) sequence of equivariant symplectic blow-ups.
This toric variety, constructed in the proof of 
Theorems \ref{thm:toric} or \ref{thm:one}, has a generic 
K\"ahler metrics in the sense of Corollary \ref{generic}. 
By Lemma \ref{thm:fulton}, this K\"ahler toric variety
is itself obtained by a sequence of equivariant symplectic
blow-ups from $\CP^2$ or from a Hirzebruch surface, again with a
generic metric.
\end{pf}

\begin{Remark}
Theorem \ref{big} is a refinement of results of Ahara, Hattori,
and Audin.
\end{Remark}

{
\begin{Example}
[Different blow-ups lead to isomorphic toric varieties]
\labell{ex:different}

Figure \ref{different} shows the polygons of three toric 
varieties. The two on the right are ruled manifolds that are 
non-isomorphic as $S^1$-spaces. Both of them, when blown up once, 
yield the same space, on the left.  
(One can obtain the polygon on the left from either one of the
polygons on the right by chopping the corner at the marked vertex.) 
\end{Example}

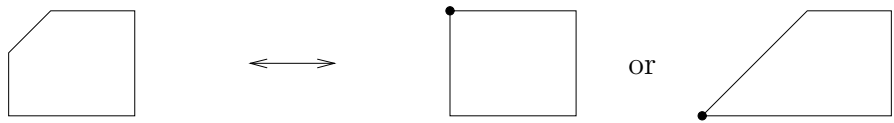
\begin{figure}
$$
\setlength{\unitlength}{0.0110in}
\begin{picture}(420,65)(0,-10)
\path(123,27)(115,25)(123,23)
\path(115,25)(155,25)
\path(147,23)(155,25)(147,27)
\path(210,50)(270,50)(270,0)(210,0)(210,50)
\put(210,50){\circle*{4}}
\path(330,0)(380,50)(420,50)(420,0)(330,0)
\put(330,0){\circle*{4}}
\path(60,50)(20,50)(0,30)(0,0)(60,0)(60,50)
\put(295,20){\mbox{or}}
\end{picture}
$$
\caption{Different blow-ups can yield the same space}
\labell{different}
\end{figure}
}

\subsection{Minimal spaces}
\labell{subsec:min-spaces}
In this section we describe the minimal spaces which occur
in Theorem \ref{big}. In particular, we will see that they are
all K\"ahler.

The first minimal spaces are $\CP^2$, which we all know and love, 
with multiples of the Fubini-Study symplectic form, and with $S^1$-actions 
obtained as inclusions $S^1 \hookrightarrow T^2$ followed by the toric
actions $(a,b) : [w_0,w_1,w_2] \mapsto [w_0,a w_1,b w_2]$.
The corresponding Delzant polygons are obtained from the one on the left of 
Figure \ref{CP2-Hirz} by applying transformations in SL$(2,\Z)$.
The graphs for all these $S^1$ spaces are shown in Figures
\ref{CP2-A}, \ref{CP2-B}, and \ref{CP2-B} turned up-side-down.
(The figures may contain some additional edges, labeled $1$,
which correspond to free gradient spheres, and which do not
occur in the graph.)
In figure \ref{CP2-A}, $m$ and $n$ are relatively prime positive
integers, and $\alpha$ and $\beta$ are positive real numbers.
In Figure \ref{CP2-B}, $\alpha$ and $\alpha + \lambda$
are moment map labels, and $\lambda$ is the area label.

The next minimal spaces are Hirzebruch surfaces. To each integer $n$
there corresponds the Hirzebruch surface
$$ \{ ([z_0,z_1,z_2],[w_1,w_2]) \in \CP^2 \times \CP^1 \ | \ 
      z_1w_2^n = z_2 w_1^n \},$$
with the toric action
$$(a,b) : ([z_0,z_1,z_2],[w_1,w_2]) \mapsto
([a z_0, z_1, b^n z_2],[w_1, b w_2]),$$
with symplectic forms induced by multiples
of the Fubini-Study forms on $\CP^2$ and on $\CP^1$. The corresponding
Delzant polygons are obtained from those on the right of 
Figure \ref{CP2-Hirz} by applying transformations in SL$(2,\Z)$.
The graphs for the corresponding $S^1$-spaces are those shown in 
Figure \ref{Hirz-graphs} and those obtained by turning them up-side-down.
In the third of these graphs, $\alpha$ is a moment map label, $s$
an area label, and $g$ a genus label.
(Again, Figure \ref{Hirz-graphs} may contain edges labeled $1$,
which do not occur in the graph, and which correspond
to free gradient spheres in the space.)

\begin{Remark} 
Of the Hirzebruch surfaces with $S^1$ actions,
the only ones whose K\"ahler metrics are not generic
(in the sense of Corollary \ref{generic})
are those in the middle of Figure \ref{Hirz-graphs}, when $c=1$. 
However, each of these, as a Hamiltonian $S^1$-space is isomorphic
to another toric varieties, whose K\"ahler metrics is generic.
Indeed, take one such a graph. Assume that $d < n-d$; otherwise, turn
the graph up-side-down.  Then, after erasing the edges labeled $1$,
we get the same graph as we do from the one on the left 
with $c=1$, with $d$ as before, and with $n$ replaced by $n-2d$.
This corresponds to a space whose K\"ahler metric is generic.
By the Uniqueness Theorem \ref{thm:uniqII},
these K\"ahler manifolds are isomorphic as Hamiltonian $S^1$-spaces.
\end{Remark}

\begin{figure}
$$
\begin{picture}(150,20)(0,0)
\path(15,10)(15,50)(55,10)(15,10)
\put(15,10){\circle*{4}}
\put(15,50){\circle*{4}}
\put(55,10){\circle*{4}}
\put(0,0){$\ss (0,0)$}
\put(55,0){$\ss (r,0)$}
\put(-10,50){$\ss (0,r)$}
\end{picture}
\quad
\begin{picture}(140,74)(0,-10)
\put(15,50){\blacken\ellipse{4}{4}}
\put(115,50){\blacken\ellipse{4}{4}}
\put(35,10){\blacken\ellipse{4}{4}}
\put(15,10){\blacken\ellipse{4}{4}}
\path(15,50)(115,50)(35,10) (15,10)(15,50)
\put(-15,50){$\ss (0,s)$}
\put(40,0){$\ss (r,0)$}
\put(120,50){$\ss (r+ns,s)$}
\put(0,0){$\ss (0,0)$}
\end{picture}
$$
\caption{Delzant polygons of $\CP^2$ and of a Hirzebruch surface}
\labell{CP2-Hirz}
\end{figure}
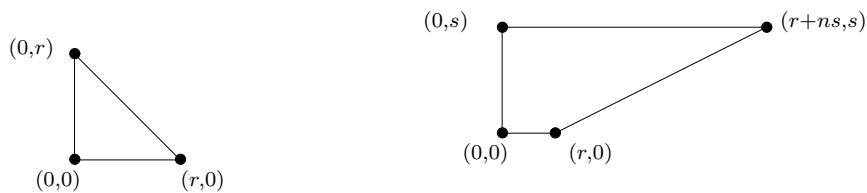

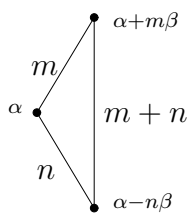
\begin{figure}
$$
\setlength{\unitlength}{0.0100in}
\begin{picture}(72,121)(0,-10)
\put(45,103){\blacken\ellipse{4}{4}}
\put(45,3){\blacken\ellipse{4}{4}}
\put(15,53){\blacken\ellipse{4}{4}}
\path(45,103)(15,53)(45,3)(45,103)
\put(12,73){$m$}
\put(15,18){$n$}
\put(50,48){$m+n$}
\put(0,53){$\ss \alpha$}
\put(55,98){$\ss \alpha + m\beta$}
\put(55,3){$\ss \alpha - n\beta$}
\end{picture}
$$
\caption{Graph for ${\protect \Bbb C} {\protect \Bbb P} ^2$, 
isolated fixed points.}
\labell{CP2-A}
\end{figure}

\begin{figure}
$$
\begin{picture}(60,99)(0,-10)
\put(40,20){\blacken\ellipse{40}{10}}
\put(40,77){\blacken\ellipse{4}{4}}
\put(15,75){$\ss \alpha + \lambda$}
\put(13,15){$\ss \alpha$}
\put(55,5){$\ss \lambda \ , \ g=0$}
\end{picture}
$$
\caption{Graphs for ${\protect \Bbb C}{\protect \Bbb P}^2$,
with a fixed surface.}
\labell{CP2-B}
\end{figure}
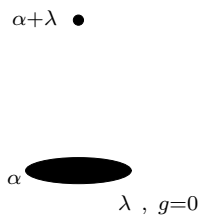

\begin{figure}
$$
\setlength{\unitlength}{0.0125in}
\begin{picture}(115,121)(0,-10)
\put(65,95){\blacken\ellipse{4}{4}}
\put(25,55){\blacken\ellipse{4}{4}}
\put(65,15){\blacken\ellipse{4}{4}}
\put(85,35){\blacken\ellipse{4}{4}}
\path(65,95)(85,35)(65,15) (25,55)(65,95)
\path(65,95)(85,35)(65,15) (25,55)(65,95)
\put(35,70){$c$}
\put(35,30){$d$}
\put(80,65){$nc+d$}
\put(65,100){$\ss r+rc+sd+nsc$}
\put(90,30){$\ss \alpha + rc$}
\put(0,55){$\ss \alpha +sd$}
\put(55,10){$\ss \alpha$}
\put(75,15){$c$}
\end{picture}
\quad \quad \quad \quad
\begin{picture}(90,134)(0,-10)
\put(60,110){\blacken\ellipse{4}{4}}
\put(60,10){\blacken\ellipse{4}{4}}
\put(20,30){\blacken\ellipse{4}{4}}
\put(20,70){\blacken\ellipse{4}{4}}
\path(60,10)(60,110)(20,70) (20,30)(60,10)
\put(5,45){$c$}
\put(65,50){$c$}
\put(30,10){$d$}
\put(0,85){$nc-d$}
\put(65,110){$\ss \alpha+rc$}
\put(60,0){$\ss \alpha$}
\put(-5,20){$\ss \alpha+sd$}
\put(-12,65){$\ss \alpha+sd+$}
\put(-9,60){$\ss rc-ns $}
\end{picture}
\quad \quad
\begin{picture}(99,119)(0,-10)
\put(55,95){\blacken\ellipse{4}{4}}
\put(35,55){\blacken\ellipse{4}{4}}
\put(35,15){\blacken\ellipse{40}{10}}
\path(55,95)(35,55)
\put(60,95){$\ss \alpha+r+ns$}
\put(10,55){$\ss \alpha+r$}
\put(5,15){$\ss \alpha$}
\put(40,0){$\ss s \ , \  g=0$}
\put(35,75){$n$}
\end{picture}
$$
\caption{Graphs for Hirzebruch surfaces}
\labell{Hirz-graphs}
\end{figure}
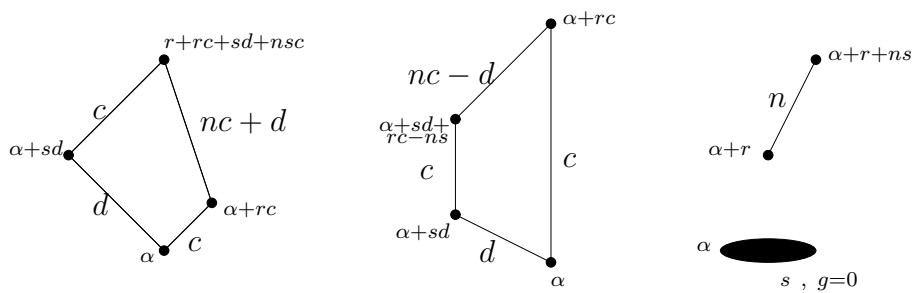

The remaining minimal $S^1$-spaces each have two fixed surfaces
and no interior fixed points.
These will turn out to be {\em ruled manifolds}, which are
symplectic 4-manifolds that are fibered over a closed surface, 
the fibers being two-spheres.
Ruled manifolds play an important role in McDuff's work \cite{md:ruled}.
In the context of Hamiltonian $S^1$-spaces, they come equipped with
fiberwise circle actions:

\begin{Definition} \labell{def:ruled}
A {\em ruled manifold\/} is an $S^2$-bundle over a closed
surface, with a circle action that fixes the base
and rotates each fiber, and with an invariant symplectic form 
and a moment map.
\end{Definition}

\begin{Remark}
Definition \ref{def:ruled} is the one that we use
in the context of Hamiltonian $S^1$-spaces. We would like to 
stress that in other contexts (complex surfaces, non-equivariant
symplectic topology), the definitions are slightly different. 
\end{Remark}

The north and south pole of each fiber fit together into two fixed
surfaces, which are both diffeomorphic to the base surface; 
at all other points, the action is free.

\begin{Lemma} \labell{no-int}
Any compact four dimensional Hamiltonian $S^1$ space
with two fixed surfaces and no interior fixed points
is a ruled manifold.
Moreover, it admits a compatible K\"ahler structure,
and its metric is generic in the sense of Corollary \ref{generic}.
\end{Lemma}

\begin{pf}
Genericity of the metric follows trivially from the
absence of interior fixed points.

By the Uniqueness Theorem \ref{thm:uniqII}, 
the space is determined by the labels of its graph:
the genus of the base space, 
the values of the moment map on the two fixed surfaces,
and the normalized symplectic areas of these surfaces.
We denote these numbers by $g$, $y_\min$, $y_\max$, $a_\min$, 
and $a_\max$, respectively. 
These satisfy the following conditions:
\labell{ruled-parameters}
\begin{enumerate}
\item
$g$ is a non-negative integer;
\item
$a_\min$ and $a_\max$ are positive real numbers;
\item
$y_\min$ and $y_\max$ are real numbers, and $y_\max$ is greater
than $y_\min$;
\item
the number $e_\min = - (a_\max - a_\min) / (y_\max - y_\min)$ 
is an integer.
\end{enumerate}
Conditions 1--3 are clear. Condition 4 follows, e.g., from the 
Guillemin-Lerman-Sternberg formula \eqref{eq:rho}, by which
$a_\max = a_\min - e_\min (y_\max - y_\min)$,
where $e_\min$ is the self intersection of the minimal surface
of the moment map.
The corresponding graphs are shown in Figure \ref{fig:no-int},
in which $a$ and $a+s$ are moment map labels, $r$ and $r+s$
are area labels, and $g$ is the genus.
(If the first graph corresponds to a space,
then $a_\min$, $a_\max$, $e_\min$, $y_\min$, and $y_\max$
are, respectively, $r$, $r+s$, $-n$, $\alpha$, and $\alpha+s$.
Similarly, if the second graph corresponds to a space,
then $a_\min$, $a_\max$, $e_\min$, $y_\min$, and $y_\max$
are, respectively, $r+ns$, $r$, $n$, $\alpha$, and $\alpha+s$.)

\begin{figure}
$$
\begin{picture}(100,75)(0,-10)
\put(35,55){\blacken\ellipse{40}{10}}
\put(35,15){\blacken\ellipse{40}{10}}
\put(10,5){$\ss \alpha$}
\put(45,0){$\ss r \ , \  g$}
\put(0,45){$\ss \alpha+s$}
\put(45,45){$\ss r+ns$}
\end{picture}
\quad
\begin{picture}(100,60)(0,-10)
\put(35,55){\blacken\ellipse{40}{10}}
\put(35,15){\blacken\ellipse{40}{10}}
\put(10,5){$\ss \alpha$}
\put(45,0){$\ss r+ns \ , \  g$}
\put(0,45){$\ss \alpha+s$}
\put(45,45){$\ss r$}
\end{picture}
$$
\caption{Two fixed surfaces and no interior fixed points}
\labell{fig:no-int}
\end{figure}
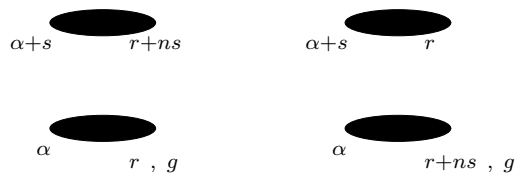

To prove the Lemma, we will show that all the graphs of 
Figure \ref{fig:no-int} come from K\"ahler ruled manifolds.

Fix a non-negative integer $n$, positive real numbers $r$ and $s$,
and a non-negative integer $g$. 
We realize a real surface of genus $g$
as an algebraic curve in $\CP^2$ that does not contain the point $[1,0,0]$:
We can take the zero set of a homogeneous polynomial $P$, 
where $P(w_0,w_1,w_2) = w_0$ if the genus is zero,
and if the genus $g$ is positive, we can take
$P(w_0,w_1,w_2) =  w_0^{2g-2} w_1^2 - \prod_{j=1}^{2g} (w_2 - \eta_j w_1)$,
where $\eta_1 , \ldots, \eta_{2g}$ are distinct nonzero complex numbers.
Then take the manifold
$$\begin{array}{l}
M = \{ ([z_0,z_1,z_2],[w_0,w_1,w_2]) \in \CP^2 \times \CP^2 \ | \  \\
\quad \quad \quad \quad \quad
  P(w_0,w_1,w_2) =0 \mbox{ and } z_1 w_2^n = z_2 w_1^n \} .
\end{array}$$
with the symplectic form
$\omega = \pi_1^* (\ol{r} \omega_{FS}) + \pi_2^* (s \omega_{FS})$,
where $\pi_1$ and $\pi_2$ are the projections to the first and
second $\CP^2$'s, $\omega_{FS}$ is the Fubini Study form on $\CP^2$,
normalized so that $\frac{1}{2\pi} \int_{\CP^1} \omega_{FS} =1$
for $\CP^1 \subset \CP^2$, and $\ol{r}$ is determined by
$\frac{1}{2\pi} \int_{P=0} \ol{r} \omega_{FS} = r$.
Consider the two $S^1$ actions on $M$, given by 
$ ( [\lambda z_0 , z_1 , z_2] , [w_0 , w_1 , w_2] )$
and by 
$ ( [z_0 , \lambda z_1 , \lambda z_2] , [w_0 , w_1 , w_2] )$.
Clearly, $M$ is K\"ahler manifold and the actions are holomorphic.
It is not hard to see that the corresponding graphs are given
by Figure \ref{fig:no-int}, with the same parameters $n$, $r$, $s$,
and $g$ as those in the construction of $M$. Finally,  the projection
$([z_0,z_1,z_2] , [w_0,w_1,w_2]) \mapsto [w_0,w_1,w_2]$
exhibits $M$ as a $\CP^1$ bundle over a surface of genus $g$
with a fiberwise circle action, so $M$ is a ruled manifold.
\end{pf}

\section{Completing the classification; our spaces are K\"{a}hler.}
\labell{sec:existence}

To complete our classification, it remains to specify
which sizes of equivariant symplectic blow-ups can be
performed on a compact four dimensional Hamiltonian $S^1$-space.

The effect of the such a blow-up on the graph of the space
was shown in Figures \ref{blowup-A}--\ref{blowup-D},
where $\lambda$ was the size of the blow-up.
These figures define a blow-up operation on {\em graphs\/}.
If $\lambda$ is large, the resulting graph will fail
to correspond to a space for obvious reasons:
the moment map labels along a chain of edges
might no longer be monotone, or a moment map label at an interior
vertex might become larger or smaller than that on the
maximal or minimal vertex, or the area label associated
to a ``fat vertex" might become non-positive.
We will show that if these problems with the graph do not occur, 
the blow-up can actually be performed on the space.

To make a precise statement, we define a partial ordering
on the vertices of the graph:
for any two vertices $v,w$ with moment map labels $\Phi(v)$ and 
$\Phi(w)$, we declare that $v < w$ if and only if 
$\Phi(v) < \Phi(w)$ and, additionally, either 
$\Phi(v)$ or $\Phi(w)$ is extremal, or that 
$v$ and $w$ are connected by a chain of edges 
along which $\Phi$ is monotone. 

If we now perform a $\lambda$ blow-up on the {\em graph\/},
the vertices of the blown-up graph aquire a natural
partial ordering, which coincides with the one described above
if $\lambda$ is small. We will always work with this same partial 
ordering, even if $\lambda$ is big.

\begin{Definition} \labell{monotone}
The $\lambda$-blown-up graph is {\em monotone\/} if 
\begin{enumerate}
\item
its moment map labels are strictly monotone 
with respect to the partial ordering on its vertices, and
\item
its area-labels are positive.
\end{enumerate}
\end{Definition}

\begin{Proposition} \labell{can blowup}
A four dimensional compact Hamiltonian $S^1$-space 
admits an equivariant symplectic
blow-up of size $\lambda$ at some fixed point 
if and only if the corresponding blown-up graph is monotone.
\end{Proposition}

\begin{Example} \label{constraints}
The graph on the right of Figure \ref{Hirz-graphs}
corresponds to a symplectic Hirzebruch surface
with an $S^1$ action that has one fixed surface
at the minimum of the moment map, one interior fixed
point, and an isolated maximum for the moment map.
Figure \ref{fig-constraints} shows the corresponding 
blown-up graphs.
\begin{itemize}
\item
If we blow up at a point on the minimal surface, 
the area of this surface decreases, and we create
a new interior fixed point.
There are two restrictions
on the amount, $\lambda$, of the blow-up:
the area of the surface must remain positive,
and the moment map value at the new interior fixed point
must be smaller than the moment map value at the maximum.
The corresponding restrictions on $\lambda$ are
$\lambda < s$ and $\lambda < r+ns$.
Since $\min(s, r+ns) = s$, 
we can blow up at the minimal surface 
by any amount $\lambda < s$. 
\smallskip
\item
If we blow up at the interior fixed point,
it gets replaced by two new interior fixed points.
Their moment map values must remain strictly between the moment map 
values at the minimum and at the maximum. This gives the
restrictions $\lambda < r$ and $\lambda < s$.\\
\item
If we blow up at the maximum, 
it gets replaced by two fixed points.
One of them becomes the new maximum.
Its moment map value must stay larger than the 
moment map values at the old interior fixed point;
this gives the restriction $\lambda< ns$.
This new maximum is connected by an edge to a new interior
fixed point, which is  
further connected by an edge to the old interior fixed point.
In order to respect
the partial ordering, the moment map value at the
new interior fixed point must also be larger than
the moment map value at the old interior fixed point.
This imposes the restriction $\lambda < s$.
Since $\min (ns , s) =s$, we can blow up at the maximum
by any amount $\lambda < s$.
\end{itemize}
\end{Example}

\begin{figure}
$$
\begin{picture}(99,119)(0,-10)
\put(35,95){\blacken\ellipse{4}{4}}
\put(35,55){\blacken\ellipse{4}{4}}
\put(35,15){\blacken\ellipse{40}{10}}
\path(35,95)(35,55)
\put(40,95){$\ss \alpha+r+ns$}
\put(40,55){$\ss \alpha+r$}
\put(5,15){$\ss \alpha$}
\put(40,0){$\ss s - \lambda, \ g=0$}
\put(25,70){$\ss n$}
\put(55,45){\blacken\ellipse{4}{4}}
\put(60,45){$\ss \alpha + \lambda$}
\end{picture}
\quad
\begin{picture}(99,119)(0,-10)
\put(35,95){\blacken\ellipse{4}{4}}
\put(35,65){\blacken\ellipse{4}{4}}
\put(35,45){\blacken\ellipse{4}{4}}
\put(35,15){\blacken\ellipse{40}{10}}
\path(35,95)(35,65)(35,45)
\put(40,95){$\ss \alpha+r+ns$}
\put(40,65){$\ss (\alpha+r) + n \lambda$}
\put(40,45){$\ss (\alpha+r) - \lambda$}
\put(5,15){$\ss \alpha$}
\put(15,55){$\ss n+1$}
\put(40,0){$\ss s,\ g=0$}
\put(25,80){$\ss n$}
\end{picture}
\quad
\begin{picture}(99,119)(0,-10)
\put(35,85){\blacken\ellipse{4}{4}}
\put(35,70){\blacken\ellipse{4}{4}}
\put(35,55){\blacken\ellipse{4}{4}}
\put(35,15){\blacken\ellipse{40}{10}}
\path(35,85)(35,70)(35,55)
\put(40,85){$\ss (\alpha+r+ns) - \lambda$}
\put(40,70){$\ss (\alpha+r+ns) - n \lambda$}
\put(40,55){$\ss \alpha+r$}
\put(5,15){$\ss \alpha$}
\put(40,0){$\ss s,\ g=0$}
\put(15,77){$\ss n-1$}
\put(25,63){$\ss n$}
\end{picture}
$$
\caption{Blow-ups of a Hirzebruch surface for Example \ref{constraints}}
\labell{fig-constraints}
\end{figure}
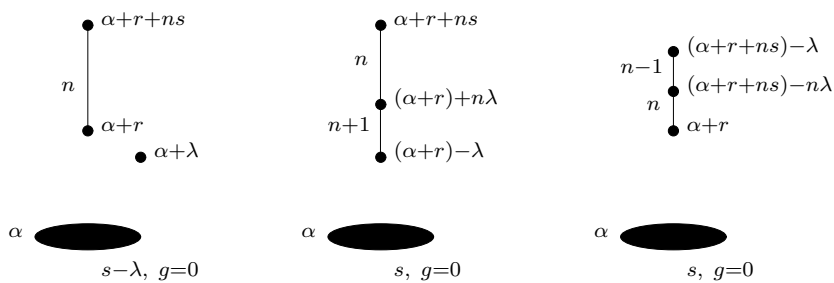

\begin{Example} 
[A space with one fixed surface is not necessarily toric]
\labell{ex:X2}
In Example \eqref{ex:X}, 
the blow-up can be made of size exactly $1$,
by Proposition \ref{can blowup},
The resulting space would have three fixed points on the same 
level set moment for the moment map,
and so it would have three non-trivial chains of gradient
spheres for {\em any\/} choice of compatible metric.
Moreover, it would not admit any second circle action
that commutes with the first. Hence a space with one
fixed surface is not necessarily toric.
\end{Example}

At the same time that we will prove Proposition \ref{can blowup}, 
we will show that our spaces are K\"ahler:

\begin{Theorem} \labell{kahler}
Every compact Hamiltonian $S^1$ space is K\"ahler, i.e., admits
a complex structure such that the action is holomorphic and the
symplectic form is K\"ahler.
Moreover, this can be done in such a way that the associated
Riemann metric is generic in the sense of Corollary \ref{generic}.
\end{Theorem}

Theorem \ref{big} enables us to prove Theorem \ref{kahler} 
and Proposition \ref{can blowup} by induction on the number of blowups: 
we already saw, in \S\ref{subsec:min-spaces},
that each minimal space admits a K\"ahler structure
with a generic metric. 
We will show, for a space that is K\"ahler and whose metric is generic,
that if its $\lambda$-blown-up graph is monotone, 
there exists an invariant K\"ahler form on the equivariant 
complex blow-up that gives rise to this blown-up graph
and whose metric is generic.

It follows that the original space admits a $\lambda$-blow-up,
because, by the Uniqueness Theorem \ref{thm:uniqII}, 
the original space is isomorphic
to the $\lambda$-blow-down of the new space,
since their graphs are the same.

We will need to understand the arrangement of the invariant complex 
curves inside a complex surface with a holomorphic $S^1$ action.
This will be
similar to the arrangement of the gradient spheres in 
a Hamiltonian $S^1$-space. Here we can't speak of the gradient
flow of the moment map (we don't yet have a moment map nor a 
metric), but instead we have the flow generated by $-J \xi_M$,
where $\xi_M$ generates the circle action and $J : TM \to TM$ 
is the complex structure; in the presence of an invariant
K\"ahler form and a moment map, this flow is the same as the
gradient flow for the moment map.

\begin{Definition} \labell{notation}
Let $M$ be a complex surface with a holomorphic $S^1$ action.
A fixed point in $M$ is called {\em maximal\/}
if its isotropy weights are non-positive, {\em minimal\/}
if its isotropy weights are non-negative, and {\em interior\/}
if one of its isotropy weights is positive and the other is negative.
Let $C$ be a connected invariant complex curve in $M$.
If $C$ is not fixed by the action, then either the action on $C$
is free and $C$ is a torus, or the action on $C$ contains two fixed
points and $C$ is a two-sphere. In the latter case, 
$C$ is called a {\em gradient complex curve\/},
and the two fixed points on it are called its 
{\em north and south poles\/}.
We distinguish them in the following way: 
if $\xi_M$ is the vector field that generates the circle
action, the trajectories of the vector field $-J \xi_M$
in $C$ approach the north pole as time goes to infinity
and approach the south pole as time approaches negative infinity.
A {\em chain of gradient complex curves\/} is a sequence of such
curves $C_1, \ldots, C_l$ in which the north pole of $C_i$
is the south pole of $C_{i+1}$ for $i=1,\ldots,l-1$,
the south pole of $C_1$ is minimal, and the north pole of $C_l$ is
maximal.
\end{Definition}

\begin{Proposition} \labell{arrangement}
Let $M$ be a complex surface with a holomorphic $S^1$ action,
obtained from $\CP^2$, a Hirzebruch surface, or a ruled surface
by a sequence of equivariant complex blow-ups at fixed points.
The set of maximal points consists of either one point
or one connected complex curve that is fixed by the action.
The same holds for the set of minimal points.
Every point in $M$ that is not an isolated extremum
belongs to exactly one chain of gradient complex curves.
\end{Proposition}

\begin{pf}
This is easily seen by induction on the number of equivariant 
blowups. 
\end{pf}

Our next goal is to identify the K\"ahler cone of such a complex surface:

\begin{Proposition} \labell{kahler-cone}
Let $M$ be a complex surface with a holomorphic $S^1$ action, 
obtained from $\CP^2$, a Hirzebruch surface, or a ruled surface
by a sequence of equivariant complex blowups at fixed points.
Let $\Omega$ be a cohomology class on $M$
that is positive on every invariant complex curve.
Then $\Omega$ contains a K\"ahler form.
\end{Proposition}

We prove Proposition \ref{kahler-cone} in Appendix \ref{app:kahler-cone}.

To proceed with the proof of Proposition \ref{can blowup},
we construct a cohomology class out of the blown-up graph:

\begin{Lemma} \labell{omega-tilde}
Let $(M,\omega,\Phi)$ be a four dimensional compact Hamiltonian 
$S^1$-space, with a compatible K\"ahler metric, that is generic
in the sense of Corollary \ref{generic}. 
Take a $\lambda$-blow-up of its graph at some vertex.
Let $p \in M$ be a fixed point that corresponds to that vertex.
If $p$ lies on a fixed surface, choose it to not be the end-point
of any non-trivial chain of gradient complex curves.
Let $\pi : \tM \to M$ be the equivariant complex blow-up at $p$.
Let $\tE_p = \pi \inv(p)$ be the exceptional divisor,
let $\Xi_p$ be its Poincar\'e dual in $H^2(\tM)$,
and let $\Omega$ be the cohomology class of $\omega$.
Let 
$$
	\tOmega = \pi ^* \Omega - \lambda \Xi_p.
$$
Then
\begin{enumerate}
\item
The labels of the $\lambda$-blown-up graph correspond to the
cohomology class $\tOmega$ in the following way:
\begin{enumerate} 
\item
For any fixed surface $\tC$ in $\tM$, the corresponding area label
is equal to $\frac{1}{2\pi} \tOmega(\tC)$.
\item
For any invariant two-form $\tomega$ in the class $\tOmega$
for which the action is Hamiltonian,
there exists a moment map $\tPhi$ whose values at the fixed points
are equal to the corresponding moment map labels of the 
$\lambda$-blown-up graph.
\end{enumerate} 
\item
$\tOmega$ is positive on invariant complex curves
if and only if the $\lambda$-blown-up graph is monotone. 
\end{enumerate}
\end{Lemma}

\begin{pf}
Let us first show how Part 2 follows from Part 1. 
An invariant complex curve is either a gradient complex curve
or a fixed complex curve.
The second condition of monotonicity is equivalent, by part 1(a),
to the positivity of $\tOmega$ on any fixed complex curve.
For a gradient complex curve $C$ with north pole $p$ and south pole $q$
and with a stabilizer of order $k$,
\begin{equation} \labell{pos-C} 
 \int_C \tomega = \frac{1}{k} (\Phi(p) - \Phi(q)),
\end{equation}
by the proof of Lemma \ref{cohomology}.
By part 1(b), the first condition for monotonicity 
is equivalent to the right hand side of \eqref{pos-C} being positive 
for every such $C$, hence it is equivalent to $\tOmega$ being positive
on all gradient complex curves.

Part 1 follows from the fact that for any complex curve $\tC$ in $\tM$,
\begin{equation} \labell{tOmega}
	\frac{1}{2\pi} \tOmega (\tC) = \frac{1}{2\pi} \Omega(C)
					- \lambda \tC \cdot \tE_p
\end{equation}
where $C = \pi(\tC)$ is either a curve or a point.

Let us prove 1(a).
Suppose that $\tC$ is fixed. If $\tC$ does not meet $\tE_p$, its area label
is the same as the corresponding area label in the graph before the blow-up.
This label is $\Omega(C)$, which is equal to $\tOmega(\tC)$ by \eqref{tOmega}
with $\tC \cdot \tE_p = 0$. 
If $\tC = \tE_p$, then by Figure \ref{blowup-D}, its area label is $\lambda$,
which is equal to $\tOmega(\tC)$ by \eqref{tOmega} with 
$\Omega(C) = \Omega(\text{point}) = 0$ 
and $\tC \cdot \tE_p = \tE_p \cdot \tE_p = -1$. 
If $\tC$ meets $\tE_p$ but is not equal to it, then the point $p$ at which
we blew up is either the north or south pole of $C$.
By figure \eqref{blowup-B} with $v=\Omega(C)$, the area label of $\tC$ 
is then $\Omega(C)-\lambda$, which is equal to $\tOmega(\tC)$ by 
\eqref{tOmega} with $\tC \cdot \tE_p = 1$.

To prove 1(b), by \eqref{pos-C} it is enough to show that for any complex
gradient sphere $\tC$ with a stabilizer of order $k$, the difference
the difference between the moment map labels corresponding to the north and
south poles of $\tC$ is equal to $k\tOmega(\tC)$. If $\tC$ does not meet 
$\tE_p$, this difference is $\Phi(p) - \Phi(q) = 
tPhi(p) - \tPhi(q) = k \Omega(C) = k \tOmega(\tC)$.
If $\tC = \tE_p$, by Figure \eqref{blowup-A} with $k=m+n$, this difference
is $k\lambda$, which is equal to $k \tOmega(\tC)$ by \eqref{tOmega} with
$\Omega(C) = \Omega(\text{point}) =0$ 
and $\tC \cdot \tE_p = \tE_p \cdot \tE_p =1$.
If $\tC$ meets $\tE_p$ but is not equal to it, by Figure \ref{blowup-A},
when we blow up, the moment map at one pole o f$\tC$ moves closer to the  
other by an amount of $k\lambda$ (with $k=m$ or $k=n$). Since the old
difference was $k\Omega(C)$, the new difference is $k\Omega(\tC)$ by 
\eqref{tOmega} with $\tC \cdot \tE_p =1$.
\end{pf}

\begin{pf*}{Proof of Proposition \ref{can blowup} and of Theorem 
\ref{kahler}}

We already proved Theorem \ref{kahler} for a minimal space.
Let us suppose that Theorem \ref{kahler} holds for a space,
and suppose that a $\lambda$-blow-up of the corresponding
graph is monotone. (See Definition \ref{monotone}.)

Let $\tOmega$ be the cohomology class constructed 
in Lemma \ref{omega-tilde} on the corresponding 
complex blow-up. By part 2 of Lemma \ref{omega-tilde},
$\tOmega$ is positive on invariant complex curves.
By Proposition \ref{kahler-cone}, $\tOmega$ contains
a K\"ahler form.
By averaging with respect to the holomorphic circle action,
we get an {\em invariant\/} K\"ahler form, $\tomega$,
in the class $\tOmega$. By Lemma \ref{exist-mom}, 
the $S^1$-action on $(\tM,\tomega)$ is Hamiltonian.
By Part 1 of Lemma \ref{omega-tilde}, we can choose
a moment map $\tPhi$ such that the graph of the 
Hamiltonian $S^1$-space $(\tM,\tomega,\tPhi)$
is the $\lambda$-blown-up graph of $(M,\omega,\Phi)$.

The $\lambda$-blow-down of this space exists,
and has the same graph as that of $(M,\omega,\Phi)$.
By the Uniqueness Theorem \ref{thm:uniqII},
the $\lambda$-blow-down of $(\tM,\tomega,\tPhi)$ 
is isomorphic to $(M,\omega,\Phi)$; hence, the 
$\lambda$-blow-up of $(M,\omega,\Phi)$ exists
and is isomorphic to $(\tM,\tomega,\tPhi)$;
in particular, it is K\"ahler.

This proves Proposition \ref{can blowup} for $(M,\omega,\Phi)$
and Theorem \ref{kahler} for any blow-up $(\tM,\tomega,\tPhi)$
of $(M,\omega,\Phi)$. Applying Theorem \ref{big}, we finish 
by induction on the number of blow-ups.
\end{pf*}

{
\subsection{Algorithm}

Theorem \ref{big} and Proposition \ref{can blowup}
provide algorithms for constructing all the graphs of compact,
four dimensional Hamiltonian $S^1$-spaces: start with 
one of the graphs of Figures \ref{CP2-A}--\ref{fig:no-int},
possibly turned up-side-down,
and perform a sequence of $\lambda$-blow-ups,
as described in Figures \ref{blowup-A}--\ref{blowup-D},
possibly turned up-side-down.
At each stage, the size $\lambda$ of the blowup must be
small enough so that the graph remains monotone
(see Definition \ref{monotone}), i.e., the moment map labels
must increase along each branch, they must be maximal and minimal
at the top and bottom vertices, and the area labels must remain
positive.

\labell{algorithm}
If the extremal sets for the moment map are both surfaces,
there exists an alternative algorithm, in which
we can choose the real labels all at once, at the end:
\begin{enumerate}
\item[(i)]
Start with two ``fat" vertices,
corresponding to the minimum and maximum of $\Phi$,
with arbitrary equal genus labels.
\item[(ii)]
Perform a sequence of blow-ups 
as in Figures \ref{blowup-A}, \ref{blowup-B}, and \ref{blowup-B}
turned up-side-down, while only keeping track of integer labels.
\item[(iii)]
Assign moment map labels $y_\min < y_\max$ to the two fat
vertices. To all other vertices, assign moment map labels
that are between $y_\min$ and $y_\max$ and that are increasing along 
each branch of edges.
\item[(v)]
Choose positive real labels $a_\min,a_\max$
subject to the condition $a_\min - a_\max = b$,
where $b$ is determined in the following way.

Let $m_p,n_p$ be the isotropy weights at the interior
fixed point $p$. The sum
$a := \sum_p { 1 \over m_p n_p} $
is an integer, by Lemma \ref{kho}
(in which $k_1 = k_l =1$ whenever the minimum and maximum
of the moment map are attained on surfaces).
Equation \eqref{eq:rho} imposes the condition
$$
        0 = - (a + e_\min + e_\max) y + (a_\min - a_\max -b)
$$
for all $y > y_\max$,
where $a$ is as above and where $b$ is a certain function
of the integer and real labels and of $e_\min$ and $e_\max$.
Choose any integers $e_\min, e_\max$ with $e_\min + e_\max = -a$;
this determines $b$.
\end{enumerate}

}

\appendix

\section{Local normal forms}
\labell{sec:normal}

\begin{Definition}
Let $(M,\omega)$ be a symplectic manifold and $N \subseteq M$
a submanifold. Suppose that the two-form $\omega$
has constant rank on $TN$. The {\em symplectic normal bundle\/}
of $N$ in $M$ is the quotient bundle
$ 
	(TN)^\omega / (TN \cap (TN)^\omega),
$
where $(\cdot)^\omega$ denotes the symplectic ortho-complement.
This is a symplectic vector bundle over $N$,
i.e., a vector bundle equipped with a fiberwise symplectic form.
\end{Definition}

We recall the equivariant constant rank embedding theorem
(see \cite{ma}).

\begin{Theorem} \labell{normal}
Let $(M_1,\omega_1)$ and $(M_2,\omega_2)$ 
be symplectic manifolds with actions of a compact Lie group $G$. 
Suppose that we have $G$-equivariant embeddings of a manifold $N$ 
in $M_1$ and in $M_2$,
such that the pullbacks of $\omega_1$ and $\omega_2$ to $N$
are equal and have constant rank, and such that the symplectic
normal bundles of $N$ in $M_1$ and in $M_2$ are isomorphic as $G$-equivariant
symplectic vector bundles. Then the map identifying the image of $N$
in $M_1$ with the image in $M_2$ extends to an equivariant symplectomorphism
between neighborhoods of these images.
\end{Theorem}


It is extremely useful to know that symplectic vector bundles
are the same as complex vector bundles, up to isomorphism.
We recall the proof of this well known fact, because we did not find
a careful reference that also kept track of a group action.

\begin{Definition}
Let $E \to M$ be a $G$-equivariant symplectic vector bundle.
A {\em compatible complex structure\/} on $E$ is a $G$-equivariant 
automorphism $J : E \to E$ with $J^2 = -\text{id}$ 
such that the equation
\begin{equation} \labell{compatible}
\brak{u,v} = \omega(u,Jv) 
\end{equation}
defines an inner product (symmetric positive definite bilinear form) 
on each fiber.  
This inner product is then called a {\em compatible metric}.
\end{Definition}

\begin{Lemma} \labell{standard}
Any $G$-equivariant symplectic vector bundle, $E \to M$, admits
a compatible complex structure.
Moreover, given any compatible complex structure over a neighborhood
of an invariant closed subset $K$ of $M$, there exists a
compatible complex structure on all of $E$ that coincides with
the given one over some smaller neighborhood of $K$.
Finally, any two compatible complex structures on $E$ are isomorphic 
as $G$-equivariant complex vector bundles over $M$.
\end{Lemma}

\begin{pf}
(Following \cite[p.8]{w:lectures}.)  
Start with any fiberwise metric on $E$. 
Average it to obtain a $G$-invariant metric, $\brak{\,,\,}$.
The relation $\brak{u,v} = \omega(u,Av)$ defines a fiberwise automorphism
$A$ of $E$ such that $-A^2$ is symmetric and positive definite.
Let $P = \sqrt{-A^2}$ and $J = A P^{-1}$.
Then $J$ is invariant and is compatible with $\omega$.

Given $J$ over a neighborhood of $K$, extend its corresponding metric
to {\em some\/} metric on $M$, without changing it on some smaller 
neighborhood of $K$. The procedure of the previous paragraph then gives
a compatible
metric that coincides with the original metric over a neighborhood of $K$.

If we have two compatible $J$'s, their corresponding metrics can be 
connected by a smooth path, for instance, by taking convex combinations.
Applying the above construction simultaneously to all
the metrics in this path produces a smooth path of equivariant
complex structures, $J_t$.  Finally, we can construct a $G$-equivariant
isomorphism from $(E,J_0)$ to $(E,J_1)$ in the following way:

Consider the $G$-equivariant bundle over $M \times [0,1]$ 
whose fiber over $(m,t)$ consists of the space of complex linear isomorphisms 
$\varphi : (E_m,J_0) \to (E_m,J_t)$.
The identity map provides a smooth $G$-equivariant
section of this bundle over the subset $M \times \{ 0 \}$.  This section
extends to a smooth $G$-equivariant section of the whole bundle.
(For instance, lift the vector field $-\frac{\partial}{\partial t}$
to the bundle, average to produce a lifting that is $G$-equivariant, 
and take the trajectories of this vector field that start at the section 
over $\{ 0 \} \times M$.)  Over $M \times \{ 1 \}$, this section
provides an equivariant isomorphism from $(E,J_0)$ to $(E,J_1)$.
\end{pf}

The converse is also true: 

\begin{Lemma} \labell{standardII}
Every $G$-equivariant complex vector bundle over $M$
admits a $G$-equivariant compatible (fiberwise) symplectic structure,
and this structure is unique up to isomorphism of $G$-equivariant
symplectic vector bundles.
\end{Lemma}

\begin{pf}
From any metric on the fibers of $E$ we can obtain a compatible metric,
and hence a compatible symplectic structure,
by averaging with respect to $J$ and $G$. 

Given two compatible symplectic structures, $\omega_0$ and $\omega_1$,
the convex combinations of their corresponding metrics
form a smooth path of compatible metrics, hence of compatible
symplectic structures, $\omega_t$.
This defines a $G$-equivariant bundle over $M \times [0,1]$,
whose fiber over $(m,t)$ consists of the space
of symplectic linear transformations from $(E_0,\omega_t)$
to $(E_m,\omega_t)$. The natural section over $t=0$ extends to a global
section, whose restriction to $t=1$ gives the required isomorphism between
$(E,\omega_0)$ and $(E,\omega_1)$.
\end{pf}

Putting together Lemma \ref{standard} and Lemma \ref{standardII}, 
we get

\begin{Corollary} \labell{s=c}
Every equivariant symplectic vector bundle admits a compatible
equivariant complex structure, unique up to isomorphism.
Two equivariant symplectic vector bundles are isomorphic 
if and only if they are isomorphic as equivariant complex vector bundles.
\end{Corollary}

Now that we understand that equivariant symplectic vector bundles are the 
same as equivariant complex vector bundles, we turn to
corollaries of Theorem \ref{normal}.
We fix a compact four dimensional symplectic manifold 
$(M,\omega)$ with a Hamiltonian circle action.
We first describe the neighborhood of one fixed point:

\begin{Corollary} \labell{stab=S1}
Let $p \in M$ be a fixed point.  Then there exist complex coordinates $z,w$ on
a neighborhood of $p$ in $M$, and unique integers $m$ and $n$, such that
\begin{equation}
\labell{normalI}
\parbox{5.0in}{
\begin{itemize}
\item[(i)] 	the circle action is 
		$\lambda \cdot (z,w) = (\lambda^m z, \lambda^n w),$ 
\item[(ii)] 	the symplectic form is 
	$\omega = {i \over 2} (dz \wedge d\ol{z} + dw \wedge d\ol{w}),$ 
\item[(iii)] 	the moment map is $\Phi(z,w) = 
		\Phi(p) + {m\over 2} |z|^2 + {n\over 2} |w|^2$.
\end{itemize}
}
\end{equation}
The integers $m,n$ are called the {\em isotropy weights} at $p$. 
They are relatively prime, because if a compact Lie group acts effectively 
on a connected manifold, it acts effectively on any invariant open subset. 
\end{Corollary}

By \eqref{normalI}, every connected component of the fixed point set is
either an isolated fixed point or is a symplectic surface, and if it
is a surface, it must be a local maximum or minimum for $\Phi$.  
By a theorem of Atiyah, Guillemin and Sternberg \cite{g-s:convexity,at}, a
local extremum is also a global extremum and is attained on a connected set.
(See \cite[lemma 5.1]{g-s:convexity}.)
These properties were summarized in Lemma \ref{fixed}.

We now describe the neighborhood of a whole fixed surface:

\begin{Corollary} \labell{normal-near-B}
Let $B \subset M$ be a two dimensional component of the fixed point set.
Then the neighborhood of $B$ in $M$ is determined up to
equivariant symplectomorphism by the genus and total area of 
the symplectic surface $(B,\omega)$,
its self intersection in $M$, and whether it is a minimum or maximum
for the moment map.
\end{Corollary}

We also need to work with finite stabilizers. Denote by $Z_k$
the cyclic subgroup of $S^1$ of order $k$.

\begin{Corollary} \labell{stab=Zk}
Let $p \in M$ be a point whose stabilizer is $Z_k$ for some $k \geq 2$. 
Then there exists an integer $l$, relatively prime to $k$, 
and an interval $I$ around $0$, and a disk $D^2$ in $\C$, such that
a neighborhood of $S^1 \cdot p$ in $M$ is $S^1$-equivariantly
symplectomorphic to the model $S^1 \times_{Z_k} (I \times D^2)$,
in which $[e^{i\theta} \nu,  h, z] = [e^{i\theta}, h, \nu^l z]$
for every $\nu \in Z_k$ and $(e^{i\theta},h,z) \in S^1 \times I \times D^2$,
with the symplectic form 
$dh \wedge d\theta + \frac{i}{2} dz \wedge d\zbar$,
and with the left $S^1$-action.
\end{Corollary}

Corollaries \ref{stab=S1} and \ref{stab=Zk} imply that each connected
component of the closure of the set of points
with stabilizer $Z_k$ is a closed symplectic surface
with a Hamiltonian action of the quotient circle $S^1 / Z_k$.
Delzant's classification \cite{de} implies that such a surface 
is equivariantly symplectomorphic to the sphere
$\{ (x_1,x_2,x_3) \in \R^3 \ | \ x_1^2 + x_2^2 + x_3^2 = 1 \}$,
with a constant multiple of its standard area form,
with the circle action that becomes
$\lambda \cdot (z,x_3) = (\lambda^k z , x_3)$
when we identify $\R^3$ with $\C \times \R$ by setting
$z = x_1 + ix_2$, and with the moment map being a constant 
plus a multiple of the height function $(x_1,x_2,x_3) \mapsto x_3$.
We call such a sphere a {\em $Z_k$-sphere\/}.
These facts were summarized in Lemma \ref{Zk}.

\begin{Remark} \labell{l}
If $C \subset M$ is a $Z_k$-sphere, its stabilizer group $Z_k$ 
acts on the fibers of the normal bundle to $C$ by
$$ \nu : z \mapsto \nu^l z  \quad , \quad \nu \in Z_k,$$ 
where $l$ is the same integer that occurs in Corollary \ref{stab=Zk}.
This integer is determined by the isotropy weights at the poles of $C$:
if the isotropy weights at the north pole are $m$ and $-k$
and those at the south pole are $k$ and $-n$,
then both $m$ and $-n$ are congruent to $l$ modulo $k$.
\end{Remark}

The neighborhood of a $Z_k$-sphere is determined by the following

\begin{Corollary} \labell{nbhd-Zk}
Let $C \subset M$ be an invariant symplectic sphere
that is not fixed by the circle action.
Then its neighborhood is determined, up to an equivariant
symplectomorphism, by its total area, its self intersection, 
its stabilizer, and the weights of the circle action on the fibers 
over the north and south poles.
\end{Corollary}

The following normal form will be applied in section \ref{sec:metrics}
when $C$ is a free gradient sphere:

\begin{Corollary} \labell{stab=e}
Consider a free orbit in $M$, and let $C \subset M$ be a two-dimensional, 
symplectic, $S^1$-invariant submanifold of $M$, containing this orbit.
By Theorem \ref{normal}, a neighborhood of the orbit inside $C$
is equivariantly symplectomorphic to $S^1 \times I$,
where $I$ is an open interval, $S^1$ acts on the left factor,
and the symplectic structure is $dh \wedge d\theta$,
where $\theta \mod 2\pi$ is a coordinate on $S^1$
and $h$ is a coordinate on $I$.
Applying Theorem \ref{normal} again, a neighborhood of the orbit in $M$ 
is equivariantly symplectomorphic to $S^1 \times I \times D^2$,
with the symplectic form $dh \wedge d\theta + dx \wedge dy$,
where $x$ and $y$ are coordinates on the disc $D^2$,
and $C$ intersects this neighborhood in $x=y=0$.
\end{Corollary}

We complete this section with a local model that follows from the slice theorem
for compact group actions, not from Theorem \ref{normal}. 
We fix a compatible metric on $M$; the gradient flow for the moment map 
then combines with the $S^1$-action into an action of $S^1 \times \R \cong \Cc$.
We provide a model for a neighborhood of a $\Cc$-orbit that keeps track of the gradient 
flow, the circle action, and the moment map, and ignores the symplectic form:

{
\begin{Lemma} \labell{CI}
Let $C$ be a two-dimensional $\Cc$-orbit in $M$.
Its moment image is an open interval; denote it $I = \Phi(C)$.
\begin{enumerate}
\item
Suppose that the stabilizer of $C$ is trivial.  Consider the model
\begin{equation} \labell{S1 times I times D2}
 S^1 \times I \times D^2,
\end{equation}
with $S^1$ acting on the left factor.
Then an invariant neighborhood of $C$ in $M$
is equivariantly diffeomorphic to this model,
in such a way that the moment map is the projection to $I$
and the gradient flow only effects the $I$ coordinate.
\item
Suppose that the stabilizer of $C$ is $\Z_k$, with $\nu \in Z_k$ acting 
on the normal bundle by multiplication by $\nu^l$, as in Remark \ref{l}.  
Consider the model
\begin{equation} \labell{model Zk}
  S^1 \times_{Z_k} (I \times D^2),
\end{equation}
in which $[e^{i\theta} \nu , h , z ] = [e^{i\theta} , h , \nu^l z ]$
for all $\nu \in Z_k$.
Then an invariant neighborhood of $C$ in $M$ is equivariantly 
diffeomorphic to this model, in such a way that the moment map
is $h$ and the gradient flow only effects $h$.
\end{enumerate}
\end{Lemma}

\begin{pf}
Part 2 would follow from part 1 by passing to a $k$-fold covering 
and working $Z_k$-equivariantly.  To prove part 1, fix an orbit in $C$.
A neighborhood of this orbit inside its level set is $S^1$-equivariantly
diffeomorphic to $S^1 \times D^2$.
If we choose this neighborhood small enough so that it does
not meet any other non-trivial gradient sphere,
the inclusion map of $S^1 \times D^2$ into $M$ 
extends uniquely to a map from the model
\eqref{S1 times I times D2} to $M$ that satisfies our requirements.
\end{pf}
}

\section{Diffeomorphisms of the two-sphere}
\labell{sec:smale}

Smale proved \cite{s} that the topological space of orientation preserving
diffeomorphisms of $S^2$ deformation retracts onto the subspace
$SO(3)$ of rigid rotations.
In \S\ref{sec:uniqueness} we needed a variant of Smale's theorem:
we needed to show that if $\Sigma$ is a two dimensional orbifold,
homeomorphic to $S^2$, and with at most two singular points,
and if $k : \Sigma \to \Sigma$ is a diffeomorphism,
then there exists a diffeotopy $k_t : \Sigma \to \Sigma$
with $k_1 = k$ and $k_0 =$identity.
Ahara and Hattori \cite[\S 10]{ah-ha}
sketched how to reduce this to the following Lemma, 
on diffeomorphisms of the two-sphere:

\begin{Lemma} \labell{smale+}
Let $k : S^2 \to S^2$ be diffeomorphism. 
\begin{enumerate}
\item
Suppose that $k$ fixes a neighborhood of the north pole.
Then there exists a diffeotopy $k_t :S^2 \to S^2$,
with $k_0 = k$ and $k_1 = \text{identity}$, such that
each $k_t$ fixes a neighborhood of the north pole.
\item
Suppose that $k$ acts by a rotation on a neighborhood of the north pole
and on a neighborhood of the south pole.
Then there exists a diffeotopy $k_t :S^2 \to S^2$,
with $k_0 = k$ and $k_1 = \text{identity}$, such that
each $k_t$ acts by a rotation on a neighborhood of each pole.
\end{enumerate}
\end{Lemma}

Smale's proof actually provides a continuous family $k_t$
of diffeomorphisms that satisfy $1$; however, his proof does not
guarantee a smooth dependence on the parameter $t$, which is what
we mean by {\em diffeotopy}.
As for Part 2 above, in spite of what Ahara and Hattori seem to imply,
this result cannot be strengthened to allow to fix neighborhoods
of the two points, and we don't see how to deduce this result from
Smale's proof, even if we don't demand smoothness in $t$.
Despite these reservations, we are confident, just like Ahara and Hattori 
were, that both parts of Lemma \ref{smale+} are well known, although we
didn't find them in the literature. 

Parts 1 and 2 easily follow from the following variants:

\begin{Lemma}
\begin{enumerate}
\item[1'.]
Let $S$ be the space of diffeomorphisms of a square that fix a
neighborhood of the boundary. Then each $k \in S$ can be connected 
to the identity by a diffeotopy $k_t$ in $S$.
\item[2'.]
Let $S'$ be the space of diffeomorphisms of the cylinder 
$S^1 \times I$ that act by a rotation on the neighborhood 
of each component of the boundary. Then each $k \in S'$ can be 
connected to the identity by a diffeotopy $k_t$ in $S'$.
\end{enumerate}
\end{Lemma}

It is easy to modify Smale's argument to obtain a diffeotopy $k_t$, 
as required in part 1'.  \footnote{Moreover, 
it is easy to do this for a whole 
family of $k$'s, parametrized by a finite dimensional manifold.
We don't see how to do this simultaneously for the {\em whole\/} space 
of diffeomorphisms, but we don't need such a strong result.}
We now sketch a proof of part 2', following Smale's ideas.

\begin{pf*}{Sketch of proof of part 2'}
Consider the cylinder $S^1 \times I$, with coordinates
$\theta\mod 2\pi$ and $h \in [0,1]$.
Let $k \in S'$ be a diffeomorphism of the cylinder
which is a rotation near each boundary component, i.e., 
which satisfies
$k(\theta, h) = (\theta + a, h)$ for $h$ near $0$
and  $k(\theta, h) = (\theta + b, h)$ for $h$ near $1$,
for some constants $a$ and $b$.
Then $k_* \dd{h} = \dd{h}$ near the bottom and near the top
of the cylinder.

We trivialize the tangent bundle of the cylinder using the global frame
$(\dd{h},\dd{\theta})$. The vector field $k_* \dd{h}$ then becomes
a function from $S^1 \times I$ to $\R^2 \ssminus \{ 0 \}$
which is constant on neighborhoods of the bottom and of the top of 
the cylinder. It follows that this function is homotopically trivial;
specifically, there exists a smooth one parameter family of vector fields, 
$\xi_t$, with $\xi_0 = k_* \dd{h}$ and $\xi_1 = \dd{h}$, such that 
$\xi_t = \dd{h}$ on neighborhoods of the bottom and of the top of the
cylinder, for all $t$.
For each $t$, the trajectories of $\xi_t$ that ascend from the bottom 
of the cylinder fill the entire cylinder.
We can reparametrize the $\xi_t$'s so that all these trajectories
have length $1$, for each $t$.
Define $k'_t$ to be the diffeomorphism for which 
$(k'_t)_* \dd{h} = \xi_t$. This defines a diffeotopy between $k$
and the identity which is a rotation near the bottom of the cylinder
for all $t$. Near the top of the cylinder, each $k_t$
is given by a diffeomorphism of $S^1$ times the identity map on $I$;
unfortunately, this diffeomorphism might not be a rotation.

We easily fix this by noting that any family of diffeomorphisms
of $S^1$ is diffeotopic to a family of rotations.
Indeed, after composing with a family of rotations
we get a family of diffeomorphisms of $S^1$ that lifts to a family
of periodic diffeomorphisms $g : \R \to \R$ which satisfy $g(0) =0$. 
We can diffeotope these to the identity via
$g_t(\theta) = (1-s) g(\theta) + s \theta$, where 
$s = s(t)$ is zero for $t$ near $0$ and is one for $t$ near $1$.  
\end{pf*}

\section{Computing a K\"ahler cone}
\labell{app:kahler-cone}

In this appendix we prove Proposition \ref{kahler-cone}:

\medskip
\noindent
{\bf Proposition \ref{kahler-cone}.} 
{\em
Let $M$ be a complex surface with a holomorphic $S^1$ action, 
obtained from $\CP^2$, a Hirzebruch surface, or a ruled surface
by a sequence of equivariant complex blowups at fixed points.
Let $\Omega$ be a cohomology class on $M$
that is positive on every invariant complex curve.
Then $\Omega$ contains a K\"ahler form.
}
\medskip

We will use Nakai's criterion:

\begin{Lemma}[Nakai's criterion] \labell{nakai}
A cohomology class $\Omega$ on a closed complex surface contains a
K\"{a}hler form if and only if it satisfies the following TWO conditions:
\begin{enumerate}
\item[(i)]
$\Omega \cdot \Omega [M] >0$
\item[(ii)]
$\Omega[C] >0$ for every complex curve $C$.
\end{enumerate}
\end{Lemma}

\begin{pf}
Nakai's criterion is well known for integral cohomology classes;
see, e.g., \cite[p.~127]{b-p-v}.
The criterion for rational classes follows immediately.
Now assume that $\Omega$ is a real cohomology class that
satisfies conditions (i) and (ii). Then these conditions
hold on some neighborhood of $\Omega$ in $H^2(M,\R)$.
The class $\Omega$ can be expressed as a convex combination
of rational cohomology classes which lie in this
neighborhood.  Every such class is represented by a K\"ahler form.
Since a convex combination of K\"ahler forms is again K\"ahler
(if the complex structure is fixed),
$\Omega$ can be represented by a K\"ahler form.
\end{pf}

To verify the first condition in Nakai's criterion, we will use moment
maps for closed two-forms which are not necessarily symplectic:

\begin{Lemma} \labell{exist-mom}
Let $M$ be complex surface obtained as above, and
let $\Omega$ be any cohomology class on $M$.
Then there exists a closed invariant one-form $\omega$
in the class $\Omega$, and for any such $\omega$ there exists
a moment map, i.e., a function $\Phi : M \to \R$, such that
\begin{equation} \labell{mom}
 d\Phi = -\iota(\xi_M) \omega
\end{equation}
where $\xi_M$ is the vector field that generates the action.
\end{Lemma}

\begin{pf}
We can take $\omega$ to be the average with respect to the circle 
action of any two-form that represents $\Omega$.
The one-form on the right hand side of \eqref{mom} is closed.
We need to show that it is exact. It is enough to show that 
its integral over any loop in $M$ is zero.
By induction on the number of blowups, one can show that 
every loop in $M$ is homologous to a loop in the minimal set
of the moment map.   But on this set $\xi_M=0$, so the one form
on the right hand side of \eqref{mom} is zero.
\end{pf}

Let us now check condition (i) of Nakai's criterion.

\begin{Lemma} \labell{nakai-1}
Let $M$ be a complex manifold with a holomorphic circle action,
obtained from $\CP^2$, a Hirzebruch surface, or a ruled manifold,
by a sequence of equivariant complex blowups at fixed points.
Let $\Omega$ be a cohomology class that is positive on invariant
complex curves.  Then $\Omega \cdot \Omega [M] >0$.
\end{Lemma}

\begin{pf}
Let $\omega$ be a closed invariant two-form representing $\Omega$,
and let $\Phi$ be a corresponding moment map.
As in \S\ref{subsec:DH}, let $a_\min$ and $a_\max$
denote $1/2\pi$ times the integrals of $\omega$ over the corresponding
fixed curves, and let $y_\min$, $y_\max$, and $y_p$, 
denote the values of $\Phi$ at the corresponding fixed points.
Since $\Omega$ is positive on fixed complex curves,
$a_\min$ and $a_\max$ are positive.
Since $\Omega$ is positive on invariant complex curves
that are not fixed, 
the value of $\Phi$ at the north pole of such a curve
is greater than the value of $\Phi$ at the south pole (as in \eqref{pos-C}).
Since every fixed point $p$ sits in a chain of such curves 
that starts at the minimal set of $\Phi$ and ends
at the maximal set of $\Phi$, we have $y_\min < y_p < y_\max$
for every interior fixed point $p$.
By the last part of Lemma \ref{rem:convex},
the push-forward of Liouville measure has a non-negative
density function which is nowhere zero. Hence the total measure, 
which is equal to ${1 \over 2} \Omega \cdot \Omega [M]$, is positive.
\end{pf}

It remains to check condition (ii) of Nakai's criterion:

\begin{Lemma} \labell{nakai-2}
Let $M$ be a complex manifold with a holomorphic circle action,
obtained from $\CP^2$, a Hirzebruch surface, or a ruled surface
by a sequence of equivariant complex blowups at fixed points. 
Let $\Omega$ be a cohomology class that is positive on invariant
complex curves. Then $\Omega$ is positive on any complex curve.
\end{Lemma}

It is enough to prove that every complex curve
is homologous to a positive combination of 
{\em invariant\/} complex curves:

\begin{Proposition} \labell{positive}
Let $M$ be a complex manifold with a circle action,
obtained from $\CP^2$, a Hirzebruch surface, or a ruled manifold
by a sequence of equivariant blow ups at fixed points. 
Then every complex curve in $M$ is homologous to a positive linear 
combination of invariant complex curves.
\end{Proposition}

The invariant complex curves generate the second homology group,
but they are not linearly independent; this makes the proof slightly
tricky.

\begin{pf*}{Proof of Proposition \ref{positive}
for a blow-up of a ruled surface}
Let $M$ be a complex surface with an $S^1$ action,
obtained from a ruled surface by a sequence of equivariant
complex blow-ups at fixed points.
Then $M$ contains two connected fixed complex curves, 
one consisting of minimal points and the other of maximal points,
and all the other invariant curves in $M$ are arranged 
in chains of gradient complex curves
(see Definition \ref{notation}).

Denote the fixed complex curves by $B_\min$ and $B_\max$.
Denote the spheres in a chain by $E_1, E_2, \ldots$, and denote
the orders of their stabilizers by $k_1, k_2, \ldots$.
{\em We suppress a second index that would denote which chain we're in.}
There are infinitely many gradient complex curves 
whose south pole lies on $B_\min$ and whose north pole lies on $B_\max$;
denote one of them by $F$. See Figure \ref{induction1}.
These invariant complex curves generate the homology of $M$:

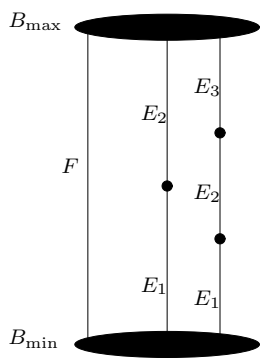
\begin{figure}
$$
\begin{picture}(284,149)(0,-10)
\put(60,125){\blacken\ellipse{70}{10}}
\put(60,5){\blacken\ellipse{70}{10}}
\put(60,65){\blacken\ellipse{4}{4}}
\put(80,45){\blacken\ellipse{4}{4}}
\put(80,85){\blacken\ellipse{4}{4}}
\path(30,125)(30,5)
\path(60,125)(60,5)
\path(80,125)(80,5)
\put(20,70){$\ss F$}
\put(50,90){$\ss E_2$}
\put(50,25){$\ss E_1$}
\put(70,100){$\ss E_3$}
\put(70,60){$\ss E_2$}
\put(70,20){$\ss E_1$}
\put(0,125){$\ss B_\max$}
\put(0,5){$\ss B_\min$}
\end{picture}
$$
\caption{Chains of gradient spheres in a space with two fixed surfaces}
\labell{induction1}
\end{figure}

\begin{Lemma} \labell{combination}
Let $C$ be a complex curve in $M$.  Then there exist 
integer coefficients $\alpha_\min$, $\alpha_\max$,
$\alpha_F$, $\alpha_i$, such that 
$$ [C] = \alpha_\min [B_\min] + \alpha_\max [B_\max] + \alpha_F [F]
   + \sum_{\text{chains}} \sum_i \alpha_i [E_i]$$
where $[ \cdot ]$ denotes homology class in $H_2(M,\Z)$.
Moreover, there exist such coefficients with $\alpha_1=0$ 
in each chain and with $\alpha_\min=0$. 
\end{Lemma}

\begin{pf}
If $M$ is itself a ruled surface, it is easy to check that 
its homology is generated by $B_\max$ and $F$. 
Otherwise, $M$ can be obtained from a ruled surface
by a sequence of blow ups; moreover, this can be done
in such a way that we never blow up a point at $B_\min$
(see Theorem \ref{thm:neat}).
Let $M$ be such a surface, and, arguing inductively, assume that 
$B_\max$, $F$, and the $E_i$'s, $i\geq 2$, generate
the homology group $H_2(M,\Z)$. 
Let $\tM$ be obtained from $M$ by an equivariant complex blow-up 
at a point $p$ that is not on $B_\min$.
A Meyer-Vietoris argument shows that the homology 
group $H_2(\tM,\Z)$ is generated by the proper transforms
$\tilde{B}_\max$, $\tF$, and $\tE_i$, and by the exceptional divisor, 
$\tE_p$. Since $p$ was not minimal, $\tE_p$ is not the first in its 
chain.  This proves the Lemma for $\tM$.
\end{pf}

It is enough to prove Proposition \ref{positive} for irreducible curves,
because every curve is the sum of irreducible curves. For invariant
curves, the proposition is tautological. Therefore, it is enough to
work with complex curves in $M$ that are Irreducible and are not
invariant.
Let $C$ be such a curve.  By Lemma \ref{combination}, we can write
\begin{equation} \labell{Csim}
 [C] = \alpha_\max [B_\max] + \alpha_F [F]
     + \sum_{\text{chains}} \sum_i \alpha_i [E_i]
\end{equation}
with $\alpha_1=0$ in each chain. We will show that all the coefficients
in this expression are non-negative.

The intersection number of any two irreducible complex curves
that are not equal to each other is non-negative.
Applying this to $C$ and $F$ (which are different because $C$ is not
invariant), we get
\begin{equation} \labell{CF}
C \cdot F \geq 0.
\end{equation}
Similarly,
\begin{equation} \labell{CBmin}
C \cdot B_\min \geq 0,
\end{equation}
and 
\begin{equation} \labell{CEi}
C \cdot E_i \geq 0 \quad \mbox{for all $i$, in each chain.} 
\end{equation}
We will express these intersection numbers
in terms of the coefficients of \eqref{Csim}.

The intersection number of any two different curves
among $B_\max$, $F$, and the $E_i$'s is $1$ if the curves intersect
and $0$ otherwise.
As for the self intersections, Lemma \ref{mg} implies that
the first and last spheres in a chain satisfy
$E_1 \cdot E_1 = -k_2$ and $E_l \cdot E_l = -k_{l-1}$, 
and the others satisfy
$E_i \cdot E_i = -\frac{k_{i-1} + k_{i+1}}{k_i}$,
where $k_i$ is the order of the stabilizer of the $i$th sphere
in a chain.

Substituting \eqref{Csim} in \eqref{CF} gives
$$ \alpha_\max \geq 0.$$
Substituting \eqref{Csim} in \eqref{CBmin}, and using
the assumption that $\alpha_1=0$ in each chain, gives
$$ \alpha_F \geq 0.$$
It remains to show that $\alpha_i \geq 0$ for all $i$, in each chain.
We will show 
\begin{equation} \labell{ineq}
 \frac{\alpha_{i+1}}{k_{i+1}} \geq \frac{\alpha_i}{k_i} \geq 0
\end{equation}
by induction on $i$.

Setting $i=1$ in \eqref{CEi} and substituting \eqref{Csim},
we get 
$C \cdot E_1 = \alpha_2 E_2 \cdot E_1 + \alpha_1 E_1 \cdot E_1$
$ = \alpha_2 - k_2 \alpha_1 \geq 0$. 
Since $k_1=1$ and $\alpha_1=0$, this gives
$\alpha_2/k_2 \geq \alpha_1/k_1 \geq 0$.
Taking \eqref{CEi} with $E_i$ that is neither the first nor last
in its chain, and substituting \eqref{Csim}, we get
$$ \begin{array}{l}
  \alpha_{i-1} (E_{i-1} \cdot E_i) + \alpha_i (E_i \cdot E_i)
   + \alpha_{i+1} (E_{i+1} \cdot E_i) = \\
  \alpha_{i-1} 
  - \frac{k_{i-1} + k_{i+1}} {k_i} \alpha_i 
  + \alpha_{i+1} \geq 0;
\end{array}$$
equivalently,
$$ 
	\alpha_{i+1} \geq \frac{k_{i+1}}{k_i} \alpha_i
		        + k_{i-1} \frac{\alpha_i}{k_i} - \alpha_{i-1}.
$$
Substituting the induction hypothesis
$\frac{\alpha_i}{k_i} \geq \frac{\alpha_{i-1}}{k_{i-1}} \geq 0$
and dividing by $k_{i+1}$, this gives
$\frac{\alpha_{i+1}}{k_{i+1}} \geq \frac{\alpha_i}{k_i} \geq 0$.
\end{pf*}

This completes the proof of Proposition \ref{positive}
when $M$ has two fixed complex curves.
To deal with one or no fixed curves,
we will use {\em backward\/} induction on blowups!
Here is the inductive step:

\begin{Lemma} \labell{op-induction}
Let $M$ be a complex surface with a holomorphic $S^1$ action,
and let $\pi : \tM \to M$ be an equivariant complex blowup at
a fixed point. Suppose that every complex curve in $\tM$
is homologous to a positive combination of invariant complex curves.
Then every complex curve in $M$ is homologous to a positive combination
of invariant complex curves. 
\end{Lemma}

\begin{pf}
Let $C$ be a complex curve in $M$. Its pre-image, $\pi\inv(C)$, is 
a complex curve in $\tM$ (more precisely, a union of complex curves).
By assumption, there exist invariant complex curves $C_j$ in $\tM$ 
such that $[\pi \inv(C)] = \sum \alpha_j [C_j]$ in $H_2(\tM)$
with $\alpha_j \geq 0$.
By applying $\pi$ to this equality, since homology pushes forward, 
and since each image $\pi(C_j)$ is either an invariant complex curve 
in $M$ or a single point, we get
$[C] = \sum \alpha_j [\pi(C_j)]$ in $H_2(M)$, which
expresses $C$ as a non-negative linear combination of invariant curves.
\end{pf}

To complete the proof of Proposition \ref{positive},
let $M$ be a complex surface with a holomorphic circle action,
obtained from either $\CP^2$ or a Hirzebruch by a sequence
of equivariant complex blowups at fixed points. 

\begin{Lemma} \labell{lem:blowup-max}
By a finite sequence of equivariant complex blow-ups,
we can get from $M$ to a complex surface $\tM$
that can be obtained from a ruled surface by a sequence
of equivariant complex blow-ups at fixed points.
\end{Lemma}

\begin{pf}
We first show how to replace an isolated minimum by a minimal $\CP^1$;
an isolated maximum can be treated similarly.
The isotropy weights $m,n$ at an isolated minimum $p$ are relatively 
prime positive integers. If $m>n$, blowing up replaces $p$ by
two new fixed points, one of which is interior, and the other
is minimal with isotropy weights $m-n$ and $n$. Repeating this,
we eventually arrive at a minimal points with weights $m=n=1$.
Blowing it up, it gets replaced by a minimal $\CP^1$.

In this way we obtain a surface $\tM$ that contains 
two connected fixed complex curves, one consisting 
of minimal points and the other of maximal points.
All the other invariant curves in $M$ are arranged 
in chains of gradient complex curves.
By the same reasoning as in the proof of 
Theorem \ref{thm:neat}, $M$ can be obtained from
a ruled surface by a sequence of equivariant complex blow-ups.
(This uses the fact that a complex curve of self intersection
$-1$ can always be blown down; see \cite[III,\S 4]{b-p-v}.)
\end{pf}

\begin{pf*}{Proof of the rest of Proposition \ref{positive}}
This follows immediately
from Lemmas \ref{op-induction} and \ref{lem:blowup-max}.
\end{pf*}

\begin{pf*} {Proof of Lemma \ref{nakai-2}}
Lemma \ref{nakai-2} follows immediately from proposition \ref{positive}.
\end{pf*} 

\begin{pf*} {Proof of Proposition \ref{kahler-cone}}
By Lemma \ref{nakai-1},
$\Omega^2$ is positive on the fundamental class $[M]$.
By Lemma \ref{nakai-2}, $\Omega$ is positive on complex curves.
By Nakai's criterion, Lemma \ref{nakai}, $\Omega$ contains a K\"ahler form.
\end{pf*}

\end{document}